\documentclass[10pt,a4paper]{article}
\usepackage[utf8]{inputenc}
\usepackage[T1]{fontenc}
\usepackage{amsmath, amssymb}
\usepackage{amsthm}
\usepackage{graphicx}
\usepackage[dvipsnames]{xcolor}
\usepackage{pgf,tikz,tikz-cd}
\usetikzlibrary{arrows,decorations.pathreplacing,calc}
\usepackage{array}
\usepackage{verbatim}
\usepackage{multicol}
\usepackage{mdframed}
\usepackage{pstricks}
\usepackage{pstricks-add}
\usepackage{tabularx}
\usepackage{fancybox}
\usepackage{stmaryrd}
\usepackage{comment}
\usepackage{caption}
\usepackage{mathrsfs}
\usepackage{enumitem}
\usepackage[all, cmtip]{xy}
\usepackage{pdfpages}
\usepackage{supertabular}
\usepackage{faktor}
\usepackage{todonotes}
\setuptodonotes{textcolor=blue, linecolor=blue, bordercolor=blue, backgroundcolor=white, size=footnotesize}
\usepackage{microtype}
\usepackage{hyperref}
\hypersetup{colorlinks=true, linkcolor=BlueViolet, filecolor=magenta, urlcolor=violet}
\usepackage{txfonts}
\usepackage[scaled=0.92]{helvet}
\usepackage[left=3.2cm,right=3.2cm,top=2.6cm,bottom=2.8cm]{geometry}
\usepackage{titlesec}
\titleformat{\section}{\large\bfseries}{\thesection}{1em}{}
\titleformat{\subsection}{\normalsize\bfseries}{\thesubsection}{1em}{}
\titleformat{\subsubsection}{\normalsize\itshape}{\thesubsubsection}{1em}{}
\titlespacing*{\section}{0pt}{3ex plus 1ex minus .2ex}{1.5ex plus .2ex}
\titlespacing*{\subsection}{0pt}{2.5ex plus 1ex minus .2ex}{1ex plus .2ex}
\usepackage{fancyhdr}
\fancypagestyle{plain}{\fancyhf{}}
\theoremstyle{plain}
\newtheorem{theo}{test}[section]
\newtheorem{theorem}[theo]{Theorem}
\newtheorem{prop}[theo]{Proposition}
\newtheorem{lemma}[theo]{Lemma}
\newtheorem{cor}[theo]{Corollary}

\newtheorem*{theoremA}{Theorem A}
\newtheorem*{theoremB}{Theorem B}
\newtheorem*{theoremC}{Theorem C}
\newtheorem*{theoremstar}{Theorem}

\theoremstyle{definition}
\newtheorem{definition}[theo]{Definition}

\theoremstyle{remark}
\newtheorem{rmq}[theo]{Remark}
\newtheorem{ex}[theo]{Example}

\DeclareMathOperator{\sing}{Sing}
\DeclareMathOperator{\Ker}{ker}

\DeclareMathOperator{\Id}{Id}

\DeclareMathOperator{\ord}{ord}
\DeclareMathOperator{\Idl}{Idl}
\DeclareMathOperator{\op}{op}

\DeclareMathOperator{\Hom}{Hom}
\DeclareMathOperator{\Proj}{Proj}
\DeclareMathOperator{\Projan}{Projan}
\DeclareMathOperator{\Int}{int}
\DeclareMathOperator{\Ass}{Ass}

\providecommand{\Qcal}{\mathcal{Q}}

\newcommand{\m}{\mathfrak{m}}
\newcommand{\n}{\mathfrak{n}}

\newcommand{\p}{\mathfrak{p}}

\newcommand{\Ocal}{\mathcal{O}}

\newcommand{\Mcal}{\mathcal{M}}

\newcommand{\Fcal}{\mathcal{F}}
\newcommand{\Ical}{\mathcal{I}}
\newcommand{\Rcal}{\mathcal{R}}

\newcommand{\inj}{\hookrightarrow}

\newcommand{\Xbar}{\overline{X}}

\newcommand{\N}{\mathbb{N}}
\newcommand{\Z}{\mathbb{Z}}
\newcommand{\Q}{\mathbb{Q}}
\newcommand{\R}{\mathbb{R}}
\newcommand{\C}{\mathbb{C}}

\author{Fran\c{c}ois Bernard \and Antoni Rangachev}
\date{}
\title{Hölder continuity of meromorphic functions and equisingularity}
\begin{document}
\maketitle
\thispagestyle{plain}
\begin{abstract}
For a complex analytic variety and an exponent between zero and one, we study the sheaf of meromorphic functions that are locally Hölder continuous with that exponent for the outer metric. It sits between the structure sheaves of the seminormalization and of the Lipschitz saturation of Pham and Teissier. We prove that these sheaves are coherent and that only finitely many of them are distinct, the exponents at which they change being rational numbers governed by the Rees valuations of the ideal of a certain diagonal. Each exponent defines a Hölder saturation of the variety. We show that a variety whose Hölder saturation is smooth for some exponent larger than one third is itself smooth and that one third is sharp. For reduced curve germs we compute the Hölder saturations from the characteristic exponents of the branches and give explicit thresholds, in terms of classical invariants, at which they become the seminormalization, the Lipschitz saturation, or a curve of the same multiplicity. We prove that homeomorphic curves are bi-Hölder homeomorphic for some positive exponent, whereas bi-Hölder homeomorphisms with exponent close enough to one preserve the multiplicity. We then introduce the relative Hölder saturation of a family and prove that it specializes to the saturations of the fibers over a dense Zariski-open subset of the base. This leads to two notions of equisingularity for families of complex analytic varieties, Hölder equisingularity, defined by a bi-Hölder trivialization, and the stronger Hölder equisaturation, defined by the triviality of the relative saturation. For families of reduced curve germs we show that they interpolate between topological, Whitney and Lipschitz equisingularity, with explicit exponents depending only on the invariants of the special fiber.
\end{abstract}

\tableofcontents

\section*{Introduction}
{\let\thefootnote\relax\footnotetext{\textit{2020 Mathematics Subject Classification.} Primary 32S15, 32S05; Secondary 13B22, 14H20, 14B05.\\ \textit{Key words and phrases.} Hölder saturation, Lipschitz saturation, seminormalization, rational powers of ideals, Rees valuations, Hölder equisingularity, Whitney equisingularity, Milnor number, multiplicity, curve singularities.}}


Let $X$ be a reduced  complex analytic variety. Its normalization $\overline X\to X$ is a finite bimeromorphic map, and every finite bimeromorphic map onto $X$ factors through it. Between $X$ and $\overline X$ sit finite bimeromorphic modifications of $X$ obtained by adjoining to $\Ocal_X$ locally bounded meromorphic functions. Two of these modifications are classical: Adjoining the continuous meromorphic functions gives the seminormalization $X^+$ of Andreotti and Norguet, and adjoining the Lipschitz ones gives the Lipschitz saturation $X^s$ of Pham and Teissier. Both $X^+\to X$ and $X^s\to X$ are homeomorphisms, and both are isomorphisms when $X$ is normal. In general they are not: for the cusp $y^2=x^3$ both $X^+$ and $X^s$ are the normalization, which is smooth. The seminormalization depends only on how the points of $\overline X$ are glued over the points of $X$, so it is a topological invariant of $X$. The Lipschitz saturation records which meromorphic functions are Lipschitz for the outer metric, and for plane curves this already determines the equisingularity type by Zariski and Pham--Teissier. The two differ exactly when $X$ carries continuous meromorphic functions which are not Lipschitz. This paper fills the gap with the $\alpha$-Hölder saturations obtained by adjoining the meromorphic functions that are Hölder continuous with exponent $\alpha$.  Furthermore, these constructions turn out to have a natural application to the equisingularity theory  of families of curves.

\textbf{Hölder meromorphic functions.} Let $\Mcal_X$ be the sheaf of meromorphic functions on $X$. By the Riemann extension theorem, the locally bounded meromorphic functions form the structure sheaf of $\overline X$, and by Andreotti and Norguet \cite{AndreottiNorguet1967SeminormalizationAnalytic} the continuous ones form the structure sheaf $\Mcal^c_X$ of $X^+$. For a real number $\alpha\in(0,1]$ denote by $\Mcal^\alpha_X$ the subsheaf of $\Mcal^c_X$ of meromorphic functions which are locally $\alpha$-Hölder for the outer metric, that is, the metric induced by a local embedding of $X$ in an affine space. For $\alpha=1$ its analytic spectrum is the Lipschitz saturation $X^s$ \cite{PhamTeissier1969LipFractions}, which coincides with Zariski's saturation \cite{Zariski1968StudiesIII} for hypersurfaces \cite{PhamTeissier1969LipFractions,Lipman1975RelativeLipSat}. For the cusp $X=V(y^2-x^3)\subset(\C^2,0)$, the function $y/x$ becomes $t$ on the normalization $t\mapsto(t^2,t^3)$, so it is continuous on $X$ but not Lipschitz, and $y/x\in\Mcal^\alpha_{X,0}$ if and only if $\alpha\leq1/3$ (Example~\ref{Ex_Cusp}). The bound comes from the points with parameters $t$ and $-t$, which are at distance $2|t|^3$, while the values of $y/x$ at them differ by $2|t|$.

Our first results say that the sheaves $\Mcal^\alpha_X$ are coherent and that only finitely many of them are distinct.

\begin{theoremA}\label{Th_IntroA}
Suppose $X$ is a sufficiently small representative of a complex analytic germ, a compact complex analytic variety, or a complex algebraic variety. Then $\Mcal^\alpha_X$ is a coherent $\Ocal_X$-algebra for every $\alpha\in(0,1]$, and there exist finitely many rational numbers $0<\alpha_1<\dots<\alpha_r=1$ such that
\[
\Mcal^\alpha_X\simeq\Mcal^{\alpha_{i+1}}_X\quad\text{for }\alpha\in(\alpha_i,\alpha_{i+1}]
\qquad\text{and}\qquad
\Mcal^\alpha_X\simeq\Mcal^c_X\quad\text{for }\alpha\in(0,\alpha_1].
\]
\end{theoremA}

Theorem~A is Theorem~\ref{Theo_FiniteNumberOfAlphas} for compact and algebraic varieties and Proposition~\ref{Prop_FiniteNumberOfAlphasCoherentLOCAL} for germs; the coherence is Lemma~\ref{Lem_Coherence}. Denote by $X^\alpha$ the analytic spectrum of $\Mcal^\alpha_X$ and call it the $\alpha$-Hölder saturation of $X$. The natural map $\pi_\alpha : X^\alpha\to X$ is finite, bimeromorphic and bijective, and it is locally bi-$\alpha$-Hölder (Corollary~\ref{Cor_HolderSatMapIsHolder}). Thus $X$ comes with a finite chain of bimeromorphic bi-Hölder modifications
\[
\overline X\longrightarrow X^{\alpha_1}=X^+\longrightarrow X^{\alpha_2}\longrightarrow\dots\longrightarrow X^{\alpha_r}=X^s\longrightarrow X,
\]
in which the first map is the normalization of the seminormalization. The proof of Theorem~A is algebraic. Let $\pi : \overline X\to X$ be the normalization and let $\Ical$ be the ideal of the diagonal, namely the kernel of $\Ocal_{\overline X\times\overline X}\to\Ocal_{\overline X\times_X\overline X}$. In a local embedding $X\subset\C^N$ with coordinates $z_1,\dots,z_N$, the ideal $\Ical$ is generated by the differences $z_j\otimes1-1\otimes z_j$. For a rational number $\alpha=p/q$, let $\Ical_\alpha$ be the rational power of $\Ical$, that is, the ideal of the elements $f$ with $f^q\in\overline{\Ical^p}$, where the bar denotes integral closure (Section~\ref{Sec_HolderSaturations}). A theorem of Lejeune-Jalabert and Teissier \cite[Theorem~7.2]{LejeuneTeissier2008ClotureInt} says that $f\in\Ical_\alpha$ if and only if $|f|\leq C\sup_j|z_j\otimes1-1\otimes z_j|^\alpha$ near every point (Lemma~\ref{Lem_InequalityPowerIdeal}). Hence a function $f$ on $\overline X$ is $\alpha$-Hölder on $X$ if and only if $f\otimes1-1\otimes f\in\Ical_\alpha$. Because rational powers of coherent ideals are coherent \cite[Proposition~4.11]{LejeuneTeissier2008ClotureInt}, the sheaf $\Mcal^\alpha_X$, which is the kernel of $f\mapsto f\otimes1-1\otimes f$ modulo $\Ical_\alpha$, is coherent. The finiteness comes from the Rees valuations $v_1,\dots,v_m$ of $\Ical$. If $e$ is the least common multiple of the integers $v_i(\Ical)$, then $\Ical_\alpha=\Ical_{n/e}$ for every $\alpha\in((n-1)/e,n/e]$ (Proposition~\ref{Prop_FinitenessOfRationalIdeals}). So the exponents $\alpha_i$ of Theorem~A are of the form $n/e$. For a compact or algebraic variety one uses the finiteness of the number of irreducible components to find a common denominator (Proposition~\ref{Prop_FiniteNumberOfAlphasCoherentGLOBAL}). Section~\ref{SectionSaturations} develops this construction in a general form. To an ideal operation $\phi$ on $\Ocal_{\overline X\times\overline X}$ we attach the $\phi$-saturation of $X$, which recovers the seminormalization for $\phi=\sqrt{\cdot}$, the Lipschitz saturation for the integral closure, the Arf and strict saturations of Lipman \cite{Lipman1971ArfRings}, and the Hölder saturations for the rational powers. This generality is forced on us: the $\alpha$-Hölder saturation is not an idempotent operation. The $\alpha$-Hölder saturation of $X^\alpha$ can be strictly larger than $X^\alpha$, and in general one only has $(X^\alpha)^\beta\subseteq X^{\alpha\beta}$ (Proposition~\ref{Prop_InteractionsBetweenSaturations}).

\textbf{The $1/3$ theorem.} By a theorem of Birbrair, Fernandes, L\^e and Sampaio \cite{BirbrairFernandesLeSampaio2016}, a complex analytic set which is bi-Lipschitz homeomorphic to a smooth one is smooth. This generalizes Mumford's theorem that a normal surface which is a topological manifold is smooth \cite{Mumford1961}. Sampaio \cite{Sampaio2025LogLipschitz} shows that a complex analytic set which is bi-$\alpha$-Hölder homeomorphic to a smooth one for every $\alpha<1$ is smooth, and that no single exponent depending only on the dimension suffices. The Hölder saturations give a version of this rigidity for bimeromorphic maps with an explicit exponent.

\begin{theoremstar}[Theorem~\ref{Th_HolderSatSmoothness}]
Let $(X,0)$ be a reduced complex analytic germ and let $\alpha\in(1/3,1]$. If $(X^\alpha,0)$ is smooth, then $(X,0)$ is smooth.
\end{theoremstar}

In particular, if $h : (X,0)\to(X',0)$ is a bimeromorphic map onto a smooth germ which is a bi-$\alpha$-Hölder homeomorphism for some $\alpha>1/3$, then $(X,0)$ is smooth. Indeed, the coordinates of $X$ composed with $h^{-1}$ are $\alpha$-Hölder meromorphic functions on the smooth germ $X'$, hence holomorphic, so $h^{-1}$ is the normalization of $X$; and a holomorphic function on $X'$ composed with $h$ is $\alpha$-Hölder on $X$, so $X^\alpha=\overline X$ is smooth. The exponent $1/3$ is sharp, because the normalization of the cusp is bi-$\frac13$-Hölder.

\textbf{Curves.} Section~\ref{Sec_ComputingHolderSat} computes the Hölder saturations of reduced curve germs. Let $(X,0)$ be a reduced curve germ. Denote by $m=m(X,0)$ its multiplicity, by $r=r(X,0)$ its number of branches, by $\delta=\delta(X,0)$ its delta invariant, that is, the length of $\Ocal_{\overline X}/\Ocal_X$, and by $\mu=\mu(X,0)$ its Milnor number in the sense of Buchweitz and Greuel \cite{BuchweitzGreuel1980ConstancyOfMilnorNbr}, so that $\mu=2\delta-r+1$. If $X$ is a branch with characteristic exponents $\beta_0=m<\beta_1<\dots<\beta_g$ and $e_i=\gcd(\beta_0,\dots,\beta_i)$, then $X^\alpha$ is the monomial curve whose exponents are the integers $p$ for which there is an $i$ with $p\geq\alpha\beta_i$ and $e_i\mid p$ (Proposition~\ref{Prop_ComputeSemigroupOfTheHolderSat}); the saturations of a reducible curve are described in Theorem~\ref{Theo_HolderSatReducibleCurve} by the same rule on each branch together with a condition on the contact of pairs of branches. From these descriptions we derive bounds that involve only the numerical invariants. For a reduced curve germ the following holds (Corollary~\ref{Cor_HolderSatSeminormBoundReducedCurves}, Corollary~\ref{Cor_BoundAlphaIsLipReduciblecase} and Proposition~\ref{Prop_HolderSatPreservesMultiplicityReduced}):
\begin{enumerate}[label=\textup{(\roman*)},leftmargin=*]
\item if $\alpha\leq\dfrac{1}{\mu-m+3}$, then $X^\alpha=X^+$ is the seminormalization;
\item if $\alpha>\dfrac{\mu-m+2}{\mu-m+3}$, then $X^\alpha=X^s$ is the Lipschitz saturation;
\item if $\alpha>\dfrac{m-r}{m-r+2}$, then $(X^\alpha,0)$ and $(X,0)$ have the same multiplicity.
\end{enumerate}
For the cusp, (i) and (iii) give the exponent $1/3$, which is sharp by Example~\ref{Ex_Cusp}.

\textbf{Hölder homeomorphisms.} Homeomorphisms between complex analytic sets are of Hölder nature. We show that a homeomorphism between curves can always be replaced by a bi-$\alpha$-Hölder homeomorphism for some $\alpha>0$, and that bi-$\alpha$-Hölder homeomorphisms with $\alpha$ sufficiently large preserve algebraic invariants such as the multiplicity. The following result illustrates this phenomenon for curves; its generalization to higher dimensions will appear elsewhere.

\begin{theoremB}\label{Th_IntroB}
\begin{enumerate}[label=\textup{(\roman*)},leftmargin=*]
\item \emph{Local case.} Suppose $(C_1,0)$ and $(C_2,0)$ are reduced complex analytic curve germs in $(\C^N,0)$ with the same number of branches. Then there exist $\alpha\in(0,1]\cap\Q$ and a bi-$\alpha$-Hölder homeomorphism $(C_1,0)\to(C_2,0)$.
\item \emph{Global case.} Suppose $C_1$ and $C_2$ are reduced projective curves in $\mathbb P^N_\C$ and $h : C_1\to C_2$ is a homeomorphism for the Euclidean topology. Then there exist $\alpha\in(0,1]\cap\Q$ and a bi-$\alpha$-Hölder homeomorphism $\tilde h : C_1\to C_2$, where the Hölder condition refers to the metrics induced by the Fubini--Study metric on $\mathbb P^N_\C$.
\item Suppose $F : (C_1,0)\to(C_2,0)$ is a bi-$\alpha$-Hölder homeomorphism of reduced curve germs. Let $M$ be the maximum of the multiplicities $m(C_1,0)$ and $m(C_2,0)$. If $\alpha>\sqrt{M/(M+1)}$, then $m(C_1,0)=m(C_2,0)$.
\end{enumerate}
\end{theoremB}

Parts (i) and (ii) are Theorem~\ref{Th_HomeoImpliesHolderHomeo} and part (iii) is Corollary~\ref{Cor_MetricMultiplicityContinuation}. The proof of (i) and (ii) rests on Theorem~A: the seminormalization $C_i^+\to C_i$ is a homeomorphism, and it is bi-$\varepsilon_i$-Hölder for the smallest threshold $\varepsilon_i=\alpha_1$ of Theorem~A; the seminormalizations are bi-Lipschitz homeomorphic to each other because their singularities are ordinary multiple points, which are determined by their number of branches. In particular, for curves the existence of a bi-$\alpha$-Hölder homeomorphism for some unspecified $\alpha>0$ is a purely topological condition. Part (iii) rests on an elementary cluster argument. The $m$ sheets of a branch of multiplicity $m$ over a transverse parameter come together in clusters whose sizes and mutual distances are read off from the characteristic sequence, and two sheets cannot be joined by a short path on the curve. A Hölder homeomorphism cannot tear a cluster apart, so the image branch must have enough sheets to accommodate it. The square root comes from the use of the Hölder estimate for both $F$ and $F^{-1}$. Part (iii) is to be compared with a result of Fernandes, Sampaio and Silva \cite{FernandesSampaioSilva2018}, who show that for plane branches a bi-$\alpha$-Hölder homeomorphism with $\alpha$ close to $1$ preserves the characteristic exponents.

\textbf{Equisingularity.} The second part of the paper is devoted to applications of the Hölder saturations to equisingularity. We introduce the notion of $\alpha$-Hölder equisaturation for families of complex analytic varieties as a unifying notion between Lipschitz and Whitney equisingularity. In this paper we treat families of curves; the general case will appear elsewhere.

\textit{Family setup.} Let $p : (X,0)\to(Y,0)$ be a flat morphism of complex analytic germs, where $Y$ is smooth of dimension $k$, $X$ is reduced and equidimensional, and the fibers $X_y$ are reduced curves. We say that $p$ is centered if there exists a section $\sigma : (Y,0)\to(X,0)$ with $\sing(p)=\sigma(Y)$, where $\sing(p)$ is the relative singular locus of $p$. In this case we identify $Y$ with $\sigma(Y)$ and choose an embedding $(X,0)\subset(\C^n\times\C^k,0)$ with $\sigma(y)=(0,y)$. Set $m_0=m(X_0,0)$, $r_0=r(X_0,0)$, $\delta_0=\delta(X_0,0)$ and $\mu_0=\mu(X_0,0)$.

\textit{Whitney equisingularity.} Let $\{x_i\}$ be a sequence of smooth points of $X$ and $\{y_i\}$ a sequence of points of $Y$, both converging to $0$. Consider the secant lines $\overline{x_iy_i}$, the tangent spaces $T_{x_i}X$ and their limits $l$ and $T$, which exist after passing to subsequences. We say that $p$ is Whitney equisingular if $l\subset T$ for every such pair of sequences. Whitney equisingularity is equivalent to the constancy of $\mu(X_y,0)$ and $m(X_y,0)$ along $Y$, and also to the existence of a strong simultaneous resolution of $p$. For $k=1$ this is due to Briançon, Galligo and Granger \cite{BrianconGallicoGranger1980DeformationsEqu}, who use results of Buchweitz and Greuel \cite{BuchweitzGreuel1980ConstancyOfMilnorNbr}; for arbitrary $k$ it is due to Bengu\c{s}-Lasnier, Gaffney and the second author \cite{BengusGaffneyRangachev2026}.

\textit{Topological equisingularity.} We say that $p$ is topologically equisingular if there exists a homeomorphism $\Phi : (X,0)\to(X_0,0)\times(Y,0)$ with $p=\mathrm{pr}_2\circ\Phi$. For centered families, topological equisingularity is equivalent to the constancy of the Milnor number $\mu(X_y,0)$: for $k=1$ this is due to Buchweitz and Greuel \cite{BuchweitzGreuel1980ConstancyOfMilnorNbr}, and for arbitrary $k$ to \cite{BengusGaffneyRangachev2026}. Since Whitney equisingularity implies the constancy of $\mu$, it implies topological equisingularity.

\textit{Hölder equisingularity.} We say that $p$ is $\alpha$-Hölder equisingular if there exists a homeomorphism $\Phi : (X,0)\to(X_0,0)\times(Y,0)$ with $p=\mathrm{pr}_2\circ\Phi$ such that $\Phi$ and $\Phi^{-1}$ are $\alpha$-Hölder for the outer metrics. For $\alpha=1$ this is the Lipschitz equisingularity studied by Mostowski \cite{Mostowski1985}, Parusi\'nski \cite{Parusinski1994}, Parusi\'nski and P\u{a}unescu \cite{ParusinskiPaunescu2017}, and Neumann and Pichon \cite{NeumannPichon2014,NeumannPichon2022}. Hölder equisingularity is a metric condition on the total space. To tie it to the algebra of Hölder saturations we use a relative version of them.

\textit{Hölder equisaturation.} Let $\pi : \overline X\to X$ be the normalization and let $\Ical$ be the ideal of the diagonal in $\Ocal_{\overline X\times_Y\overline X}$. The relative $\alpha$-Hölder saturation $X^\alpha_Y$ is the finite space over $X$ whose structural algebra consists of the functions $h$ on $\overline X$ with $h\otimes1-1\otimes h\in\Ical_\alpha$. A function $h$ belongs to it if and only if it is $\alpha$-Hölder on every fiber with a Hölder constant independent of the fiber (Proposition~\ref{Prop_RelativeHolderEstimate}). We say that $p$ is $\alpha$-Hölder equisaturated if $X^\alpha_Y$ is a product over $Y$, that is, $X^\alpha_Y\simeq F\times Y$ over $Y$ for some curve germ $F$ (Definition~\ref{Def_AlphaEquisaturation}). For $\alpha=1$ this is the Lipschitz equisaturation of \cite{Lipman1975RelativeLipSat,Bernard2025RelativeLipSat}. An $\alpha$-Hölder equisaturated family is $\alpha$-Hölder equisingular: the composition of $F\times Y\simeq X^\alpha_Y\to X$ with the inverse of $F\to X_0$ is a bi-$\alpha$-Hölder trivialization (Proposition~\ref{Prop_ProductHolderTriviality}). Section~\ref{Sec_RelativeSpecialization}, which is the most technical part of the paper, shows that the relative Hölder saturation specializes to the Hölder saturations of the fibers over a dense Zariski-open subset of $Y$: for $y$ in such a subset, $(X^\alpha_Y)_y\simeq(X_y)^\alpha$, and for finitely many exponents the subset can be chosen simultaneously (Theorem~\ref{Th_GenericHolderSpecialization}). The proof reduces to showing that the normalization of the graded algebra $\Ocal_W[\Ical t,\Ical t^b]$, whose degree $a$ part is the rational power $\Ical_{a/b}$ of the diagonal ideal on $W=\overline X\times_Y\overline X$, commutes with restriction to a general fiber. Away from the finitely many pairs of points of $\overline X$ lying over a singular point of a fiber this is an explicit computation; at those pairs it follows from the generic normality of the fibers of a normal space over a smooth base. Thus for an $\alpha$-Hölder equisaturated family the factor $F$ is the $\alpha$-Hölder saturation of the generic fiber, and it sits between $X_0$ and $X_0^\alpha$. The next result summarizes the main results of Section~\ref{Sec_HolderEquisaturation}.

\begin{theoremC}\label{Th_IntroC}
Let $p : (X,0)\to(Y,0)$ be a centered flat family of reduced curve germs over a smooth base. The following holds.
\begin{enumerate}[label=\textup{(\roman*)},leftmargin=*]
\item If $p$ is $\alpha$-Hölder equisaturated for some $\alpha>\dfrac{\mu_0-m_0+2}{\mu_0-m_0+3}$, then $p$ is Lipschitz equisaturated. Conversely, if $p$ is Lipschitz equisaturated, then it is $\alpha$-Hölder equisaturated for every $\alpha\in(0,1]$. If $p$ is $\alpha$-Hölder equisaturated for some $\alpha>\dfrac{m_0-r_0}{m_0-r_0+2}$, then $p$ is Whitney equisingular. If $p$ is $\alpha$-Hölder equisaturated for some $\alpha\in(0,1]$, then $p$ is $\alpha$-Hölder equisingular, hence topologically equisingular, and $\mu(X_y,0)$ is constant.
\item If $p$ is $\alpha$-Hölder equisingular for some $\alpha>\sqrt{m_0/(m_0+1)}$, then $p$ is Whitney equisingular.
\item If $p$ is topologically equisingular, then $p$ is $\alpha$-Hölder equisaturated for every $\alpha\leq\dfrac{1}{2\delta_0+1}$.
\end{enumerate}
\end{theoremC}

The first assertion of (i) is Proposition~\ref{Prop_HolderEquisaturationImpliesLipschitz} together with the Lipschitz bound for curves recalled above; the second is Proposition~\ref{Prop_InteractionsBetweenSaturations} applied to the relative saturation; the third is Theorem~\ref{Th_EquisaturationWhitney}, where a sharper bound in terms of the characteristic exponents of the branches of $X_0$ is given; the fourth is Proposition~\ref{Prop_ProductHolderTriviality} together with Corollary~\ref{Cor_MuSomeAlpha}. Part (ii) follows from Theorem~B\,(iii) applied to the restrictions of a trivialization to the fibers, which gives the constancy of the multiplicity of $X_y$ at $0$, together with the constancy of the Milnor number along topologically equisingular families. Part (iii) is Theorem~\ref{Th_UniformSmallExponent}; it refines the topological triviality theorem of \cite{BengusGaffneyRangachev2026} by an explicit Hölder exponent, and Example~\ref{Ex_RelativeSeminormalExponent} shows that the exponent $1/(2\delta_0+1)$ can be attained. Note that (iii) and the last assertion of (i) say that topological equisingularity is the same as $\alpha$-Hölder equisaturation for some $\alpha$, and that Whitney equisingularity lies strictly between: Remark~\ref{Prop_NoUniversalWhitneyExponent} gives, for every $\alpha\in(0,1]$, a Whitney equisingular family of branches which is not $\alpha$-Hölder equisaturated.

\textbf{Organization.} Section~\ref{SectionSaturations} develops the general theory of $\phi$-saturations and proves the coherence of the saturation of a coherent sheaf. Section~\ref{Sec_HolderSaturations} introduces the rational powers of ideals and the Hölder saturations, proves Theorem~A, identifies the Hölder saturations with the sheaves of Hölder meromorphic functions, and proves the $1/3$ theorem. Section~\ref{Sec_ComputingHolderSat} computes the Hölder saturations of curve germs from the characteristic sequences of their branches and derives the numerical bounds above. Section~\ref{Sec_HolderHomeoCurves} proves Theorem~B. Section~\ref{Sec_RelativeSpecialization} introduces the relative Hölder saturation and proves that it specializes generically to the Hölder saturations of the fibers. Section~\ref{Sec_HolderEquisaturation} introduces $\alpha$-Hölder equisaturation, compares it with Lipschitz, Whitney and topological equisingularity, and proves Theorem~C.

\textbf{Acknowledgments}
	The authors would like to thank  Bernard Teissier for many stimulating and helpful discussions.  The authors were supported by the Swiss National Science Foundation grant no.~228841, ``Topological, Metric, and Algebraic Equisingularity Theory'', under the program PROMYS and by the Ministry of Education and Science of Bulgaria.

\textbf{AI Disclosure}
The second author gave a short course on $\alpha$-Hölder saturations during the conference “Metric Geometry of Singularities,” held at the University of Chicago Center in Paris in June 2021. Most of the research presented in this manuscript was carried out by the authors at IMJ-PRG during 2024–2025, while the final part (Section~\ref{Sec_HolderEquisaturation}) was completed at IMI-BAS in 2026.

AI tools (Claude Fable 5.1 and ChatGPT Astra 6.1) were used for proofreading and improving the exposition of the manuscript. They also contributed to the simplification and development of several arguments. In particular, AI found a substantially simpler analytic proof of Theorem~\ref{Th_HolderSatSmoothness}, replacing an algebraic argument involving localizations at height-one prime ideals. It suggested the cluster argument used in Proposition~\ref{Prop_ClusterLemma} as an alternative to the original approach based on the work of~\cite{FernandesSampaioSilva2018}, and identified a simplification in the proof of Lemma~\ref{Lem_RationalPowersGeneralFibres}. AI also assisted in constructing some of the examples in Section~\ref{Sec_HolderEquisaturation}. All mathematical arguments have been carefully checked by the authors, who take full responsibility for the content of the final manuscript.
\section{General properties of saturations}\label{SectionSaturations}

In this section we set up a general theory of saturations. Let $R\to A\to B$ be ring morphisms and let $\phi$ be an operation on the ideals of $B\otimes_RB$, such as taking the radical, the integral closure, or a rational power. The $\phi$-saturation of $A$ in $B$ over $R$ is the set of elements $b\in B$ such that $b\otimes1-1\otimes b$ lies in the image under $\phi$ of the ideal of the diagonal. For $R=\C$, $A=\Ocal_{X,0}$ and $B=\Ocal_{\overline X,0}$ this recovers the seminormalization for the radical, the Lipschitz saturation of Pham and Teissier \cite{PhamTeissier1969LipFractions} for the integral closure, the Arf and strict saturations of Lipman \cite{Lipman1971ArfRings}, and, for the rational powers, the Hölder saturations of this paper. Two features of the Hölder saturations force us to work in this generality. There is one saturation for every exponent $\alpha\in(0,1]$, so we need to know how saturations for different operations are related, and the Hölder saturation is not idempotent, so we need to know what happens when saturations are composed (Proposition~\ref{Prop_InteractionsBetweenSaturations}). We first treat rings, then sheaves on complex analytic varieties, where the main point is that the saturation of a coherent sheaf for a coherent ideal operation is coherent (Lemma~\ref{Lem_Coherence}). We believe that this framework can be of use beyond the present paper, for instance for the Lipschitz saturation of ideals of Gaffney, or for the seminormalization of ideals of Gaffney and Vitulli.\\

Let \textbf{Poset} denote the category of partially ordered sets whose morphisms are monotone functions. We denote by $\Idl$ the functor
\[
\begin{array}{rccc}
     \Idl : & \textbf{Ring}^{\op} & \to & \textbf{Poset}   \\
     & R & \mapsto & \{ \text{ideals of }R\}
\end{array}
\]
where the partial order is given by inclusion of ideals.

\begin{definition}
    Let $R$ be a ring. An element $\phi$ of $ \Hom_{\textbf{Poset}}\big(\Idl(R), \Idl(R)\big) $ will be called an "ideal operation" on $R$. For every ideal $I$, we will write $I^\phi := \phi(I)$.
\end{definition}

\begin{rmq}
    Since $\phi$ is a morphism in \textbf{Poset}, it preserves the order relation, that is, for all ideals $I\subseteq J$, we have $I^\phi\subseteq J^\phi$.
\end{rmq}

\begin{definition}
    Let $R\to A\to B$ be ring morphisms and let $\phi$ be an ideal operation on $B\otimes_R B$. Then the $\phi$-saturation of $A$ in $B$ relatively to $R$ is defined by
    \[
    A^\phi_{B,R} := \Big\{ b\in B \mid b\otimes_R 1 - 1\otimes_R b \in \Ker(\psi)^\phi\Big\}
    \]
    where $\psi : B\otimes_R B \to B\otimes_A B$ is the canonical map.
\end{definition}

\begin{rmq}\label{Rmq_ReecritureDefSaturation}
    If we consider the morphism
    $$\begin{array}{rccc}
        \Delta : & B & \to & B\otimes_R B \\
        & b & \mapsto & b\otimes 1 - 1\otimes b
    \end{array}
    $$
    then we have $\ker(\psi) = \langle \Delta(g(a)) \mid a\in A\rangle$ where $g$ is the morphism $A\to B$. Hence we get
    $$
    A^\phi_{B,R} = \Delta^{-1}\Big(\big\langle \Delta\circ g(a) \mid a\in A\big\rangle^\phi\Big) = \Delta^{-1}\Big(\big\langle \Delta\circ g(A)\big\rangle^\phi\Big).
    $$
\end{rmq}

\begin{lemma}\label{Lem_SatProcessIsOrderPreserving}
    Let $R\to A\to B$ be ring morphisms and let $\phi_1$ and $\phi_2$ be ideal operations. Then, if $\phi_1 \leq \phi_2$, i.e. $I^{\phi_1}\subseteq I^{\phi_2}$ for all ideal $I$, then 
    $$ A^{\phi_1}_{B,R} \subseteq A^{\phi_2}_{B,R}.$$
\end{lemma}

\begin{proof}
    If $\phi_1 \leq \phi_2$, then $\ker(\psi)^{\phi_1} \subseteq \ker(\psi)^{\phi_2}$. So 
    $$A^{\phi_1}_B = \Delta^{-1}\big(\ker(\psi)^{\phi_1}\big) \subseteq \Delta^{-1}\big(\ker(\psi)^{\phi_2}\big) = A^{\phi_2}_B.$$
\end{proof}

If $A\to B$ is an integral morphism, then we recover the classical notions of seminormalization and Lipschitz saturation.

\begin{definition}
    Let $A\to B$ be an integral morphism of R-algebras. Then $A^+_B := A^{\sqrt{.}}_B$ is called the seminormalization of $A$ in $B$ and $A^s_B := A^{\overline{(.)}}_B$ is called the Lipschitz saturation of $A$ in $B$.
\end{definition}

\begin{rmq}\label{Rmq_IdpluspetitquePhi}
    A trivial example of an ideal operation is the constant map $\phi : I \mapsto B\otimes_R B$, for which $A^\phi_{B,R} = B$. Another natural example is $\phi = \Id$, in which case $A^{\Id}_{B,R}$ is the strict saturation of $A$ in $B$, a notion that appears, for instance, in Lipman's work \cite{Lipman1971ArfRings}. One can observe that the strict saturation always satisfies $g(A)\subseteq A^{\Id}_{B,R}$. More generally, for any ideal operation $\phi$, we have $g(A) \subseteq A^\phi_{B,R}$ if and only if $A^{\Id}_{B,R} \subseteq A^\phi_{B,R}$. In particular, for a general ideal operation, the inclusion $g(A) \subseteq A^\phi_{B,R}$ does not necessarily holds, even though $A^{\phi}_{B,R}$ is by definition the $\phi$-saturation of
    $$ R\to g(A) \inj B.$$
    Take for example $R=\C$ and $A=B=\C[x]$ with the ideal operation $\phi : I\mapsto \{0\}$, for which $A^\phi_{B,R} = R$. Henceforth, we will assume that all ideal operations satisfy the condition $\Id\leq \phi$. This ensures that we can always consider the following commutative diagram:
    \[\begin{tikzcd}
    	R & A & {g(A)} & B \\
    	&& {A^\phi_{B,R}}
    	\arrow[from=1-1, to=1-2]
    	\arrow[two heads, from=1-2, to=1-3]
    	\arrow[hook, from=1-3, to=1-4]
    	\arrow[hook, from=1-3, to=2-3]
    	\arrow[hook, from=2-3, to=1-4]
    \end{tikzcd}\]
    In particular, we obtain a diagram $R\to A^\phi_{B,R}\to B$, which allows us to iterate the saturation process.
    \end{rmq}

\begin{lemma}\label{Lem_CompositionDesPhi}
    Let $A\to B$ be a morphism of $R$-algebras. Then 
    $$ (A^{\phi_1}_B)^{\phi_2}_B \subseteq A^{\phi_2 \circ \phi_1}_B $$
\end{lemma}

\begin{proof}
    By the previous remark, we can assume that $A\inj B$ is an extension of $R$-algebras, and we have $A \subseteq A^{\phi_1}_{B,R} \subseteq B$. By definition $(A^{\phi_1}_B)^{\phi_2}_B = \Delta^{-1}\Big(\big\langle \Delta(c) \mid c\in A^{\phi_1}_B\big\rangle^{\phi_2}\Big)$ and $$\Big\langle \Delta(c) \mid c\in A^{\phi_1}_B\Big\rangle = \Big\langle \Delta(c) \mid c\in \Delta^{-1}\big(\langle \Delta(a) \mid a\in A\rangle^{\phi_1}\big)\Big\rangle \subseteq \big\langle \Delta(a) \mid a\in A\big\rangle^{\phi_1}.$$
    So $(A^{\phi_1}_B)^{\phi_2}_B \subseteq \Delta^{-1}\Big(\big(\langle \Delta(a) \mid a\in A\rangle^{\phi_1}\big)^{\phi_2}\Big) = \Delta^{-1}\Big(\langle \Delta(a) \mid a\in A\rangle^{\phi_2\circ\phi_1}\Big) = A^{\phi_2 \circ \phi_1}_B.$
\end{proof}

We recover some classical properties of Lipschitz saturation and seminormalization within the general framework of $\phi$-saturations, but restrict ourselves to those properties that will be needed in the remainder of the paper.

\begin{definition}
    Let $R$ and $B$ be rings. We define the category $\mathcal{C}_{B,R}$ whose objects are commutative diagrams of the form $R\to A\to B$ and whose morphisms are ring homomorphisms $A_1\to A_2$ making the following diagram commute:
\[\begin{tikzcd}
	R & {A_1} & B \\
	R & {A_2} & B
	\arrow[from=1-1, to=1-2]
	\arrow["\Id", from=1-1, to=2-1]
	\arrow[from=1-2, to=1-3]
	\arrow[from=1-2, to=2-2]
	\arrow["\Id"', from=1-3, to=2-3]
	\arrow[from=2-1, to=2-2]
	\arrow[from=2-2, to=2-3]
\end{tikzcd}\]
In other words, a morphism in $\mathcal{C}_{B,R}$ is given by a sequence $R\to A_1\to A_2\to B$.
\end{definition}

We consider a saturation functor.

\begin{prop}\label{Prop_Functoriality}
    Let $R$ and $B$ be rings and let $\phi$ be an ideal operation. Then the mapping
    \[
    \begin{array}{rccc}
        \text{Sat :} & \mathcal{C}_{B,R} & \to & \textbf{Arr}(\mathcal{C}_{B,R}) \\
        & A & \mapsto & \Big( A \to A^\phi_{B,R}\Big)
    \end{array}
    \]
    is a functor.
\end{prop}

\begin{proof}
    By Remark \ref{Rmq_IdpluspetitquePhi}, the mapping Sat is well-defined. Now let $h : A_1\to A_2$ be a morphism of $\mathcal{C}_{B,R}$. Then we have a commutative diagram:
\[\begin{tikzcd}
	R & {A_1} & {A_2} & B \\
	R & {A_1^\phi} & {A_2^\phi} & B
	\arrow[from=1-1, to=1-2]
	\arrow[from=1-1, to=2-1]
	\arrow[from=1-2, to=1-3]
	\arrow[from=1-2, to=2-2]
	\arrow[from=1-3, to=1-4]
	\arrow[from=1-3, to=2-3]
	\arrow[from=1-4, to=2-4]
	\arrow[from=2-1, to=2-2]
    \arrow["\subseteq"{marking, allow upside down}, draw=none, from=2-2, to=2-3]
	\arrow["\subseteq"{marking, allow upside down}, draw=none, from=2-3, to=2-4]
\end{tikzcd}\]
Indeed, if we write $g_1 :A_1\to B$ and $g_2 :A_2\to B$, with $g_1 = g_2\circ h$ then
$$\Delta(g_1(A_1)) = \Delta(g_2\circ h(A_1)) \subseteq \Delta(g_2(A_2))$$ and by Remark \ref{Rmq_ReecritureDefSaturation} and the fact that $\phi$ is order-preserving, we deduce that $(A_1)^\phi_{B,R} \subseteq (A_2)^\phi_{B,R}$. The identity and composition properties of the functor Sat then follow directly.
\end{proof}

\begin{prop}[Universal property]
    Let $A\xrightarrow{\,\text{{\color{white}f}}\,} C\xrightarrow{\,\text{g}\,} B$ be morphisms of $R$-algebras and let $\phi$ be an ideal operation. We say that $C$ is $\phi$-saturated in $B$ if $C = C^\phi_B$. In this case, we have the inclusion
    $$A^\phi_B\subseteq C.$$
    If moreover $\phi\circ\phi = \phi$, then $A^\phi_B$ is the smallest $\phi$-saturated $R$-algebra lying between $A$ and $B$, and we have the equivalence
    $$ A^\phi_B = C^\phi_B \iff g(C)\subseteq A^\phi_B.$$
\end{prop}

\begin{proof}
    By Proposition \ref{Prop_Functoriality}, we have $A^\phi_B \subseteq C^\phi_B$. So if $C$ is saturated, we have $A^\phi_B \subseteq C$. Assume now that $\phi\circ\phi = \phi$, then $A^\phi_B$ is $\phi$-saturated by Lemma \ref{Lem_CompositionDesPhi} and Lemma \ref{Lem_SatProcessIsOrderPreserving}. If $g(C)\subseteq A^\phi_B$, then by functoriality of Sat, we get $C^\phi_B  = g(C)^\phi_B\subseteq (A^\phi_B)^\phi_B = A^{\phi}_B$. Conversely, if we assume that $A^\phi_B = C^\phi_B$, then we get
    $$ g(C)\inj C^\phi_B = A^\phi_B.$$
\end{proof}

Let us prove now that, for some ideal operations, saturations and localization commute. We say that an ideal operation $\phi$ on a ring $R$ commutes with localization if for every multiplicative set $S\subset R$ the ideal operation extends on $S^{-1}R$ and if for every ideal $I$ of $R$, we have $S^{-1}I^\phi = (S^{-1}I)^\phi$.

\begin{prop}\label{Prop_SatLocalizCommute}
Let $R$ be a ring, let $A\to B$ be a morphism of $R$-algebras, let $\phi$ be an ideal operation which commutes with localization and let $S$ be a multiplicative subset of $A$. Then
$$ 
S^{-1}(A_B^\phi)=\bigl(S^{-1}A\bigr)^\phi_{S^{-1}B}.
$$
\end{prop}

\begin{proof}
For convenience, we will write $A_S := S^{-1}A$ and $B_S := S^{-1}B$. Let $T:=\{s\otimes t\mid s,t\in S\}$, then we have $B_S\otimes_{R}B_S\simeq T^{-1}(B\otimes_R B)$. Denote by $I_{\Delta,S}$ the diagonal ideal associated with
$A_S \to B_S$. We first observe that $I_{\Delta,S}=T^{-1}I_\Delta$. Indeed, the inclusion $T^{-1}I_\Delta \subseteq I_{\Delta,S}$ is immediate. For the other one, if $a\in A$ and $s\in S$, then
$$
\Delta\left(\frac{a}{s}\right)
=
\frac{(1\otimes s)\Delta(a)-(1\otimes a)\Delta(s)}
{s\otimes s}
\in T^{-1}I_\Delta.
$$
Since we have assumed that $\phi$ commutes with localization, we obtain $I^\phi_{\Delta,S} = (T^{-1}I_{\Delta})^\phi = T^{-1}I_{\Delta}^\phi$.
Let $b\in A_B^\phi$ and $s\in S$. Using the same identity, we have
$$
\Delta\left(\frac{b}{s}\right)
=
\frac{(1\otimes s)\Delta(b)-(1\otimes b)\Delta(s)}
{s\otimes s}.
$$
Since $\Delta(b)\in I_{\Delta}^\phi$ and $\Delta(s)\in I_\Delta\subseteq I_{\Delta}^\phi$, we get
$\Delta(b/s)\in T^{-1}I_{\Delta}^\phi=I_{\Delta,S}^\phi$, and hence
$b/s\in(A_S)^\phi_{B_S}$. This proves the inclusion
$S^{-1}(A_B^\phi)\subseteq(A_S)^\phi_{B_S}$.

Conversely, let $b/s\in(A_S)^\phi_{B_S}$. The preceding identity, together
with $\Delta(s)\in I_\Delta \subset I_\Delta^\phi$ and $\Delta(b/s)\in I_\Delta^\phi$, shows that $\Delta(b)\in T^{-1}I_\Delta^\phi$. Consequently,
there exist $u,v\in S$ such that $(u\otimes v)\Delta(b)\in I_\Delta^\phi$. We then have
$$
(uv\otimes 1)\Delta(b)
=
(u\otimes v)\Delta(b)
+
(u\otimes 1)\Delta(v)\Delta(b).
$$
The first term belongs to $I_\Delta^\phi$, while the second belongs to
$I_\Delta\subseteq I_\Delta^\phi$. Thus $(uv\otimes 1)\Delta(b)\in I_\Delta^\phi$. Since
$\Delta(uv)\in I_\Delta$, it follows that
$$
\Delta(uvb)
=
(uv\otimes 1)\Delta(b)
+
(1\otimes b)\Delta(uv)
\in I_\Delta^\phi.
$$
Therefore $uvb\in A_B^\phi$, and
$$
\frac{b}{s}=\frac{uvb}{uvs}\in S^{-1}(A_B^\phi).
$$
This proves the reverse inclusion.
\end{proof}

For ideal operations that commutes with localization, we generalize \cite[Proposition 2.6]{GrecoTraverso80}. In particular, we get that $\phi$-saturation can be check at height 1 prime ideals if $A$ is $(S_2)$.

\begin{theorem}
Let $R$ be a ring, let $A\inj B$ be a extension of $R$-algebras with $A$ Noetherian, let $\phi$ be an ideal operation which commutes with localization. Then the following conditions are equivalent:

\begin{enumerate}
    \item $A=A_B^\phi$.

    \item For every prime ideal $\mathfrak p$ of $A$, one has
    $A_{\mathfrak p}=(A_{\mathfrak p})^\phi_{B_{\mathfrak p}}$.

    \item For every maximal ideal $\mathfrak m$ of $A$, one has
    $A_{\mathfrak m}=(A_{\mathfrak m})^\phi_{B_{\mathfrak m}}$.

    \item For every $\mathfrak p\in\operatorname{Ass}_A(B/A)$, one has
    $A_{\mathfrak p}=(A_{\mathfrak p})^\phi_{B_{\mathfrak p}}$.
\end{enumerate}

Assume moreover that $B$ is contained in the total ring of fractions of $A$. Then the preceding conditions are equivalent to:

\begin{enumerate}
    \item[5.] For every prime ideal $\mathfrak p$ of $A$ such that
    $\operatorname{depth}(A_{\mathfrak p})=1$, one has
    $A_{\mathfrak p}=(A_{\mathfrak p})^\phi_{B_{\mathfrak p}}$.
\end{enumerate}

If, moreover, $A$ satisfies Serre's condition $(S_2)$, the preceding conditions are equivalent to:

\begin{enumerate}
    \item[6.] For every prime ideal $\mathfrak p$ of height one, one has
    $A_{\mathfrak p}=(A_{\mathfrak p})^\phi_{B_{\mathfrak p}}$.
\end{enumerate}
\end{theorem}

\begin{proof}
    By proposition \ref{Prop_SatLocalizCommute}, the $\phi$-saturation commutes with localization, so condition $1)$ implies all the other conditions. In particular, we have $1)\implies 2)\implies 3)$. Now, if we assume condition $3)$, we have 
    $$ (A^\phi_B)_{\mathfrak m} = \bigl(A_{\mathfrak m}\bigr)^\phi_{B_{\mathfrak m}} = A_{\mathfrak m} $$
    for every maximal ideal $\mathfrak m$ of $A$. Since two $A$-modules are equal if and only if their localizations at every maximal ideal of $A$ are equal (see \cite[Tag~00HN]{StacksProject}), we obtain $A^\phi_B=A$. Thus $3)$ implies $1)$. The implication $4) \implies 1)$ goes as follow: Suppose $4)$ is satisfied and assume that $A\neq A^\phi_B$. Then one can consider an element $\p \in \Ass_A(A^\phi_B/A)$. By \ref{Prop_SatLocalizCommute}, we have that $(A^\phi_B)_\p$ is the $\phi$-saturation of $A_\p$ in $B_\p$. Hence, by $4)$, we have $(A^\phi_B/A)_\p \simeq (A^\phi_B)_\p/A_\p = 0$, which is a contradiction with $A^\phi_B/A \neq 0$ and $\p\in \Ass_A(A^\phi_B/A)$.

    Assume now that $B$ is contained in the total ring of fractions $K$ of $A$ and let us prove $(5)\Rightarrow(1)$. Suppose that $A\neq A^\phi_B$ and let $\p\in\Ass_A(A^\phi_B/A)$. There exists $b\in A^\phi_B$ such that $\p=\text{Ann}_A(b+A)$. Since $A^\phi_B\subseteq B\subseteq K$, we can write $b=y/x$ with $x,y\in A$ and $x$ a non-zero-divisor. We then have
    \[
    \p =\{a\in A\mid ay/x\in A\}
    \]
    Hence $\p$ is the annihilator of the class of $y$ in $A/xA$, and therefore $\p\in\Ass_A(A/xA)$. Since $x$ is a non-zero-divisor, it follows that $\text{depth}(A_{\p})=1$. By condition (5) and Proposition \ref{Prop_SatLocalizCommute}, we obtain $(A^\phi_B)_{\p}=(A_{\p})^\phi_{B_{\p}}=A_{\p}$, which contradicts $\p\in\Ass_A(A^\phi_B/A)$. Hence $A=A^\phi_B$.

    Finally, if $A$ is $(S_2)$, then $\text{ht}(\p)=1$ is equivalent to $\text{depth}(A_\p) =1$. Hence $5)\iff 6)$.
\end{proof}

\begin{rmq}
    Let $X$ be a reduced complex algebraic variety and let $\pi^+:X^+\to X$ be its seminormalization. If $X$ satisfies Serre's condition $(S_2)$, then $X^+$ also satisfies $(S_2)$. Indeed, an $(S_2)$ variety has the following connectedness property: removing a closed subset of codimension at least two from a connected sufficiently small neighborhood does not disconnect it \cite[\S6]{Kollar2016Maps}, see also \cite{Hartshorne1962Connectedness}. Since $\pi^+$ is a homeomorphism, $X^+$ has the same
    connectedness property as $X$. Finally, $X^+$ is seminormal, and for a seminormal complex variety this connectedness property is equivalent
    to Serre's condition $(S_2)$ \cite[\S6]{Kollar2016Maps}. Hence $X^+$ is $(S_2)$.
\end{rmq}

\begin{ex}
    This is not the case for the Lipschitz saturation. Let $R=\C$, $A=\C[x^7, y^2, xy^9, x^4y^8]$ and $B=\C[x,y]$. In other words, $A$ is the coordinate ring of the affine toric variety with semigroup $\Gamma$ generated by $(7,0),(0,2),(1,9)$ and $(4,9)$. One can easily compute the semigroup $\Gamma^s$ of its Lipschitz saturation using \cite{FGST2020LipschitzAnAlgebraicApproch}. See that $A$ is $(S_2)$ because, if we denote by $\sigma = R_{\geq 0}.\Gamma$ the cone generated by $\Gamma$, then we have $$\Gamma = \bigcap_{\tau \text{ face of } \sigma} \Gamma+\Z(\Gamma\cap\tau)$$
    which mean that $A = \cap_{\text{ht}(p)=1} A_\p$. However, this is not the case for $\Gamma^s$ because the element $(2,8)$ belongs to $\Gamma^s + \Z(7,0)$ and to $\Gamma^s + \Z(0,2)$ while $(2,8)\notin \Gamma^s$.
\end{ex}

\subsection{Saturation of a coherent sheaf.}\label{SubsectionCoherence}

In order to define the notion of saturation of analytic varieties, we must verify that a saturation produces a coherent sheaf.

\begin{definition}
    Let $\phi$ be an ideal operation and $X$ a complex analytic variety. We say that $\phi$ is \textit{coherent} if, for every coherent sheaf of ideal $\Ical\subset \Ocal_X$, the sheaf $\Ical^\phi$ such that $\Gamma(U,\Ical^\phi) := \Gamma(\Ical,U)^\phi$ on every open set $U\subset X$, is well-defined and is a coherent $\Ocal_X$-module.
    
  
\end{definition}

Let us see that a saturation obtained from a coherent ideal operation, constructs coherent algebras.

\begin{lemma}\label{Lem_Coherence}
    Let $X$ be an analytic variety and let $\pi : \overline{X} \to X$ be its normalization. Let $\phi$ be a coherent ideal operation. Then, the sheaf $\Ocal_{X}^\phi$ defined by
    $\Ocal_{X}^\phi(U) := \Ocal_{X}(U)^\phi_{\Ocal_{\overline{X}}(\pi^{-1}(U))}$
    for every open set $U$ of $X$, is a coherent $\Ocal_X$-algebra.
\end{lemma}

\begin{proof}
We consider the diagonal ideal $\Ical$ which is defined as the kernel of the morphism
\[
\begin{tikzcd}
	{\Ocal_{\overline{X}\times \overline{X}}} & {\Ocal_{\overline{X}\times_X \overline{X}}}
	\arrow[from=1-1, to=1-2]
\end{tikzcd}
\]
and which is coherent as an $\Ocal_{\overline{X}\times \overline{X}}$-module. By assumption, the ideal sheaf $\Ical^{\phi}$ is a coherent $\Ocal_{\overline{X}\times \overline{X}}$-module. Recall that we assume ideal operations to satisfy $\Ical\subset \Ical^\phi$. Then, by \cite[Consequence 6, p.238]{grauert2012coherent}, $\Ical^\phi/\Ical$ is a coherent $\Ocal_{\overline{X}\times \overline{X}}/\Ical$-module. Moreover,
$\Ocal_{\overline{X}\times \overline{X}}/\Ical \simeq \Ocal_{\overline{X}\times_X \overline{X}}$.

Since the normalization morphism is finite, $\Ocal_{\overline{X}}$ is a coherent $\Ocal_X$-module. Moreover, the natural morphism $\overline{X}\times_X\overline{X}\to X$ is finite, and therefore $\Ocal_{\overline{X}\times_X\overline{X}}$ is also a coherent $\Ocal_X$-module by \cite[p.240]{grauert2012coherent}. By definition of coherence, $\Ical^\phi/\Ical$ is locally the cokernel of a morphism of the form
$\Ocal_{\overline{X}\times_X \overline{X}}^{\oplus n} \to
\Ocal_{\overline{X}\times_X \overline{X}}^{\oplus m}$.
Since the cokernel of a morphism of coherent $\Ocal_X$-modules is coherent, we get that $\Ical^\phi/\Ical$ is a coherent $\Ocal_X$-module. Hence
$\Ocal_{\overline{X}\times_X \overline{X}}/(\Ical^\phi/\Ical)$ is a coherent $\Ocal_X$-module, and we can consider the map of coherent $\Ocal_X$-modules
\[
\begin{tikzcd}
	{\Delta : \Ocal_{\overline{X}}} &
	{\Ocal_{\overline{X}\times_X \overline{X}}} &
	{\Ocal_{\overline{X}\times_X \overline{X}}/(\Ical^\phi/\Ical)}
	\arrow[from=1-1, to=1-2]
	\arrow[two heads, from=1-2, to=1-3].
\end{tikzcd}
\]
It is defined by $p\mapsto p\otimes1-1\otimes p$. We can see that $\Delta$ is a morphism of $\Ocal_X$-modules because, for every $f,g\in\Ocal_{\overline{X}}$ and every $r\in\Ocal_X$, we have $\Delta(f+g)=\Delta(f)+\Delta(g)$ and
$\Delta(rf)=(r\otimes1)\Delta(f)+(1\otimes f)\Delta(r)=(r\otimes1)\Delta(f)$,
since $\Delta(r)=0$ on $\overline{X}\times_X\overline{X}$. Then the kernel of $\Delta$ is a coherent $\Ocal_X$-module, and it corresponds exactly to the $\Ocal_X$-algebra $\Ocal_X^\phi$.
\end{proof}

\begin{rmq}
With the notation of the proof, let $q:\overline{X}\times_X\overline{X}\to X$ be the natural morphism. The preceding argument shows that one can write the $\phi$-saturation as
\[
\Ocal_X^\phi
=
\ker\left(
\pi_*\Ocal_{\overline X}
\xrightarrow{\ \Delta\ }
q_*\left(
\Ocal_{\overline X\times_X\overline X}/(\Ical^\phi/\Ical)
\right)
\right),
\]
where $\Delta(f)=f\otimes1-1\otimes f$.
\end{rmq}

\section{Hölder saturations}\label{Sec_HolderSaturations}

In the rest of the paper we focus on the Hölder saturations, which are the saturations of Section~\ref{SectionSaturations} for the ideal operation of taking rational powers. The section has three parts. We first recall the rational powers of an ideal, introduced by Lejeune-Jalabert and Teissier \cite{LejeuneTeissier2008ClotureInt} in complex analytic geometry, and prove that an ideal has only finitely many distinct rational powers, with jumps at rational numbers governed by its Rees valuations. This gives the coherence of the Hölder saturations and the finiteness of the exponents, which is Theorem~A of the introduction. We then identify the Hölder saturations with the sheaves of Hölder meromorphic functions, using the characterization of rational powers by inequalities due to Lejeune-Jalabert and Teissier. Finally we prove the $1/3$-theorem: a reduced germ whose $\alpha$-Hölder saturation is smooth for some $\alpha>1/3$ is itself smooth. We refer to \cite[Section~10.5]{SwansonHuneke2006IntegralClosure} for the properties of rational powers in commutative algebra.

\begin{definition}
    Let $R$ be a Noetherian ring and let $I\subset R$ be an ideal. Let $\alpha \in \Q_{>0}$ and consider $p,q\in \N$ such that $\alpha = p/q$. Then the rational power $I_{\alpha}$ of $I$ is defined by 
    \[
    I_\alpha := \big\{ r\in R \mid r^q \in \overline{I^p} \big\}
    \]
    where $\overline{I^p}$ denote the integral closure of the ideal $I^p$.
\end{definition}

Of course, if $\alpha = 1$ we recover the notion of integral closure of an ideal. We can also recover the notion of radical of an ideal by considering $I_{> 0} := \cup_{\alpha >0} I_\alpha$. This simply comes from the fact that the condition $r^n\in I$ implies $r\in I_{\frac{1}{n}}$ (see \cite[Example 10.5.4]{SwansonHuneke2006IntegralClosure} for details).\\ 

Let us state some properties of rational powers of ideals that will be useful in the present work.

\begin{prop}[\cite{SwansonHuneke2006IntegralClosure}Proposition 10.5.2]
    Let $R$ be a Noetherian ring, $I\subset R$ an ideal and $\alpha,\beta \in \Q_{>0}$. Then
    \begin{enumerate}
        \item[(1)] The definition of $I_\alpha$ does not depend on the choice of $p$ and $q$.
        \item[(2)] If $\alpha \leq \beta$, then $I_{\beta} \subseteq I_{\alpha}$.
        \item[(3)] The ideal $I_\alpha$ is integrally closed.
    \end{enumerate}
\end{prop}

Taking the rational power of an ideal is an ideal operation in the sense of Section \ref{SectionSaturations}. Moreover, for $\alpha \leq \beta \leq 1$, we have the chain of inclusions
\[
I \subseteq \overline{I} = I_1 \subseteq I_\beta \subseteq I_\alpha \subseteq I_{>0} = \sqrt{I}
\]
This means in particular that, when $\alpha \leq 1$, this ideal operation verifies $\Id \leq (.)_\alpha$ as in Remark \ref{Rmq_IdpluspetitquePhi}. Hence, in this situation, we can consider the saturation associated to $(.)_\alpha$.

\begin{definition}
    Let $A\to B$ be a morphism of algebras and let $\alpha$ be a rational number in $(0,1]$. The $\alpha$-Hölder saturation $A^\alpha_B$ of $A$ in $B$ is the saturation $A^{(.)_\alpha}_B$ associated to the ideal operation $I\to I_\alpha$.
\end{definition}

From the properties satisfied by rational powers of ideals, we can deduce the following important relations between Hölder saturations.

\begin{prop}\label{Prop_InteractionsBetweenSaturations}
    Let $A\to B$ be a morphism of algebras such that $B\otimes B$ is Noetherian. Then Hölder saturations of $A$ in $B$ verify the following properties :
    \begin{enumerate}
        \item[(1)] For every rational numbers $\alpha \leq \beta$ in $(0,1]$, we have $A^{(.)_{\beta}} \subseteq A^{(.)_{\alpha}}$.
        \item[(2)] For every rational numbers $\alpha, \beta \in (0,1]$, we have $\big(A^{(.)_\alpha}\big)^{(.)_\beta} \subseteq A^{(.)_{\alpha*\beta}}$.
        \item[(3)] For every rational number $\alpha \in (0,1]$ we have $\big(A^{(.)_\alpha}\big)^{\overline{(.)}} = \big(A^{\overline{(.)}}\big)^{(.)_\alpha} = A^{(.)_\alpha}$.
        \item[(4)] For every rational number $\alpha \in (0,1]$, we have $\big(A^{(.)_\alpha}\big)^{\sqrt{.}} = \big(A^{\sqrt{.}}\big)^{(.)_\alpha} = A^{\sqrt{.}}$.
    \end{enumerate}
    For convenience, we wrote $A^\phi$ instead $A^\phi_B$ since every saturation here is taken in $B$.
\end{prop}

\begin{proof}
    By Lemma \ref{Lem_CompositionDesPhi}, Lemma \ref{Lem_SatProcessIsOrderPreserving} and Proposition \ref{Prop_Functoriality}, Properties $(1)$ and $(3)$ are consequences of \cite[Proposition 10.5.2]{SwansonHuneke2006IntegralClosure}. Property $(4)$ can be obtained the same way but one as to show that, for every $\alpha \in (0,1]$ and any ideal $I$, we have $(\sqrt{I})_\alpha = \sqrt{I_\alpha} = \sqrt{I}$. Recall from \cite[Example 10.5.4]{SwansonHuneke2006IntegralClosure} that $\sqrt{I} = I_{>0} := \cup_{\alpha>0} I_\alpha$. The inclusion $I_\alpha \subseteq \sqrt{I}$ is clear and we get $\sqrt{I_\alpha}\subseteq \sqrt{I}$. Moreover, since $\Id \leq (.)_\alpha$, we have $\sqrt{I} \subseteq \sqrt{I_\alpha}$. Hence $\sqrt{I} = \sqrt{I_\alpha}$. We also have $\sqrt{I} \subseteq (\sqrt{I})_\alpha = (I_{\frac{1}{p_\alpha*n}})_\alpha$ for some large $n$, where $\alpha = p_\alpha/q_\alpha$. By property $(2)$ we have $(I_{\frac{1}{p_\alpha*n}})_\alpha \subseteq I_{\frac{1}{q_\alpha*n}} = \sqrt{I}$ for $n$ large enough. Hence $\sqrt{I} = (\sqrt{I})_\alpha$.

    Concerning Property $(2)$, we have to show that, for every ideal $I$ in a noetherian ring $R$, we have $(I_{\alpha})_{\beta} \subseteq I_{\alpha*\beta}$. If we write $\alpha = p_\alpha/q_\alpha$ and idem for $\beta$, then an element $f$ of $R$ belongs to $(I_{\alpha})_{\beta}$ if and only if, for every discrete valuation ring $V$ lying between $R$ and its total ring of fraction, we have $f^{q_\beta} \in I_{\alpha}^{p_\beta} V$. Since $V$ is a DVR, we can consider $g\in I_{\alpha}$ such that $I_{\alpha}V = \langle g \rangle$. Then we have $g^{q_\alpha} \in I^{p_\alpha}V$. Hence we get $f^{q_\beta q_\alpha} \in \langle g^{q_\alpha} \rangle^{p_\beta} V$ which implies $f^{q_\beta q_\alpha} \in (I^{p_\alpha})^{p_\beta} V$ and so $f\in I_{\alpha*\beta}$.
\end{proof}

\begin{rmq}
    The converse inclusion of $(2)$ is not true in general. Take $A=\C[x,y]/\langle y^2-x^3\rangle$ and let $B$ be the normalization of $A$. Then we will see in Example \ref{Ex_Cusp} that $A^{(.)_{\alpha}}_{B} = B$ if $\alpha \in (0,\frac{1}{3}]$ and $A^{(.)_{\alpha}}_{B} = A$ if $\alpha > \frac{1}{3}$. Hence 
    $$ \big(A^{(.)_{1/2}}_{B}\big)^{(.)_{1/2}}_{B} = A \text{ and }A^{(.)_{1/4}}_{B} = B.$$\smallskip
\end{rmq}

We have proven in Proposition \ref{Prop_SatLocalizCommute} that a saturation commutes with localization if the ideal operation does. Therefore, we prove the following lemma.

\begin{lemma}
    Let $R$ be a Noetherian ring, let $I$ be an ideal of $R$ and let $S$ be a multiplicative subset of $R$. Let $\alpha=p/q\in\mathbb{Q}_{> 0}$, then
    $$
    S^{-1}I_\alpha = \bigl(S^{-1}I\bigr)_\alpha.
    $$
\end{lemma}

\begin{proof}
    By \cite[Proposition 1.1.4]{SwansonHuneke2006IntegralClosure} integral closure of ideals commute with localization, so we have
    $$
    \overline{\bigl(S^{-1}I\bigr)^p} = S^{-1}\overline{I^p}.
    $$
    Let $x\in I_\alpha$ and $s\in S$. By definition, $x^q\in\overline{I^p}$, and therefore $(x/s)^q\in S^{-1}\overline{I^p}$. It follows that $x/s\in(S^{-1}I)_\alpha$. This proves the inclusion $S^{-1}I_\alpha\subseteq(S^{-1}I)_\alpha$.
    
    Conversely, let $z=r/s\in(S^{-1}I)_\alpha$, with $r\in R$ and $s\in S$. Then $r^q/s^q = z^q \in \overline{\bigl(S^{-1}I\bigr)^p} = S^{-1}\overline{I^p}$. So, we can consider $u\in S$ such that $ur^q\in\overline{I^p}$. Since $\overline{I^p}$ is an ideal, we obtain $(ur)^q = u^{q-1}(ur^q) \in \overline{I^p}$. Thus $ur\in I_\alpha$, and hence
    $$
    z=\frac{r}{s}=\frac{ur}{us}\in S^{-1}I_\alpha.
    $$
\end{proof}

From Proposition \ref{Prop_SatLocalizCommute}, we deduce that

\begin{prop}
    Let $R$ be a ring, let $A\to B$ be a morphism of $R$-algebras and let $S$ be a multiplicative subset of $A$. Then, for every $\alpha\in (0,1]\cap \Q$, we have
    $$ 
    S^{-1}(A_B^{(.)_\alpha})=\bigl(S^{-1}A\bigr)^{(.)_\alpha}_{S^{-1}B}.
    $$
\end{prop}

\subsection{Hölder saturations of complex analytic varieties}\label{Subsec_HolderSatAnalyticVarieties}

Let $X$ be a complex analytic or algebraic variety and let $\overline{X}$ be its normalization. By \cite[Proposition 4.11]{LejeuneTeissier2008ClotureInt}, the rational power operation is a coherent ideal operation as defined is Subsection \ref{SubsectionCoherence}. So by Lemma \ref{Lem_Coherence}, the sheaf $\Ocal^\alpha_{X}$ defined by $\Ocal^\alpha_{X}(U) := \Ocal_{X}(U)^\alpha_{\Ocal_{\overline{X}}(\pi^{-1}(U))}$, for every $U\subset X$, is a coherent $\Ocal_X$-algebra. Hence, we can consider the variety $X^{\alpha}$ such that $(|X^\alpha|,\Ocal_{X^{\alpha}}) := (|X|,\Ocal^{\alpha}_X)$ which is called the $\alpha$-Hölder saturation of $X$.

\subsubsection{Finiteness of Hölder exponents}\label{Subsec_FinitenessExponents}

Let $(X,x)$ be a germ of analytic space and let $\Ical\subset \Ocal_{X,x}$ be an ideal. Then one can consider the normalized blow-up
\[
Z=\Projan\left(\bigoplus_n \overline{\Ical^n}t^n\right)\to X.
\]
Let $E$ be the pull-back of the exceptional divisor in $Z$ and let $E=\cup_{i=1}^m E_i$ be its irreducible components. Then, the $E_i$'s are prime divisors on a normal space, so the rings $\Ocal_{Z,E_i}$ are DVRs. In fact, the valuations $v_{E_i}$ they induce correspond to the \textit{Rees valuations} of $\Ical$ (see \cite[Chapter 10]{SwansonHuneke2006IntegralClosure}). Since we work locally, there are only finitely many such irreducible components. We know from \cite{LejeuneTeissier2008ClotureInt} that an element $f$ belongs to the integral closure of $\Ical$ if and only if its composition by the normalized blow-up morphism belongs to $E$. Hence $f\in\overline{\Ical}$ if and only if $v_{E_i}(f)\geq v_{E_i}(\Ical)$ for every $i$. Using this crucial fact, together with the finiteness of Rees valuations, we deduce that there are only finitely many distinct rational powers $\Ical_\alpha$ for $\alpha\in(0,1]$. This proof comes from \cite[Proposition 10.5.5]{SwansonHuneke2006IntegralClosure}.

\begin{prop}\label{Prop_FinitenessOfRationalIdeals}
    Let $(X,x)$ be a germ of analytic space and let $\Ical\subset\Ocal_{X,x}$ be an ideal. Let $v_i:=v_{E_i}$ be the Rees valuations of $\Ical$. Set $e$ to be the least common multiple of $v_1(\Ical),\dots,v_m(\Ical)$. Then, for every $n\in\N^*$ and every $\alpha\in(\frac{n-1}{e},\frac{n}{e}]\cap\Q$, we have $\Ical_{\frac{n}{e}}=\Ical_\alpha\subseteq\Ical_{\frac{n-1}{e}}$.
\end{prop}

\begin{proof}
    Let $e_i=v_i(\Ical)$ and write $e=e_id_i$. Let $\alpha\in\Q_{>0}$, we want to prove that $\Ical_\alpha=\Ical_{\frac{n}{e}}$ where $n=\lceil e\alpha\rceil$. Since $\frac{n}{e}\geq\alpha$, we know that $\Ical_{\frac{n}{e}}\subseteq\Ical_\alpha$. To prove the opposite inclusion, let us consider $f\in\Ical_\alpha$. By the discussion above on Rees valuations, we have $v_i(f)\geq\alpha e_i=\alpha e/d_i$ for all $i=1,\dots,m$. So we have $d_iv_i(f)\geq\alpha e$, for every $i$, and since $d_iv_i(f)\in\Z$, we get $d_iv_i(f)\geq\lceil\alpha e\rceil=n$. Therefore $ev_i(f)\geq ne_i=nv_i(\Ical)$ for all $i=1,\dots,m$, which is equivalent to $v_i(f^e)\geq v_i(\Ical^n)$ for all $i$. This implies $f\in\Ical_{\frac{n}{e}}$.
\end{proof}

\begin{rmq}
    We have
     $$\overline{\Ical} = \Ical_1 \subseteq \Ical_{\frac{e-1}{e}}\subseteq \dots \subseteq \Ical_{\frac{1}{e}} \subseteq \Ical_{>0} = \sqrt{\Ical}.$$
     Recall that $\Ical_{>0} := \cup_{\alpha>0} \Ical_{\alpha}$. By the previous proposition, we have that every power of ideal $\Ical_\alpha$ belongs to $\Ical_{1/e}$. So $$\Ical_{1/e} = \Ical_{>0} = \sqrt{\Ical}.$$
\end{rmq}

By applying this proposition to the diagonal ideal $\Ical\subset \Ocal_{\overline{X}\times \overline{X}}$, we get that every germ of complex analytic variety admits a finite number of $\alpha$-Hölder saturations.

\begin{prop}[local case]\label{Prop_FiniteNumberOfAlphasCoherentLOCAL}
    Let $(X,x)$ be a germ of complex analytic variety. Then there exist rational numbers $0<\alpha_1<\dots<\alpha_r = 1$ such that 
    \[
    \forall i\in\{1,\dots,r-1\} \text{\hspace{0.1cm} }\forall \alpha \in (\alpha_i,\alpha_{i+1}]\cap \Q \text{\hspace{0.5cm} }\Ocal^\alpha_{X,x} \simeq \Ocal^{\alpha_{i+1}}_{X,x}
    \]
    and
    \[
    \forall \alpha \in (0,\alpha_1] \cap \Q\text{\hspace{0.5cm} }\Ocal^{\alpha}_{X,x}\simeq \Ocal_{X,x}^{\alpha_1} \simeq \Ocal^{+}_{X,x}
    \]
\end{prop}

Of course, here, the finite set $\{\alpha_1,\dots,\alpha_r\}$ depends on the point $x$. This result extends to the global case for algebraic varieties and compact analytic spaces, but one can not adapt the previous discussion for general analytic spaces since those do not necessarily admit a finite number of irreducible components.

\begin{prop}[global case]\label{Prop_FiniteNumberOfAlphasCoherentGLOBAL}
    Let $X$ be a complex algebraic variety (resp. a compact analytic space). Then there exists rational numbers $0<\alpha_1<\dots<\alpha_r = 1$ such that 
    \[
    \forall i\in\{1,\dots,r-1\} \text{\hspace{0.1cm} }\forall \alpha \in (\alpha_i,\alpha_{i+1}]\cap \Q \text{\hspace{0.5cm} }\Ocal^\alpha_{X} \simeq \Ocal^{\alpha_{i+1}}_{X}
    \]
    and
    \[
    \forall \alpha \in (0,\alpha_1] \cap \Q\text{\hspace{0.5cm} }\Ocal^{\alpha}_{X}\simeq \Ocal_X^{\alpha_1} \simeq \Ocal^{+}_{X}
    \]
\end{prop}

\begin{proof}
In the algebraic case, this follows directly from the preceding discussion, since the normalized blow-up has only finitely many exceptional irreducible components.

In the compact analytic case, every point admits a neighborhood on which the conclusion of the local proposition holds. By compactness, we may choose a finite covering by such neighborhoods. Taking the union of the corresponding finite sets of exponents gives a finite set $0<\alpha_1<\cdots<\alpha_r=1$ which works on every member of the covering, and hence on $X$.
\end{proof}



\subsection{Hölder meromorphic functions}

In this section, we extend the theory to non-necessarily rational $\alpha$s and we study the connection between Hölder saturations and Hölder meromorphic functions. We first state general facts about $\alpha$-Hölder meromorphic functions before making the connection with Hölder saturations.

\begin{definition}
     Let $X$ be a complex analytic variety, let $\Ocal_{X,0}$ be the ring of germs of holomorphic functions at $0$ and let $\Mcal_{X,0}$ be the ring of germs of meromorphic functions at $0$.
     \begin{enumerate}
         \item For any $\alpha\in\R_{> 0}$, we denote by $\Mcal^\alpha_{X,0}$ the set of elements $f\in \Mcal_{X,0}$ such that there exists a neighborhood of $0$ where $f$ can be extended as a $\alpha$-Hölder function.

        \item We denote by $\Mcal^{c}_{X,0}$ the set of elements $f\in \Mcal_{X,0}$ such that there exists a neighborhood of $0$ where $f$ can be extended as a continuous function.
     \end{enumerate}
\end{definition}


It is clear that one can consider the sheaf $\Mcal^\alpha_{X}$ (resp. $\Mcal^{c}_{X}$) of $\alpha$-Hölder (resp. continuous) meromorphic functions. For $\alpha\leq \beta$, we have $\Mcal^\beta_{X,0} \subseteq \Mcal^\alpha_{X,0}$ and, by considering a bounded neighborhood of $0$, one can see that $\Mcal^\alpha_{X,0}$ is a ring from the following inequalities
\[
\begin{array}{ll}
    |(fg)(x)-(fg)(y)| & \leq |f(x)||g(x)-g(y)| + |g(y)||f(x)-f(y)|\medskip\\
     & \leq C_1|x-y|^\alpha + C_2|x-y|^\alpha
\end{array}
\]
and
\[
\begin{array}{ll}
    |(f+g)(x)-(f+g)(y)| & \leq |g(x)-g(y)| + |f(x)-f(y)|\medskip\\
     & \leq 2|x-y|^\alpha.
\end{array}
\]
Moreover, the ring $\Mcal^\alpha_{X,0}$ is a subring of $\overline{\Ocal_{X,0}}$, since Hölder functions are locally bounded.

\begin{rmq}
    If $\alpha > 1$, let us prove that every locally $\alpha$-Hölder function $f$ on $X$ is constant. Let $x\in X$ and let $U$ be a neighborhood of $x$ in $X$ such that $f$ is a $\alpha$-Hölder function on $U$. If $(X,x) = (\C^n, 0)$, then the result is well-known. Indeed, in this situation, we can assume $U$ to be convex and, for all $x,y\in U$, we can consider the subdivision $\{x_i\}_{i=1,\dots,n}$ of the segment $[x,y]$ defined by $x_i = x + i/n(y-x)$. Then we have
    \[
    |f(x)-f(y)| \leq \sum_{i=0}^{n-1} |f(x_{i})-f(x_{i+1})| \leq \sum_{i=0}^{n-1} C |x_{i}-x_{i+1}|^\alpha = \sum_{i=0}^{n-1} C \left(\frac{|x-y|}{n}\right)^\alpha = C|x-y|^\alpha n^{1-\alpha}
    \]
    which goes to $0$ when $n$ goes to infinity. So, by taking finer and finer subdivisions of $[x,y]$, we deduce that $f(x)=f(y)$ and so $f$ is constant on $U$. Now, if $x$ is a smooth point of an analytic space $X$, we get that $(U,x)$ is holomorphically isomorphic to a neighborhood $V$ of $0$ in $\C^n$, for some $n$. Hence, the result holds in this situation. If $x$ is not smooth, we know that $f$ is constant on $U\cap X^{reg}$. We can consider a sequence of smooth points of $X$ converging to $x$ and, by continuity of $f$, we get that $f$ is constant on $U$.
\end{rmq}

For this reason, we will usually assume $\alpha \leq 1$. Now, let us give a small proof of the following lemma, which comes from \cite[Theorem 7.2]{LejeuneTeissier2008ClotureInt}.

\begin{lemma}\label{Lem_InequalityPowerIdeal}
    Let $f\in \Ocal_{X,0}$ and let $I\subseteq \Ocal_{X,0}$ be an ideal generated by $g_1,\dots,g_n\in \Ocal_{X,0}$. Let $\alpha\in \Q_{>0}$, then $f\in I_{\alpha}$ if and only if there exists $C\in \R$ such that $|f| \leq C.(\sup_i|g_i|)^\alpha$ in a neighborhood of $0$.
\end{lemma}

\begin{proof}
    Let us write $\alpha = p/q$. By definition, we have $f\in I_{p/q}$ if and only if $f^q \in \overline{I^p}$. Moreover, we have $I^p = \langle \prod_{i=1}^n g_i^{k_i} \mid \{k_i\}_i\text{ s.t. }\sum k_j = p\rangle$, so by \cite[Corollary 0.8]{PhamTeissier1969LipFractions}, we get that there exists $C\in \R$ such that
    $|f^q| \leq C.(\sup_{\{k_i\}_i\text{ s.t. }\sum k_j = p}|\prod_{i=1}^n g_i^{k_i}|)$ in a neighborhood of $0$. Then the conclusion comes from the following inequalities:
    \[
    (\sup_i |g_i|)^p = \sup_i|g_i^p| \leq \sup_{\{k_i\}_i\text{ s.t. }\sum k_j = p}|\prod_{i=1}^n g_i^{k_i}|
    \]
    and, for any such $\{k_i\}_i$,
    \[
    \prod_{i=1}^n g_i^{k_i}| \leq |\prod_{i=1}^n (\sup_i |g_i|)^{k_i} = (\sup_i |g_i|)^{p}.
    \]
\end{proof}

From this lemma, we get the link between Hölder saturations and Hölder meromorphic functions. This allows us to extend the definition of $X^\alpha$ to non-rational $\alpha$'s. 

\begin{theorem}\label{Th_HolderSatIsHolderFcts}
    Let $X$ be a complex analytic variety and let $\pi : \overline{X} \to X$ be its normalization. Then, for every rational number $0< \alpha \leq 1$, we have
    \[
    (X,\Mcal^\alpha_X) \simeq (X^\alpha,\Ocal_{X^\alpha})
    \]
\end{theorem}

\begin{proof}
    
    Let $f\in \Mcal_{X,0}^\alpha$ then $f$ a locally bounded meromorphic function, hence it is regular on the normalization. This mean that $f\circ \pi \in \bigoplus_{x\in \pi^{-1}(0)} \Ocal_{\overline{X},x}$. Moreover, since the normalization map is surjective, proper and continuous, the function $|f(x)-f(y)| / (\sup_i|x_i-y_i|)^\alpha$ is locally bounded if and only if its composition by the normalization map is locally bounded (see \cite[Section 1.1]{Bernard2025RelativeLipSat} for example). So, by Lemma \ref{Lem_InequalityPowerIdeal}, the Hölder condition on $f$ is equivalent to $f\circ \pi\otimes 1-1\otimes f\circ \pi \in I_\alpha$ where $I$ is the ideal $\langle x_i\circ\pi-y_i\circ\pi \rangle \subset \Ocal_{\overline{X}\times \overline{X}}$ which, by definition, mean that $f\circ\pi \in \Ocal_{X,0}^\alpha$.
    
\end{proof}

\begin{cor}\label{Cor_HolderSatMapIsHolder}
    Let $X$ be a complex analytic variety and let $\pi_\alpha : X^\alpha \to X$ be its $\alpha$-Hölder saturation. Then $\pi_\alpha$ is a locally bi-$\alpha$-Hölder meromorphic map.
\end{cor}


\begin{proof}
    Locally, we have 
    \[
    \begin{array}{ccccc}
     \Ocal_{X,0} &\inj& \Mcal^\alpha_{X,0} &\xrightarrow{\sim}& \Ocal_{X^\alpha,0}\\
     && f &\mapsto& {f\circ\pi_\alpha}
    \end{array}
    \]
    So, if $x := (x_1,\dots,x_n)$ is a system of local coordinates of $(X^\alpha,0)$. One can consider $p_i : (X^\alpha,0) \to (\C,0)$, the coordinate projection $x \mapsto x_i$. This is a regular function of $\Ocal_{X^\alpha,0}$, so there exists $f_i\in \Mcal^\alpha_{X,0}$ such that $f_i\circ \pi_\alpha = p_i$. Hence, the function $f:=(f_1,\dots, f_n)$ is $\alpha$-Hölder and $f\circ\pi_\alpha = \Id$. Moreover, we have $\pi_\alpha\circ (f\circ\pi_\alpha) = \pi_\alpha$ and since $f\mapsto f\circ\pi_\alpha$ is an isomorphism, we have $(\pi_\alpha\circ f)\circ\pi_\alpha = \Id\circ \pi_\alpha$ implies $(\pi_\alpha\circ f) = \Id$.
\end{proof}

For a compact complex analytic space or a complex algebraic variety $X$, we know from Theorem \ref{Th_HolderSatIsHolderFcts} and Proposition \ref{Prop_FiniteNumberOfAlphasCoherentGLOBAL} that there exists rational numbers $0<\alpha_1<\dots<\alpha_r=1$ such that 
\[
\forall i\in\{1,\dots,r-1\} \text{\hspace{0.1cm} }\forall \alpha \in (\alpha_i,\alpha_{i+1}]\cap \Q \text{\hspace{0.5cm}}\Mcal^\alpha_{X} \simeq \Mcal^{\alpha_{i+1}}_{X}
\]
and
\[
\forall \alpha \in (0,\alpha_1]\cap \Q \text{\hspace{0.5cm} }\Mcal^{\alpha}_{X}\simeq \Ocal_{X^+}
\]
where $X^+$ is the seminormalization of $X$. Let us recall from \cite{AndreottiNorguet1967SeminormalizationAnalytic} that $\Ocal_{X^+}$ is isomorphic to the sheaf of continuous meromorphic functions $\Mcal^{c}_{X}$. Now, let us consider a non-necessarily rational number $\alpha\in (\alpha_1,1]$. We have $\alpha\in (\alpha_i,\alpha_{i+1}]$ for some $i\in\{1,\dots,r-1\}$. Let $p/q \in \Q$ such that $\alpha_i < p/q < \alpha$. Then $p/q$ is a rational number which belongs to the interval $(\alpha_i,\alpha_{i+1}]$ so, by Proposition \ref{Prop_FiniteNumberOfAlphasCoherentGLOBAL}, we get
\[
\Mcal^{\alpha_{i+1}}_{X} \subseteq \Mcal^{\alpha}_{X} \subseteq \Mcal^{p/q}_{X} = \Mcal^{\alpha_{i+1}}_{X}.
\]
Hence $\Mcal_{X}^\alpha = \Mcal_{X}^{\alpha_i+1}$. For $\alpha\in(0,\alpha_1]$, we get $\Mcal^{\alpha}_{X}\simeq \Mcal^{c}_X$ by the same argument. In conclusion, we can generalize Proposition \ref{Prop_FiniteNumberOfAlphasCoherentGLOBAL} to sheaves of $\alpha$-Hölder meromorphic function for a non-necessarily rational $\alpha$.

\begin{theorem}\label{Theo_FiniteNumberOfAlphas}
Let $X$ be a compact complex analytic space or a complex algebraic variety, and, for every $\alpha\in(0,1]$, let $\Mcal_X^\alpha$ be the sheaf of $\alpha$-Hölder meromorphic functions. Then there exist finitely many rational numbers
$0<\alpha_1<\cdots<\alpha_r=1$ such that
\[
\forall i\in\{1,\ldots,r-1\}\ \forall\alpha\in(\alpha_i,\alpha_{i+1}]
\qquad
\Mcal_X^\alpha\simeq\Mcal_X^{\alpha_{i+1}}
\]
and
\[
\forall\alpha\in(0,\alpha_1]
\qquad
\Mcal_X^\alpha\simeq\Mcal_X^{c}.
\]
\end{theorem}

A striking consequence of this theorem is that every continuous meromorphic function is $\alpha$-Hölder for some small $\alpha$. Note that, by Proposition \ref{Prop_FiniteNumberOfAlphasCoherentLOCAL}, the same local Theorem holds for a germ of complex analytic variety $(X,x)$ with its ring of $\alpha$-Hölder meromorphic functions $\Mcal^\alpha_{X,x}$.

\subsection{The $1/3$-Theorem}

Now that we have establish the link between $\alpha$-Hölder saturation and bi-$\alpha$-homeomorphisms, we prove a Hölder regularity theorem. Something remarkable about this result is that it holds for germs of arbitrary dimension.

\begin{theorem}\label{Th_HolderSatSmoothness}
    Let $(X,0)$ be a reduced complex analytic germ and let $\alpha\in(1/3,1]$. If $(X^\alpha,0)$ is smooth, then $(X,0)$ is smooth.
\end{theorem}

\begin{proof}
    Let $d=\dim X$. Since $(X^\alpha,0)$ is smooth, we may write
    $\Ocal_{X^\alpha,0}\simeq\C\{z_1,\ldots,z_d\}$.
    Consider an embedding $(X,0)\subset(\C^N,0)$ and let
    $x_1,\ldots,x_N$ be the restrictions to $X$ of the ambient coordinates. Through the inclusion
    $\Ocal_{X,0}\subseteq\Ocal_{X^\alpha,0}$, write
    $x_j=f_j(z_1,\ldots,z_d)$, with
    $f_j\in\C\{z_1,\ldots,z_d\}$.

    Suppose that $(X,0)$ is singular. We first claim that the linear parts of the $f_j$ do not span the space of linear forms in $z_1,\ldots,z_d$. Otherwise, one could choose $d$ of them, say $f_{j_1},\ldots,f_{j_d}$, with linearly independent linear parts. By the holomorphic inverse function theorem, $f_{j_1},\ldots,f_{j_d}$ would form a system of local coordinates on $(X^\alpha,0)$. Hence
    $\Ocal_{X^\alpha,0}=\C\{f_{j_1},\ldots,f_{j_d}\}\subseteq\Ocal_{X,0}$.
    Since the reverse inclusion holds by definition of the Hölder saturation, this would give
    $\Ocal_{X,0}=\Ocal_{X^\alpha,0}$, and therefore $(X,0)$ would be smooth, a contradiction.

    Thus there exists $v\in\C^d\setminus\{0\}$ which is annihilated by the linear part of every $f_j$. Consider the arc $\gamma(t)=tv$ in $(X^\alpha,0)$ and set
    \[
        \varphi(t)
        :=
        \pi_\alpha(\gamma(t))
        =
        \bigl(f_1(tv),\ldots,f_N(tv)\bigr),
    \]
    where $\pi_\alpha:(X^\alpha,0)\to(X,0)$ is the natural morphism. If
    $m:=\min_j\{\ord_t(f_j(tv))\}$, then $m\geq2$.

    By Corollary \ref{Cor_HolderSatMapIsHolder}, $\pi_\alpha$ is a homeomorphism, so since $\gamma$ is injective, we get that $\varphi$ is injective. It follows that the greatest common divisor of all the exponents appearing with nonzero coefficient in the components of $\varphi$ is equal to $1$. Indeed, otherwise there would exist an integer $e>1$ such that $\varphi(t)=\varphi(\xi t)$ for every $e$-th root of unity $\xi$, contradicting the injectivity of $\varphi$.

    Let $q$ be the smallest exponent appearing with nonzero coefficient in one of the components of $\varphi$ which is not divisible by $m$. Such an exponent exists by the preceding observation, and
    $q\geq m+1\geq3$. Let $\xi$ be a primitive $m$-th root of unity. By the definition of $q$, all terms of order strictly smaller than $q$ cancel in
    $\varphi(t)-\varphi(\xi t)$, while the terms of order $q$ do not. Thus
    $\ord_t(\varphi(t)-\varphi(\xi t))=q$. In particular, there exist a constant $C_1>0$ and a neighborhood of $0$ such that
    \[
        \|\varphi(t)-\varphi(\xi t)\|
        \leq C_1|t|^q.
    \]

    On the other hand,
    $\|\gamma(t)-\gamma(\xi t)\|
    =|1-\xi|\,\|v\|\,|t|$.
    Since $\pi_\alpha^{-1}$ is $\alpha$-Hölder by Corollary \ref{Cor_HolderSatMapIsHolder}, there exists $C_2>0$ such that
    \[
        |1-\xi|\,\|v\|\,|t|
        \leq
        C_2\|\varphi(t)-\varphi(\xi t)\|^\alpha
        \leq
        C_2C_1^\alpha |t|^{\alpha q}.
    \]
    Since $|1-\xi|\,\|v\|>0$, this inequality can hold for arbitrarily small nonzero $t$ only if $\alpha q\leq1$. As $q\geq3$, we obtain $\alpha\leq1/3$, contradicting the assumption $\alpha>1/3$. Hence $(X,0)$ is smooth.
\end{proof}

\begin{rmq}
    The interval $(1/3,1]$ is sharp because we have seen in Example \ref{Ex_Cusp} that, for the cusp $(X,0)$ with local ring $\C\{t^2,t^3\}$, we have $(X^\alpha,0) = (\overline{X},0) = (\C,0)$ for every $\alpha \in (0,1/3]$.
\end{rmq}

\section{Computing Hölder saturations of curves}\label{Sec_ComputingHolderSat}

In this section, we show how to compute Hölder saturations of any reduced analytic curve.

\subsection{Hölder saturations and generic projections.}

Let $(X,0) \subset (\C^N,0)$ be a reduced and equidimensional germ of complex analytic variety. What is know as Whitney's cone $C_5(X,0)$ is the set of elements $v\in \C^N$ such that there exist sequences of pairs of distinct points $\{p_i\},\{q_i\} \subset X\backslash\{0\}$, converging to $0$, and a sequence $\{\lambda_i\} \subset \C^*$ such that 
$$ \lambda_i(p_i-q_i) \to v.$$

\begin{definition}\label{Def_GenericProjection}
    A linear projection $\pi : (\C^N,0) \to (\C^M,0)$ is called \textit{generic} with respect to $X$ if
    $$ \ker(\pi)\cap C_5(X,0) = \{0\}$$
\end{definition}

Suppose now that $(X,0) \subset (\C^N,0)$ is a germ of reduced singular curve. 
  We recall the relation between general projections and the Lipschitz saturation which is stated in \cite[Proposition VI.3.1]{BrianconGallicoGranger1980DeformationsEqu} and \cite[Proposition 1.4.11]{FGST2020LipschitzAnAlgebraicApproch}.


\begin{prop}
    Let $\pi : (\C^N,0) \to (\C^2,0)$ be a generic linear projection. Then
    \[X^s \simeq \pi(X)^s\]
\end{prop}

We saw in point $(3)$ of Lemma \ref{Prop_InteractionsBetweenSaturations} that the $\alpha$-Hölder saturation of $X$ is the same as the $\alpha$-Hölder saturation of $X^s$, so we get the following consequence.

\begin{cor}\label{Cor_HolderStableUnderProjection}
    Let $\pi : (\C^N,0)\to (\C^2,0)$ be a generic linear projection. Then 
    \[X^\alpha \simeq \pi(X)^\alpha\]
    for every $\alpha\in(0,1]$.
\end{cor}


Thus, one can always compute the Hölder saturation of a curve by first taking a generic plane projection.

\subsection{Hölder saturation for irreducible complex curves}\label{Subsec_IrreducibleCurves}

Let $(X,0)$ be an irreducible curve germ and write
$A:=\Ocal_{X,0}\subset\overline A=\C\{t\}$. If $m$ is the multiplicity of $X$, we can choose the uniformizer $t$ in such a way that there exists a transverse parameter $x\in A$ satisfying $x=t^m$. Following \cite[Chapter VI]{BrianconGallicoGranger1980DeformationsEqu}, we call such a $t$ a \textit{good uniformizer}. For an element $f=\sum_{p\geq0}a_pt^p\in A$ viewed in $\overline{A}$, let $\operatorname{Ex}(f):=\{p\in\N\mid a_p\neq0\}$. Then the \textit{set of exponents} of $A$ is
\[
E(A):=\bigcup_{f\in A}\operatorname{Ex}(f).
\]

Notice that $E(A)$ is different in general from the valuative semigroup $\Gamma(A)=\{\ord_t(f)\mid 0\neq f\in A\}$. Indeed, $E(A)$ also remembers
the exponents of the non-leading terms which appear in the elements of
$A$. Recall from \cite{BrianconGallicoGranger1980DeformationsEqu} or \cite{FGST2020LipschitzAnAlgebraicApproch}, that one can associate to $E(A)$ its characteristic sequence in the following way. Let $\beta_0:=\min(E(A)\setminus\{0\})$ and $e_0:=\beta_0$. If $e_i>1$, we define recursively $\beta_{i+1}:=\min\{p\in E(A)\mid e_i\nmid p\}$ and $e_{i+1}:=\gcd(\beta_0,\ldots,\beta_{i+1})$.
We stop when $e_g=1$. The sequence
$\{\beta_0,\ldots,\beta_g\}$ is called the \textit{characteristic
sequence} of $A$. The saturated set associated with $E(A)$ is then
\[
E(A)^s = \{0\}\cup \bigcup_{i=0}^g \{p\in\N^*\mid p\geq\beta_i,\ e_i\mid p\}.
\]

One can be visualize a saturated set of exponents as in the following figure: a succession of groups $e_0\Z\subset e_1\Z \subset \dots$ lying between the characteristic exponents.

\begin{figure}[h]
    \centering
\begin{tikzpicture}[line cap=round,line join=round,>=triangle 45,x=0.55cm,y=0.55cm]
\clip(-0.3,-1) rectangle (23,2);
\draw [->,line width=1pt] (0,0) -- (22,0);

\draw [shift={(6.05,-2.28)},line width=1pt]  plot[domain=1.01:2.17,variable=\t]({3.17*cos(\t r)},{3.17*sin(\t r)});
\draw [shift={(10.87,-4.95)},line width=1pt]  plot[domain=1.37:2.01,variable=\t]({5.89*cos(\t r)},{5.89*sin(\t r)});
\draw [shift={(21.24,-5.53)},line width=1pt]  plot[domain=1.46:2.05,variable=\t]({6.61*cos(\t r)},{6.61*sin(\t r)});
\draw (12.3,1.15) node[anchor=north west] {$\dots$};
\draw (12.3,-0.45) node[anchor=north west] {$\dots$};
\begin{scriptsize}
\draw[fill=black] (0,0) circle (2pt);
\draw[color=black] (0,-0.7) node {$0$};
\draw[fill=blue] (4,0) circle (2.5pt);
\draw[color=blue] (4,-0.7) node {$\beta_0$};
\draw[fill=blue] (8,0) circle (2.5pt);
\draw[color=blue] (8,-0.7) node {$\beta_1$};
\draw[fill=blue] (18,0) circle (2.5pt);
\draw[color=blue] (18,-0.7) node {$\beta_g$};
\draw[color=black] (6.18,1.5) node {$e_0\Z$};
\draw[color=black] (10.72,1.5) node {$e_1\Z$};
\draw[color=black] (20.50,1.5) node {$\Z$};
\end{scriptsize}
\end{tikzpicture}
\caption{Illustration of a saturated set of exponents}
\end{figure}
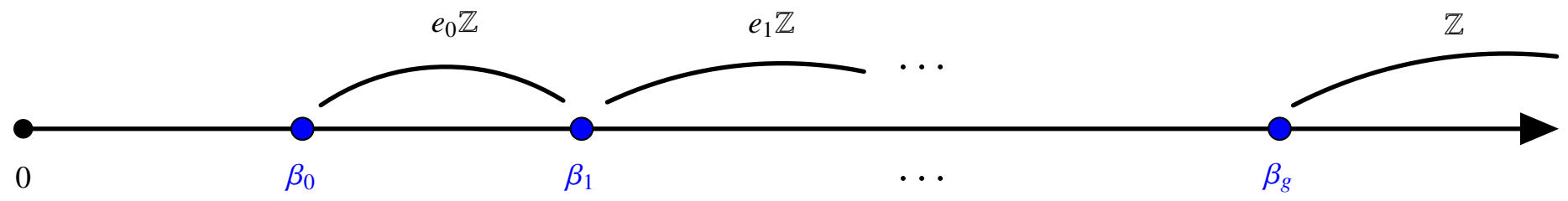

By \cite[Theorem VI.1.6]{BrianconGallicoGranger1980DeformationsEqu}, the
Lipschitz saturation of $A$ is the monomial algebra
\[
A^s=\C\{t^p\mid p\in E(A)^s\}.
\]

\begin{rmq}
If $X$ is a plane branch, the characteristic sequence defined above
coincides with the usual sequence of characteristic exponents of $X$;
see \cite[Lemma VI.3.3]{BrianconGallicoGranger1980DeformationsEqu}.
More generally, for a space branch, this sequence is the one which appears
as the sequence of characteristic exponents of a generic plane projection.
Thus the numbers $\beta_i$ make sense intrinsically, without choosing a
plane projection.
\end{rmq}

Since $(X^s)^\alpha=X^\alpha$, we can compute the Hölder saturation from
the monomial algebra $A^s$. Let
\[
E^\alpha(X,0):=
\{p\in\N\mid t^p\in\Ocal_{X,0}^\alpha\}.
\]
The following proposition shows that the set $E^\alpha(X,0)$ is obtained by translating the groups $e_i\Z$ along the axis.

\begin{prop}\label{Prop_ComputeSemigroupOfTheHolderSat}
Let $X$ be a germ of irreducible complex curve. Let $p\in\N^*$ and $\alpha\in(0,1]$. Then
\[
p\in E^\alpha(X,0)
\iff
\exists i\in\{0,\ldots,g\}\text{ such that }
p\geq\alpha\beta_i
\text{ and }
e_i\mid p.
\]
\end{prop}

\begin{proof}
    Let $\psi = (\varphi,\varphi') : (\C,0) \to (X\times X,0)$ be an arc. We consider the following numbers :
    \begin{enumerate}
        \item[] $v_0 = \inf\{ \ord\varphi, \ord\varphi'\}$
        \item[] $e = \inf\{ p\in \N^* \mid  \ord(\varphi^p - (\varphi')^p) > pv_0\}$ and $e=0$ if the inf does not exists.
        \item[] $v_1 = \ord(\varphi - \xi\varphi')$ where $\xi$ is the $e^{\text{th}}$ root of unity such that $v_1>v_0$.
    \end{enumerate}
    Then, by \cite[Lemma VI.1.1]{BrianconGallicoGranger1980DeformationsEqu}, we have, for every $p\in \N^*$,
    \[\begin{array}{l}
         \ord(\varphi^p-(\varphi')^p) = pv_0 \text{ if }e\nmid p \\
         \ord(\varphi^p-(\varphi')^p) = (p-1)v_0+v_1 \text{ if }e\mid p  
    \end{array}\]
    Now consider $i_e = \inf\{ k \text{ such that } e\nmid \beta_k\}$. By the proof of \cite[Theorem VI.1.6]{BrianconGallicoGranger1980DeformationsEqu}, we have $\ord(\psi^*(I_\Delta)) = \inf\big\{\beta_{i_e}v_0; (\beta_0 - 1)v_0 +v_1\big\}$.

    Let $p\in \N^*$ be such that $p\geq \alpha\beta_i$ and $e_i \mid p$ for some $i$. If $\beta_i \geq \beta_{i_e}$, then $\ord(\varphi^p - (\varphi')^p) \geq \alpha \beta_i v_0 \geq \alpha\beta_{i_e}v_0\geq \alpha\ord(\psi^*(I_\Delta))$.
    If $\beta_i < \beta_{i_e}$, then $e\mid p$ because $e$ divides $e_i = \gcd(\beta_0,\dots,\beta_{i})$ by definition of $i_e$, and $e_i \mid p$ by assumption. So, by the Lemma, we have
    \[
    \begin{array}{lll}
        \ord(\varphi^p - (\varphi')^p) & = (p-1)v_0+v_1&  \\
         & \geq (\alpha\beta_i-1)v_0 + v_1& \text{ by assumption }\\
         & \geq (\alpha\beta_0-1)v_0 + v_1& \text{ because }\beta_i \geq\beta_0\\
         & \geq \alpha(\beta_0-1)v_0 + \alpha v_1& \text{ because }\alpha \leq 1\\
         & \geq \alpha \ord(\psi^*(I_\Delta))
    \end{array}
    \]
    Hence, for every arc $\psi$, we have $\ord(\varphi^p - (\varphi')^p) \geq \alpha \ord(\psi^*(I_\Delta))$. Which means that $p\in E^\alpha(X,0)$.

    Conversely, let $p\in\N^*$ be such that, for every $i$, we have $p<\alpha\beta_i$ or $e_i\nmid p$. Let us consider $i$ such that $\alpha\beta_{i-1} \leq p< \alpha\beta_i$ with $\beta_{-1} = 0$. Assume $i\neq 0$ and consider the arc $\psi = (\varphi, \varphi')$ with $\varphi(t) = t^{v_0}$ and $\varphi'(t) = \xi t^{v_0} + t^{v_1}$ where $\xi$ is a primitive $e_{i-1}^{th}$-root of unity. So, for this arc, we have $e=e_{i-1}$ and $\beta_{i_e} = \beta_i$, so 
    \[
    \ord(\psi^*(I_\Delta)) = \inf\big\{\beta_{i}v_0; (\beta_0 - 1)v_0 +v_1\big\}.
    \]
    If $i\neq 0$, we moreover have $e=e_{i-1} \nmid p$ by assumption, because $\alpha\beta_{i-1}\leq p$. Then, $\ord(\varphi^p - (\varphi')^p) = pv_0< \alpha\beta_iv_0 < \alpha(\beta_0 - 1)v_0+\alpha v_1$ for $v_1$ large enough. If $i=0$, then consider $\varphi(t) = t^{v_0}$ and $\varphi'(t) = 2t^{v_0}$. We get $\ord(\varphi^p - (\varphi')^p) = pv_0 < \alpha\beta_0v_0 = \alpha \ord(\psi^*(I_\Delta))$.
\end{proof}

\begin{ex}\label{Ex_Cusp}
    Let $X$ be the cuspidal curve given by $y^2-x^3=0$. Its characteristic
sequence is $\beta_0=2$, $\beta_1=3$, with $e_0=2$ and $e_1=1$. One may
think that the $1/2$-Hölder saturation is the normalization of $X$ because
$1\geq\alpha\beta_0$. However, $e_0$ does not divide $1$. In fact,
$X^\alpha=X$ if $\alpha>1/3$ and $X^\alpha=\overline X$ if
$\alpha\leq1/3$.
\end{ex}

The previous example illustrates the more general fact that one can determine, from the characteristic sequence, when the $\alpha$-Hölder saturation coincides with the normalization.

\begin{cor}\label{Cor_HolderEqNormaIffBeta_g}
    Let $(X,0)$ be an irreducible complex curve germ and let $\beta_g$ be the greatest element of its characteristic sequence. Then
    \[
    (X^\alpha,0)=(\overline X,0)
    \iff
    \alpha\leq\frac{1}{\beta_g}.
    \]
\end{cor}

\begin{proof}
    We have $e_i = 1$ if and only if $i=g$, so the condition $e_i\mid 1$ can only be verified by $e_g$. Moreover $1\geq \alpha\beta_g$ if and only if $\alpha \leq \frac{1}{\beta_g}$, so $1$ belongs to $E^\alpha(X,0)$ if and only if $\alpha \leq \frac{1}{\beta_g}$.
\end{proof}

From this result, one can deduce a sufficient condition, depending on the $\delta$-invariant $\delta(X,0) := \dim_\C(\overline{\Ocal_{X,0}}/\Ocal_{X,0})$ and the multiplicity $m=m(X,0)$ of the germ, for the Hölder saturation to coincide with the normalization.

\begin{cor}\label{Cor_HolderSatDeltaBound}
    Let $(X,0)$ be an irreducible complex curve germ of multiplicity $m$. Then
    \[
    \alpha\leq\frac{1}{2\delta(X,0)-m+3} \implies (X^\alpha,0) = (\overline{X},0).
    \]
\end{cor}

\begin{proof}
    Let $\{\beta_0,\ldots,\beta_g\}$ be the characteristic sequence of $(X,0)$. Let $1\leq p<\beta_g$ and suppose that $p\in E(A)^s$. By definition of the
    saturated set, there exists $i\leq g-1$ such that $e_i\mid p$. Since
    $e_{g-1}\mid e_i$, it follows that $e_{g-1}\mid p$. Hence every integer $1\leq p < m$ is a gap of $E(A)^s$ and every integer $m<p<\beta_g$ which is not divisible by $e_{g-1}$ is a gap of $E(A)^s$. Since $e_{g-1}\geq 2$ by definition, we obtain
    \[
    \delta(X^s,0)
    \geq (\beta_g-1)
    -\left\lceil\frac{\beta_g-m}{e_{g-1}}\right\rceil
    \geq (\beta_g-1)
    -\left\lceil\frac{\beta_g-m}{2}\right\rceil \geq \frac{\beta_g+m-3}{2}.
    \]
    Thus $\beta_g\leq 2\delta(X^s,0)-m+3$, and since $\Ocal_{X,0}\subseteq\Ocal_{X^s,0}\subseteq\Ocal_{\overline X,0}$ implies
    $\delta(X^s,0)\leq\delta(X,0)$, we get the desired result from Corollary \ref{Cor_HolderEqNormaIffBeta_g}.
\end{proof}

We can also obtain a lower bound on $\alpha$ for which the $\alpha$-Hölder saturation will coincide with the Lipschitz saturation.

\begin{cor}\label{Cor_HolderSatEqualsLipschitzIrreducible}
Let $(X,0)$ be an irreducible complex curve germ and let $\beta_g$ be the greatest element of its characteristic sequence. Then
$$
\alpha>1-\frac1{\beta_g}
\quad\Longrightarrow\quad
(X^\alpha,0)=(X^s,0).
$$
\end{cor}

\begin{proof}
Let $\{\beta_0,\ldots,\beta_g\}$ be the characteristic sequence of $(X,0)$ and put $e_i=\gcd(\beta_0,\ldots,\beta_i)$. Let $p\in \N^*$ such that there exists $e_i$ satisfying $e_i\mid p$. Then $p\in E(A)^s$ if and only if $p\geq \beta_i$ and $p\in E(A)^\alpha$ if and only if $p\geq \alpha\beta_i$. If $\alpha>1-1/\beta_g$, then, for every $i$, we have 
$$
\alpha\beta_i>
\beta_i-\frac{\beta_i}{\beta_g}
\geq\beta_i-1.
$$
Since $p$ is an integer and $\beta_i \geq \alpha\beta_i$, we get $p\geq \alpha\beta_i > \beta_i-1$ if and only if $p\geq \beta_i$. Hence $E(A)^s = E(A)^\alpha$, which proves the desired result.
\end{proof}

From the proof of Corollary \ref{Cor_HolderSatDeltaBound}, we know that $\beta_g\leq 2\delta(X,0)-m+3$, so we can deduce directly the following additional bound. Recall that the Milnor number of the branch $(X,0)$ is $\mu(X,0)=2\delta(X,0)$.

\begin{cor}\label{Cor_BoundAlphaIsLipIrreduciblecase}
    Let $(X,0)$ be an irreducible complex curve germ of multiplicity $m$, with $\delta=\delta(X,0)$ and $\mu=\mu(X,0)$. Then
    \[
    \alpha> 1 - \frac{1}{2\delta-m+3} = \frac{\mu-m+2}{\mu-m+3} \implies (X^\alpha,0) = (X^s,0).
    \]
\end{cor}

Another consequence of Proposition \ref{Prop_ComputeSemigroupOfTheHolderSat} is that we can determine, from the characteristic sequence, the $\alpha$'s for which the $\alpha$-Hölder saturation has the same multiplicity as the original curve.

\begin{prop}\label{Prop_HolderSatWithSameMultiplicity}
    Let $(X,0)$ be an irreducible complex curve germ. Then $(X,0)$ and $(X^\alpha,0)$ have the same multiplicity if and only if 
    \[
    \alpha > \max_{i=1,\dots,g}\Big\{\frac{\beta_0-e_i}{\beta_i}\Big\}
    \]
\end{prop}

\begin{proof}
    Since, for every $i\in \{1,\dots,g\}$, we have $e_i = \gcd(\beta_0,\dots,\beta_i)$, then $\beta_0-e_i$ is the greatest element, strictly smaller than $\beta_0$, which is divisible by $e_i$. So, if there exists $i$ such that $\alpha\beta_i\leq \beta_0-e_i$, then $\beta_0-e_i\in E^\alpha(X,0)$ and, since $\beta_0-e_i<\beta_0$, this would mean that $(X,0)$ and $(X^\alpha,0)$ do not have the same multiplicity.
    
    Conversely, suppose that $\alpha \beta_i > \beta_0-e_i$ for all $i$ and let $p<\beta_0$. If there exists $i$ such that $e_i\mid p$, then we have $p\leq \beta_0 - e_i < \alpha \beta_i$ and so $p\notin E^\alpha(X,0)$. Hence $(X,0)$ and $(X^\alpha,0)$ have the same multiplicity. 
\end{proof}

\begin{rmq}
    One might think it would suffice to check the condition only for the first $\beta_i$'s, since they are closer to $\beta_0$, but this is not the case. Consider, for instance, the characteristic sequence $\beta_0=8$, $\beta_1=12$ and $\beta_2=13$. Then $e_0=8$, $e_1=4$ and $e_2=1$. For $\alpha=1/2$, we have $\alpha\beta_1>\beta_0-e_1$ but $\alpha\beta_2\leq\beta_0-e_2$.
\end{rmq}

As a consequence of Proposition \ref{Prop_HolderSatWithSameMultiplicity}, one obtains a sufficient condition depending only on the multiplicity of the curve.

\begin{cor}\label{Cor_HolderSatPreservesMultiplicity}
    Let $(X,0)$ be an irreducible complex curve germ of multiplicity $m = m(X,0)$. Then
    \[
    \alpha>\frac{m-1}{m+1}
    \quad\Longrightarrow\quad
    m(X^\alpha,0)=m(X,0).
    \]
\end{cor}

\begin{proof}
    For every $i\in\{1,\ldots,g\}$, both $m=\beta_0$ and $\beta_i$
    are divisible by $e_i$. Since $\beta_i>m$, it follows that $\beta_i\geq m+e_i$. Therefore
    \[
    \frac{m-e_i}{\beta_i}
    \leq
    \frac{m-e_i}{m+e_i}
    \leq
    \frac{m-1}{m+1},
    \]
    since $e_i\geq1$. The result follows from
    Proposition~\ref{Prop_HolderSatWithSameMultiplicity}.
\end{proof}

\begin{rmq}
    Note that one can not obtain a better uniform sufficient condition since the previous condition is an equivalence for curve germs with characteristic sequence of the form $\{m,m+1\}$.
\end{rmq}

\begin{rmq}
    Note that, even if this bound is quite close to $1$ in general, it is not a trivial bound in the sense that the $\alpha$-Hölder saturation for such an $\alpha$ is not the Lipschitz saturation in general. For example, consider a germ of curve $(X,0)$ such that $\Ocal_{X,0} = \C\{t^2,t^4,t^5\}$. It is Lipschitz saturated, the bound says that multiplicity is preserves by $\alpha$-Hölder saturation for $\alpha > 1/3$ and yet $\Ocal^\alpha_{X,0}$ is different from $\Ocal_{X,0}^s$ for every $\alpha\in (0,3/5]$.
\end{rmq}

\subsection{Hölder saturation for reducible complex curves}

In this subsection we extend the computation of the preceding one to reduced curve germs with several branches. A function on the normalization is $\alpha$-Hölder on $X$ if and only if it is $\alpha$-Hölder on each branch and the differences of its values on two branches are controlled by the distance between the branches. The first condition is the one computed in Subsection~\ref{Subsec_IrreducibleCurves}. The second one is governed by the contact between the branches, which we measure by the orders $m_{i,j,\xi}$ of the diagonal ideal along the arcs joining the sheets of two branches. Theorem~\ref{Theo_HolderSatReducibleCurve} describes the saturation in these terms. We begin with the valuative criterion, which holds for any reduced complex analytic germ, and then specialize to curves.

Let $X$ be a germ of reduced complex analytic variety and let
$X=\cup_{i=1}^n X_i$ be its decomposition into irreducible components.
Then, for any function $f:X\to\C$, the following two properties are
equivalent:
\begin{enumerate}
    \item for every $(x,y)\in X\times X$, we have
    $|f(x)-f(y)|\leq C\|x-y\|^\alpha$;
    \item for every $i,j\in\{1,\ldots,n\}$ and every
    $(x_i,y_j)\in X_i\times X_j$, we have
    \[
    |f_{|X_i}(x_i)-f_{|X_j}(y_j)|
    \leq C\|x_i-y_j\|^\alpha.
    \]
\end{enumerate}

Moreover, if $\overline X$ denotes the normalization of $X$, then $\overline X=\overline X_1\sqcup\cdots\sqcup\overline X_n$ and
\[
\Ocal_{\overline X}
=
\bigoplus_{i=1}^n \Ocal_{\overline X_i}.
\]
Indeed, for a reduced ring, the total ring of fractions is the product of
the fraction fields of its irreducible components, and its integral closure
in the total ring of fractions is the product of the integral closures of
the irreducible components. Let $f$ be a locally bounded meromorphic function on $X$. Then its pull-back to the normalization can be written $f\circ\pi=(f_1,\ldots,f_n)$ where $f_i=f_{|X_i}\circ\pi_i$. For $i,j\in\{1,\ldots,n\}$, let
\[
I_\Delta^{(i,j)} := \big\langle p_i\otimes1-1\otimes p_j \mid p=(p_1,\ldots,p_n)\in\Ocal_X \big\rangle
\]
be an ideal of $\Ocal_{\overline X_i} \widehat\otimes_{\C} \Ocal_{\overline X_j} = \Ocal_{\overline X_i\times\overline X_j}$. By applying to $I_\Delta^{(i,j)}$ the interpretation of rational powers of
ideals in terms of inequalities, see \cite[Theorem 7.2]{LejeuneTeissier2008ClotureInt}, we get that $f$ is $\alpha$-Hölder if
and only if $f_i\otimes1-1\otimes f_j \in \big(I_\Delta^{(i,j)}\big)_\alpha$ for every $i,j\in\{1,\ldots,n\}$. By the valuative criterion, this is equivalent to asking, for every $i,j$, the condition
\[
(P_{i,j}) :
\qquad
\ord_u\big(f_i\circ\varphi_i(u)-f_j\circ\varphi_j(u)\big)
\geq
\alpha\,
\ord_u\big(\varphi^*(I_\Delta^{(i,j)})\big)
\]
for every arc $\varphi=(\varphi_i,\varphi_j): (\C,0)\longrightarrow (\overline X_i\times\overline X_j,0)$. Hence, if we define $m_\varphi(f) := \ord_u\big(f_i\circ\varphi_i(u)-f_j\circ\varphi_j(u)\big)$ and $m_\varphi:=\inf_{g\in\Ocal_X}m_\varphi(g)=\ord_u\big(\varphi^*(I_\Delta^{(i,j)})\big),$
the condition becomes
\[
(P_{i,j}) :
\qquad
m_\varphi(f)\geq\alpha m_\varphi
\]
for every arc $\varphi$.

Suppose now that $X$ is a germ of reduced complex analytic curve. Then $\overline{\Ocal_{X}}=\C\{t_1\}\oplus\cdots\oplus\C\{t_n\}$.
For every $i$, let $m_i\in\N^*$ be the multiplicity of $X_i$, let
$\eta=\operatorname{lcm}(m_1,\ldots,m_n)$ and let
$\eta_i=\eta/m_i$. We say that $(t_1,\ldots,t_n)$ is a family of
\textit{good uniformizers} if there exists a transverse parameter
$x\in\Ocal_{X}$ whose image in the normalization is $x=(t_1^{m_1},\ldots,t_n^{m_n})$. For such a family, we have $t_i^{m_i}\otimes1-1\otimes t_j^{m_j}
\in I_\Delta^{(i,j)}$. Let $\xi$ be a $\eta$-th root of unity and consider the arc
\[
\begin{array}{cccc}
\varphi_{i,j,\xi}:&
(\C,0)&\longrightarrow&
(\overline X_i\times\overline X_j,0)\\
&u&\longmapsto&
(u^{\eta_i},(\xi u)^{\eta_j}).
\end{array}
\]
For an element $f=(f_1,\ldots,f_n)\in\overline{\Ocal_{X}}$, define $m_{i,j,\xi}(f):=m_{\varphi_{i,j,\xi}}(f)$ and $m_{i,j,\xi}:=m_{\varphi_{i,j,\xi}}.$

\begin{theorem}\label{Theo_HolderSatReducibleCurve}
    Let $\Ocal_X\subset \overline{\Ocal_X} = \C\{t_1\} \oplus \dots \oplus \C\{t_n\}$ be a reduced analytic algebra of dimension one with a family of good uniformizers $(t_1,\dots,t_n)$. Then, an element $f=(f_1,\dots,f_n)\in \overline{\Ocal_X}$ is an element of $\Ocal_X^\alpha$ if and only if
    \[\begin{array}{rl}
        (1) & f_i\in \Ocal_{X_i}^\alpha\text{, for all }i\in \{1,\dots,n\} \medskip\\
        (2) & m_{i,j,\xi}(f) \geq \alpha m_{i,j,\xi}\text{, for all }i\neq j\in \{1,\dots,n\}^2\text{ and for all }\eta\text{-root of unity }\xi. 
    \end{array}\]
\end{theorem}

\begin{proof}
    Let us recall that an element $f = (f_1,\dots,f_n)\in \overline{\Ocal_{X}}$ is in $\Ocal_{X}^\alpha$, if and only if, it verifies $(P_{i,j})$ for all $i,j\in \{1,\dots,n\}$. Then $(1)$ is equivalent to having $(P_{i,i})$ for all $i\in\{1,\dots,n\}$, and $(2)$ is clearly a necessary condition, since it is the same as $(P_{i,j})$ but only for the arcs of the form $\varphi_{i,j,\xi}$.
    
    Conversely, let us consider an element $f = (f_1,\dots,f_n)\in \overline{\Ocal_{X}}$ which verifies condition $(1)$ and $(2)$. Let $\varphi : (\C,0)\to \overline{X_i}\times \overline{X_j}$ be an arc and suppose that $\ord(\varphi_i^{m_i} - \varphi_j^{m_j}) = \min\big\{\ord(\varphi_i^{m_i}),\ord(\varphi_j^{m_j})\big\} = \ord(\varphi^*(I_\Delta^{(i,j)}))$. Then, consider the following inequality
    \[
    \begin{array}{lll}
        \ord\big(f_i\circ\varphi_i -f_i(0) + f_i(0)-f_j(0)+f_j(0)-f_j\circ\varphi_j\big) \geq \min\big\{ & \ord(f_i\circ\varphi_i -f_i(0)),&\\
        & \ord(f_i(0)-f_j(0)),&\\
        &\ord(f_j(0)-f_j\circ\varphi_j)&\big\}
    \end{array}
    \]
    By $(1)$, we have $f_i \in \Ocal_{X_i}^\alpha$. So $\ord_t(f_i\circ\varphi_i(t) - f_i(0)) \geq \alpha \ord(\varphi_i^{m_i} - 0) = \alpha m_i\ord(\varphi_i)$. Idem for $f_j$. Moreover, condition $(2)$ implies that $f_i(0) = f_j(0)$. Hence, we get 
    \[
    \ord(f_i\circ\varphi_i - f_j\circ\varphi_j) \geq \min\big\{ \alpha m_i\ord(\varphi_i), \alpha m_j\ord(\varphi_j)\big\} = \alpha \ord(\varphi^*(I_\Delta))
    \]
    which is exactly condition $(P_{i,j})$.
    
    Now, let us consider the case where $\ord(\varphi_i^{m_i} - \varphi_j^{m_j}) > \min\big\{\ord(\varphi_i^{m_i}),\ord(\varphi_j^{m_j})\big\}$. Following the notations of the proof of \cite[Theorem VI.2.2]{BrianconGallicoGranger1980DeformationsEqu}, we write $\varphi_i(u) = u^{v_i}E_i(u)$ with $E_i(u)\in \C\{u\}^{\times}$, idem for $j$. Then, we consider the arc
    \[
    \begin{array}{cccc}
        \mathcal{T} :& (\C,0) &\to &\overline{X_i}\times \overline{X_j}  \\
        & u & \mapsto &(\tau(u), \tau'(u)) 
    \end{array}
    \]
    such that $\tau(u)=u^{m_iv_i}E(u^\eta)^{1/\eta_i}$ and $\tau'(u)=u^{m_jv_j}E(u^\eta)^{1/\eta_j}$. Then we have $\tau(u)^{\eta_i} = u^\eta E(u^\eta) = \varphi_i(u^\eta)$, idem for $\tau'$. Therefore, we have 
    \[
    \begin{array}{lll}
	   \eta \ord\big(f_i\circ\varphi_i(u)-f_j\circ\varphi_j(u)\big) & = & \ord\big(f_i\circ\varphi_i(u^\eta)-f_j\circ\varphi_j(u^\eta)\big)\medskip\\
	   & = & \ord\big(f_i(\tau(u)^{\eta_i})-f_j(\tau'(u)^{\eta_j})\big)
    \end{array}
    \]
    and 
    \[
    \begin{array}{lcll}
    &\ord(\varphi_i^{m_i} - \varphi_j^{m_j}) &>& \min\big\{\ord(\varphi_i^{m_i}),\ord(\varphi_j^{m_j})\big\}\medskip\\
    \iff & \ord(\varphi_i(u^\eta)^{m_i} - \varphi_j(u^\eta)^{m_j}) &>& \min\big\{\ord(\varphi_i(u^\eta)^{m_i}),\ord(\varphi_j(u^\eta)^{m_j})\big\}\medskip\\
    \iff & \ord(\tau(u)^{\eta} - \tau'(u)^{\eta}) &>& \min\big\{\ord(\tau(u)^{\eta},\ord(\tau'(u)^{\eta}\big\} = \eta v_0\\
    \end{array}
    \]
    where $v_0 = \ord(\tau(u)) = \ord(\tau'(u))$. By \cite[Lemma VI.1.1]{BrianconGallicoGranger1980DeformationsEqu}, we can consider a $\eta$-root of unity $\xi$ such that $\ord(\tau'-\xi\tau)>v_0$. Now, let us consider the expression $f_i(x^{\eta_i}) - f_j(y^{\eta_j})$ in $\C\{x,y\}$. Then, by Weierstrass division theorem, there exists $Q\in \C\{x,y\}$ and $R\in\C\{x\}[y]$ such that 
    \[
    \hspace{ -0.5cm}(*)\text{ \hspace{0.5cm}}f_i(x^{\eta_i}) - f_j(y^{\eta_j}) = (y-\xi x)Q(x,y) + R(x,y) \text{ with }\deg_y(R) < 1.
    \]
    So $R(x,y)$ does not depend on $y$ and, in particular, we have $R(x,y) = R(x,\xi x) = f_i(x^{\eta_i}) - f_j((\xi x)^{\eta_j}) = U(x)x^{m_{i,j,\xi}(f)}$ with $U\in \C\{x\}^{\times}$. Then $(y-\xi x)Q(x,y) = f_j((\xi x)^{\eta_j}) - f_j(y^{\eta_j})$ and, since, by $(1)$, we can apply $(P_{j,j})$ to $f_j$ with the arc $u\mapsto ((\xi \tau(u))^{\eta_j}, \tau'(u)^{\eta_j})$, we have
    \[
    \ord\big(f_j((\xi \tau)^{\eta_j}) - f_j((\tau')^{\eta_j})\big) \geq \alpha \ord\big( ((\xi \tau)^{\eta_j})^{m_j}, ((\tau')^{\eta_j})^{m_j} \big) = \alpha\ord\big(\tau^\eta - \tau'^\eta\big).
    \]
    So, by applying the morphism $\mathcal{T}^*$ to the expression $(*)$, we get that $\ord\big(f_i(\tau)^{\eta_i} - f_j((\tau')^{\eta_j})\big) \geq \min\big\{ m_{i,j,\xi}(f)v_0, \alpha\ord\big(\tau(u)^\eta - \tau'(u)^\eta \big) \big\}$. Moreover, by \cite[Lemma VI.2.3]{BrianconGallicoGranger1980DeformationsEqu}, we have $\eta\ord\big(\varphi^*(I_\Delta^{(i,j)})\big) = \min\big\{ m_{i,j,\xi} v_0, \ord\big(\tau(u)^\eta - \tau'(u)^\eta \big)\big\}$. Then 
    \[
    \begin{array}{llll}
	   \eta \ord\big(f_i\circ\varphi_i - f_j\circ\varphi_j\big) & \geq & \min\big\{ m_{i,j,\xi}(f) v_0, \alpha\ord\big(\tau(u)^\eta - \tau'(u)^\eta \big)\big\}&\medskip\\
	   & \geq & \min\big\{ \alpha m_{i,j,\xi} v_0, \alpha\ord\big(\tau(u)^\eta - \tau'(u)^\eta \big)\big\}&\text{ by }(2)\medskip\\
	   & \geq & \alpha \min\big\{ m_{i,j,\xi} v_0, \ord\big(\tau(u)^\eta - \tau'(u)^\eta \big)\big\} = \alpha \eta\ord\big(\varphi^*(I_\Delta^{(i,j)})\big).&
    \end{array}
    \]
    This mean that $f$ verifies $(P_{i,j})$.
\end{proof}

\begin{rmq}
    The case of a germ of curve with smooth branches is particularly simple since condition $(1)$ of Theorem \ref{Theo_HolderSatReducibleCurve} is then empty and $\eta_i = \eta = 1$ for all $i$, so $m_{i,j,\xi}(f) = \ord(f_i(u)-f_2(u))$. 
\end{rmq}

\begin{ex}\label{Ex_Parusinski}
    Let $X_1 := \big\{(y-x^2)y = 0\big\} \subset \C^2$ and let $X_2 := \big\{ z=y-x^2=0 \big\} \cup \big\{ y=x=0 \big\} \subset \C^3$. The normalization of $(X_1,0)$ is given by 
    $$\begin{array}{ccc}
         (\C,0)\cup(\C,0)& \to & (X_1,0)  \\
         t_1 & \mapsto & (t_1,t_1^2)\\
         t_2 & \mapsto & (t_2,0)
    \end{array}.$$
    Moreover the branches of $(X_1,0)$ have both multiplicity one so, in criterion $(2)$ of Theorem \ref{Theo_HolderSatReducibleCurve}, we only need to consider 
    $$m_{1,2,1} = \inf_{g\in \Ocal_{X_1}}\big\{ \ord_u( g_1\circ \varphi_{1,2,1}(u) - g_2\circ \varphi_{1,2,1}(u)) \big\} = \inf\big\{ \ord_u(u-u), \ord_u(u^2-0) \big\} =2$$
    with $\varphi_{1,2,1}(u) = (u,u)$. Since the branches are smooth, the conditions of Theorem \ref{Theo_HolderSatReducibleCurve} simply becomes 
    $$f = (f_1,f_2) \in \Ocal_{X_1}^\alpha \iff \ord_u(f_1(u)-f_2(u))\geq 2\alpha$$
    
    Now, to compute the Hölder saturation of $(X_2,0)$, one can for example consider the orthogonal projection $\pi$ to the plane $\{z-y=0\}$ which is given by $\pi :(\C^3,0)\to (\C^2,0) , (x,y,z)\mapsto (x,y/2+z/2)$. Since $\ker(\pi)$ is generated by $(0,1,-1)$, one can verify that it is a generic linear projection in the sense of Definition \ref{Def_GenericProjection}. So, by Corollary \ref{Cor_HolderStableUnderProjection}, the Hölder saturation of $X_2$ is isomorphic to the one of $\pi(X_2)$. Under this projection, the first branch becomes $\{2y-x^2=0\}$ and the second one becomes $\{x=0\}$. Then, the normalization of $(\pi(X_2),0)$ is given by 
    $$\begin{array}{ccc}
         (\C,0)\cup(\C,0)& \to & (X_2,0)  \\
         t_1 & \mapsto & (t_1,\frac{1}{2}t_1^2)\\
         t_2 & \mapsto & (0,t_2)
    \end{array}.$$
    So $m_{1,2,1} = \inf\big\{ \ord_u(u-0), \ord_u(\frac{1}{2}u^2-u) \big\} =1$. Then the conditions of Theorem \ref{Theo_HolderSatReducibleCurve} becomes
    $$f\in \Ocal_{X_2}^\alpha \iff \ord_u(f_1(u)-f_2(u))\geq \alpha$$
    
    If we consider for example the element $(t_1^2,t_2)\in \overline{\Ocal_{X_1}} = \overline{\Ocal_{X_2}}$. We have $\ord_u((t_1^2 - t_2)\circ \varphi_{1,2,1}(u) ) = \ord_u(u^2 - u)=1$, so this element belongs to $\Ocal_{X_2}^s$ but not to $\Ocal_{X_1}^s$.
\end{ex}

In Example \ref{Ex_Parusinski}, the key point is that $m_{1,2,1}$ is equal to $2$ for $X_1$ and $1$ for $X_2$, so it seems that the curves should be $\alpha$-Hölder equivalent if $\alpha\leq 1/2$. From the next corollary, which deals with the particular case of curves with smooth branches, we can see that it is indeed the case.

\begin{cor}
    Let $X$ be a germ of curve with smooth branches $X_i$ crossing at $0$. In this situation, we have $\varphi_{i,j,\xi}(u) = (u,u)$ for every $i,j$. For $f = (f_1,\dots,f_n) \in \overline{\Ocal_{X}}$, we write $m_{i,j}(f) := m_{i,j,\xi}(f) = \ord_u\{ f_i(u) - f_j(u) \}$ and $m_{i,j} := m_{i,j,\xi} = \inf_{f\in \Ocal_{X}}\{ m_{i,j}(f) \}$. Then 
    \[
    X^\alpha = X^+ \iff \alpha \leq \inf_{i\neq j}\big\{1/m_{i,j}\big\}
    \]
\end{cor}

\begin{proof}
    Let $f = (f_1,\dots,f_n) \in \overline{\Ocal_{X}}$. We have $f\in \Ocal_{X}^\alpha$ if and only if it satisfies conditions $(1)$ and $(2)$ of Theorem \ref{Theo_HolderSatReducibleCurve}. However, condition $(1)$ is always true since, for every branch, we have $\Ocal_{X_i}^\alpha \simeq \overline{\Ocal_{X_i}} \simeq \C\{t_i\}$.
    
    Assume $\alpha \leq \inf_{i\neq j}\big\{1/m_{i,j}\big\}$. We already have $\Ocal_{X}^\alpha \subseteq \Ocal^+_{X}$, so let us consider $f\in \Ocal^+_{X}$ in order to show the opposite inclusion. For all $i\neq j$, we have
    $$ 1 \geq \inf_{k\neq l}\big\{1/m_{k,l}\big\} m_{i,j} \geq \alpha m_{i,j}.$$
    Then $f$ belongs to $\Ocal_{X}^\alpha$ if and only if $m_{i,j}(f)\neq 0$. Since $f$ is an element of the seminormalization, it has constant value on the fibers of the normalization morphism, so $f_i(0) = f_j(0)$. Hence $m_{i,j}(f)\geq 1$ and we get $f\in \Ocal_{X}^\alpha$.
    
    Conversely, suppose $X^\alpha\simeq X^+$ and assume that there exist $i_0,j_0$ such that $1/m_{i_0,j_0} < \alpha$. This is equivalent to $1< \alpha m_{i_0,j_0}$ so, if we construct an element $f\in \Ocal^+_{X}$ such that $m_{i_0,j_0}(f) =1$, we will obtain a contradiction. Let $f = (0,..,0,t_{i_0},0,..,0) \in \bigoplus_{k} \C\{t_k\}$. It is an element of $\Ocal_{X}^+$ since $f(0)= (0)_{1,\dots,n}$ and $m_{i_0,j_0}(f) = \ord_u\{f_{i_0}(u) - f_j(0)\} = \ord_u\{u - 0\} = 1$.
\end{proof}

\subsection{Bounds in the case of reduced curves}\label{Subsec_BoundsReducedCurves}


In this subsection, we extend the bounds obtained previously for irreducible curves to the reducible case. We start by establishing a bound on the contact integers $m_{i,j,\xi}$ considered in Theorem \ref{Theo_HolderSatReducibleCurve} involving the $\delta$-invariant of the curve and its number of branches.

\begin{lemma}\label{Lem_Bound_mijxi_Delta}
    Let $(X,0)$ be a reduced complex curve germ with $r$ irreducible
    components. Let $\delta=\delta(X,0)$ and let us recall the notation $\eta_i = \eta/m_i$ where $\eta$ is the lcm of the multiplicities $m_i$ of the branches of $(X,0)$. For $i\neq j$, let $a_{i,j}:=\min\{\eta_i,\eta_j\}$. Then, for every $\eta$-th root of unity $\xi$, we have
    \[
    m_{i,j,\xi} \leq a_{i,j}\big(\delta-r+2\big).
    \]
\end{lemma}

\begin{proof}
    Let $\ell\in\{i,j\}$ be such that $\eta_\ell=a_{i,j}$. For every
    $k\geq1$, consider the element
    \[
    h_k=(0,\ldots,0,t_\ell^k,0,\ldots,0)
    \in\Ocal_{X^+}.
    \]
    We then have $m_{i,j,\xi}(h_k)=a_{i,j}k$. So, in particular, if $a_{i,j}k<m_{i,j,\xi}$, then $h_k\notin\Ocal_{X^s}$ by Theorem
    \ref{Theo_HolderSatReducibleCurve} applied with $\alpha=1$.

    Moreover, the classes of the elements $h_k$ such that
    $a_{i,j}k<m_{i,j,\xi}$ are linearly independent in
    $\Ocal_{X^+}/\Ocal_{X^s}$. Indeed, for a nonzero linear combination
    $h=\sum_k c_kh_k$, if $k_0$ is the smallest integer such that
    $c_{k_0}\neq0$, then $m_{i,j,\xi}(h)=a_{i,j}k_0<m_{i,j,\xi}$, and therefore $h\notin\Ocal_{X^s}$. Hence, we can deduce the following inequality
    \[
    \left\lceil\frac{m_{i,j,\xi}}{a_{i,j}}\right\rceil-1
    \leq
    \dim_\C\big(\Ocal_{X^+}/\Ocal_{X^s}\big).
    \]
    Now, since $\delta(X^s) = \dim_{\C}\big(\Ocal_{\overline{X}}/\Ocal_{X^s}\big)$ and $\dim_\C\big(\overline{\Ocal_X}/\Ocal_{X^+}\big)=r-1$, we have
    \[
    \dim_\C\big(\Ocal_{X^+}/\Ocal_{X^s}\big)
    =
    \delta(X^s)-r+1
    \leq
    \delta(X)-r+1.
    \]
    So $\left\lceil\frac{m_{i,j,\xi}}{a_{i,j}}\right\rceil \leq \delta(X)-r+2$, which gives $m_{i,j,\xi} \leq a_{i,j}\big(\delta(X)-r+2\big)$.
\end{proof}

We can now extend Corollary \ref{Cor_HolderSatDeltaBound} to an
arbitrary reduced curve germ. Recall that, if $X$ has $r$ branches,
its Milnor number is $\mu(X)=2\delta(X)-r+1$.

\begin{cor}\label{Cor_HolderSatSeminormBoundReducedCurves}
    Let $(X,0)$ be a reduced complex curve germ of multiplicity
    $m=m(X,0)$. Then
    \[
    \alpha\leq\frac{1}{\mu(X,0)-m+3}
    \quad\Longrightarrow\quad
    (X^\alpha,0)=(X^+,0).
    \]
\end{cor}

\begin{proof}
    Write $X=\cup_{i=1}^r X_i$, and denote by $m_i$ and $\delta_i$ respectively the multiplicity and the $\delta$-invariant of each $X_i$. Let $\alpha \in (0,1]$ be such that $\alpha\leq 1/(\mu(X)-m+3)$ and let $f=(f_1,\ldots,f_r)\in\Ocal_{X^+}$. We want to prove that $f$ verifies conditions $(1)$ and $(2)$ of Theorem \ref{Theo_HolderSatReducibleCurve}.
    
    First, by \cite[Lemma 1.2.2]{BuchweitzGreuel1980ConstancyOfMilnorNbr}, we know that $\delta(X) \geq \sum_{j=1}^r\delta_j+r-1$. Moreover, for every $i$, we have $\delta_i\geq m_i-1$ and so $\delta(X)-\delta_i\geq m-m_i$. Since $m-m_i\geq r-1$, we obtain
    \[
    \mu(X)-m+3-\big(2\delta_i-m_i+3\big) = 2\big(\delta(X)-\delta_i\big) -(m-m_i)-r+1 \geq0.
    \]
    Hence
    \[
    \alpha
    \leq
    \frac{1}{\mu(X)-m+3}
    \leq
    \frac{1}{2\delta_i-m_i+3}.
    \]
    By Corollary \ref{Cor_HolderSatDeltaBound}, this implies that $(X_i)^\alpha=\overline{X_i}$, so $f_i\in\Ocal_{X_i}^\alpha$ for every $i$.

    Now, let $i\neq j$ and let $\xi$ be a $\eta$-th root of unity. Since
    $f\in\Ocal_{X^+}$, we have $f_i(0)=f_j(0)$, and therefore $m_{i,j,\xi}(f)\geq \min\{\eta_i,\eta_j\}=a_{i,j}$. On the other hand, Lemma \ref{Lem_Bound_mijxi_Delta} gives $m_{i,j,\xi} \leq a_{i,j}\big(\delta(X)-r+2\big)$ which implies $\alpha m_{i,j,\xi} \leq \alpha a_{i,j}\big(\delta(X)-r+2\big)$. Since
    \[
    \mu(X)-m+3-\big(\delta(X)-r+2\big) = \delta(X)-m+2,
    \]
    and $\delta(X)\geq m-1$, we have
    \[
    \alpha
    \leq
    \frac{1}{\mu(X)-m+3}
    \leq
    \frac{1}{\delta(X)-r+2}.
    \]
    This allows us to deduce $\alpha m_{i,j,\xi} \leq a_{i,j}$. Then we get $\alpha m_{i,j,\xi}\leq a_{i,j} \leq m_{i,j,\xi}(f)$, which is exactly condition $(2)$ of Theorem \ref{Theo_HolderSatReducibleCurve}. We have proven that $f\in\Ocal_X^\alpha$ by Theorem \ref{Theo_HolderSatReducibleCurve} and so $\Ocal_{X^+}\subseteq\Ocal_X^\alpha$. The reverse inclusion follows from Proposition \ref{Prop_InteractionsBetweenSaturations}, and therefore $X^\alpha=X^+$.
\end{proof}

We can also now extend Corollary \ref{Cor_BoundAlphaIsLipReduciblecase} to an arbitrary reduced curve germ.
\begin{cor}\label{Cor_BoundAlphaIsLipReduciblecase}
     Let $(X,0)$ be a reduced complex curve germ of multiplicity $m$, with $\mu=\mu(X,0)$. Then
     $$
    \alpha > \frac{\mu-m+2}{\mu-m+3} \implies (X^\alpha,0) = (X^s,0).
     $$
\end{cor}

\begin{proof}
Write $X=\bigcup_{i=1}^r X_i$, and denote by $m_i$ and $\delta_i$
the multiplicity and the $\delta$-invariant of $X_i$. Let $(\mu-m+2)/(\mu-m+3)<\alpha\leq1$ and let $f=(f_1,\ldots,f_r)\in\Ocal_{X^\alpha}$. As in the proof of Corollary \ref{Cor_HolderSatSeminormBoundReducedCurves}, we have $2\delta_i-m_i+3\leq\mu-m+3$. So, for each singular branch, Corollary \ref{Cor_BoundAlphaIsLipIrreduciblecase} gives $f_i\in\Ocal_{X_i^s}$.

It remains to verify condition $(2)$ of Theorem \ref{Theo_HolderSatReducibleCurve}. Using the notation of this theorem, fix $i\neq j$ and one of the arcs $\varphi_{i,j,\xi}(u)=(u^{\eta_i},(\xi u)^{\eta_j})$. By definition of $m_{i,j,\xi}$, for every $g=(g_1,\ldots,g_r)\in\Ocal_X$, the series $g_i(u^{\eta_i})$ and $g_j((\xi u)^{\eta_j})$ agree in every degree strictly smaller than $m_{i,j,\xi}$. Since the restriction $\Ocal_X\to\Ocal_{X_i}$ is surjective, it follows that, for every $p\in E(\Ocal_{X_i})$, we have 
$$\eta_i p<m_{i,j,\xi} \implies \eta_j\mid\eta_i p.$$

Let us show that the same implication holds for $E(\Ocal_{X_i})^s$. Let $p>0$ belong to this saturated set and assume that $\eta_i p<m_{i,j,\xi}$. Write $\{\beta_{i,0},\ldots,\beta_{i,g_i}\}$ for the characteristic sequence of $X_i$. Then there exists $\ell$ such that $\beta_{i,\ell}\leq p$ and $e_{i,\ell} := \gcd(\beta_{i,0},\ldots,\beta_{i,\ell})\mid p$. For every $h\leq\ell$, we have $\eta_i\beta_{i,h}\leq \eta_i p<m_{i,j,\xi}$ and since $\beta_{i,h}\in E(X_i)$, we get $\eta_j\mid\eta_i\beta_{i,h}$. Thus $\eta_j\mid\eta_i e_{i,\ell}$, and consequently $\eta_j\mid\eta_i p$.

Suppose now that $m_{i,j,\xi}(f)<m_{i,j,\xi}$. Since $f_i\in\Ocal_{X_i^s}$, every exponent below $m_{i,j,\xi}$ in $f_i(u^{\eta_i})$ is therefore divisible by both $\eta_i$ and $\eta_j$. The same argument applies to $f_j((\xi u)^{\eta_j})$. Then both $\eta_i$ and $\eta_j$ divide $m_{i,j,\xi}(f)$. On the other hand, $m_{i,j,\xi}$ is the order of a difference $g_i(u^{\eta_i})-g_j((\xi u)^{\eta_j})$, so at least one of $\eta_i$ and $\eta_j$ divides $m_{i,j,\xi}$. Hence
$$
m_{i,j,\xi}-m_{i,j,\xi}(f)\geq\min\{\eta_i,\eta_j\}.
$$
By Lemma~\ref{Lem_Bound_mijxi_Delta}, we have
$m_{i,j,\xi}\leq\min\{\eta_i,\eta_j\}(\delta-r+2)$,
and therefore
$$
m_{i,j,\xi}(f)\leq m_{i,j,\xi}-\min\{\eta_i,\eta_j\}
\leq\left(1-\frac1{\delta-r+2}\right)m_{i,j,\xi}
<\alpha m_{i,j,\xi}.
$$
This contradicts $f\in\Ocal_{X^\alpha}$. We conclude that $m_{i,j,\xi}(f)\geq m_{i,j,\xi}$ for every $i\neq j$ and every $\eta$-th root of unity $\xi$. Theorem~\ref{Theo_HolderSatReducibleCurve} now gives $f\in\Ocal_{X^s}$. Thus $\Ocal_{X^\alpha}\subseteq\Ocal_{X^s}$, and the reverse inclusion always holds.
\end{proof}

We also extend Corollary \ref{Cor_HolderSatPreservesMultiplicity} in the case of reducible curves.

\begin{prop}\label{Prop_HolderSatPreservesMultiplicityReduced}
    Let $(X,0)$ be a reduced complex curve germ with $r$ irreducible components and multiplicity $m=m(X,0)$. Then
    \[
    \alpha>\frac{m-r}{m-r+2}
    \quad\Longrightarrow\quad
    m(X^\alpha,0)=m(X,0).
    \]
\end{prop}

\begin{proof}
    Write $X=\cup_{i=1}^r X_i$, and denote by $m_i$ the multiplicity of each $X_i$. Since $m=\sum_{i=1}^rm_i$ and $m_i\geq1$ for every $i$, we have $m_i\leq m-r+1$ for every $i$. Hence
    \[
    \frac{m_i-1}{m_i+1}\leq\frac{m-r}{m-r+2}<\alpha.
    \]
    By Corollary \ref{Cor_HolderSatPreservesMultiplicity}, we get $m(X_i^\alpha)=m_i$ for every $i$. Moreover, by Theorem \ref{Theo_HolderSatReducibleCurve}, we have $\Ocal_{X_i}\subseteq\Ocal_{(X^\alpha)_i}\subseteq\Ocal_{X_i^\alpha}$. Since $\Ocal_{X_i}$ and $\Ocal_{X_i^\alpha}$ have the same multiplicity $m_i$, the intermediate ring $\Ocal_{(X^\alpha)_i}$ also has multiplicity $m_i$. Therefore
    \[
    m(X^\alpha)=\sum_{i=1}^rm((X^\alpha)_i)=\sum_{i=1}^rm_i=m(X).
    \]
\end{proof}

\section{Hölder homeomorphisms between curves}\label{Sec_HolderHomeoCurves}
 In this section we study homeomorphisms between curves. This serves as a warm-up for the study of families of curves in Sections~\ref{Sec_RelativeSpecialization} and~\ref{Sec_HolderEquisaturation}. We obtain two results. The first one says that a homeomorphism between reduced curves can always be improved to a bi-$\alpha$-Hölder homeomorphism for some $\alpha>0$: two curve germs with the same number of branches are bi-$\alpha$-Hölder homeomorphic, and a homeomorphism between projective curves gives rise to a bi-$\alpha$-Hölder homeomorphism between them (Theorem~\ref{Th_HomeoImpliesHolderHomeo}). The proof rests on Theorem~\ref{Theo_FiniteNumberOfAlphas}: the seminormalization is a homeomorphism, and by that theorem it is bi-$\varepsilon$-Hölder for some $\varepsilon>0$. The second result says that bi-$\alpha$-Hölder homeomorphisms with $\alpha$ close to $1$ are rigid. They preserve the multiplicities of the branches, the threshold being governed by the first characteristic exponent, and more generally they constrain the characteristic exponents of the image branches (Proposition~\ref{Prop_ClusterLemma} and Corollary~\ref{Cor_MetricMultiplicityContinuation}). The proof is an elementary cluster argument. The $n$ sheets of a branch of multiplicity $n$ come together in clusters whose sizes and mutual distances are read off from the characteristic sequence. A Hölder homeomorphism cannot tear a cluster apart, so the image branch must have enough sheets to accommodate it.
 
Between these two results lies an interval of exponents where the situation is not understood. Let $B$ be a branch of multiplicity $n$ with characteristic sequence $\beta_0=n<\beta_1<\dots<\beta_g$. By Theorem~\ref{Th_HomeoImpliesHolderHomeo} and Corollary~\ref{Cor_HolderEqNormaIffBeta_g}, $B$ is bi-$\alpha$-Hölder homeomorphic to a disc for every $\alpha\leq1/\beta_g$, whereas by Proposition~\ref{Prop_ClusterLemma} it is not for $\alpha^2>n/\beta_g$. Thus the largest exponent for which $B$ is bi-$\alpha$-Hölder homeomorphic to a smooth germ lies in $[1/\beta_g,\sqrt{n/\beta_g}]$. Even for a branch with a single characteristic pair $(n,\beta_1)$ we can only narrow this interval down to $[n/\beta_1,\sqrt{n/\beta_1}]$ (see Remark~\ref{Rmq_ClusterThresholds}).

\subsection{Hölder maps}
 
Let us fix some terminology. Let $E\subset\C^N$ and $E'\subset\C^{N'}$ be subsets and let $\alpha\in(0,1]$. We say that a map $h : E\to E'$ is \textit{$\alpha$-Hölder} if there exists a constant $C>0$ such that $\|h(p)-h(q)\|\leq C\|p-q\|^\alpha$ for all $p,q\in E$. We say that $h$ is \textit{locally $\alpha$-Hölder} if each point of $E$ has a neighborhood on which $h$ is $\alpha$-Hölder. A \textit{bi-$\alpha$-Hölder homeomorphism} is a homeomorphism $h : E\to E'$ such that both $h$ and $h^{-1}$ are $\alpha$-Hölder. We say that a homeomorphism $h$ is \textit{locally bi-$\alpha$-Hölder} if each point of $E$ has a neighborhood $U$ such that $h|_U : U\to h(U)$ is bi-$\alpha$-Hölder. For germs these notions refer to representatives. The $1$-Hölder maps are the Lipschitz maps. On subsets of a projective space we use the metric induced by the Fubini--Study metric. Locally this metric is bi-Lipschitz equivalent to the Euclidean metric of an affine chart, so the local notions do not depend on this choice. We will use repeatedly the following elementary facts.
 
\begin{lemma}\label{Lem_HolderElementary}
\begin{enumerate}[label=\textup{(\roman*)}]
\item Suppose $h : E\to E'$ is $\alpha$-Hölder and $k : E'\to E''$ is $\beta$-Hölder. Then $k\circ h$ is $\alpha\beta$-Hölder. In particular, the composition of a bi-$\alpha$-Hölder homeomorphism with a bi-Lipschitz homeomorphism is bi-$\alpha$-Hölder. A holomorphic map between open subsets of $\C^N$ and $\C^M$ is locally Lipschitz. Hence an isomorphism of analytic germs is locally bi-Lipschitz, because it is induced by holomorphic maps of the ambient spaces in both directions.
\item Suppose $E$ is compact and $h : E\to E'$ is a locally bi-$\alpha$-Hölder homeomorphism. Then $h$ is bi-$\alpha$-Hölder.
\item Suppose $E$ is bounded and $h : E\to E'$ is $\alpha$-Hölder. Then $h$ is $\beta$-Hölder for every $\beta\leq\alpha$. In particular, if $h$ is a homeomorphism of bounded sets such that $h$ is $\alpha$-Hölder and $h^{-1}$ is $\alpha'$-Hölder, then $h$ is bi-$\min(\alpha,\alpha')$-Hölder.
\end{enumerate}
\end{lemma}
 
\begin{proof}
(i) Suppose $\|h(p)-h(q)\|\leq C\|p-q\|^\alpha$ and $\|k(y)-k(y')\|\leq C'\|y-y'\|^\beta$. Then $\|k(h(p))-k(h(q))\|\leq C'C^\beta\|p-q\|^{\alpha\beta}$. The remaining assertions are clear.
 
(ii) Cover $E$ by finitely many open sets $U_1,\dots,U_m$ such that $h|_{U_i}$ is bi-$\alpha$-Hölder with a constant $C$ independent of $i$. Let $\eta>0$ be a Lebesgue number of this covering. If $\|p-q\|<\eta$, then $p$ and $q$ lie in a common $U_i$, so $\|h(p)-h(q)\|\leq C\|p-q\|^\alpha$. If $\|p-q\|\geq\eta$, then $\|h(p)-h(q)\|\leq\operatorname{diam}(E')\leq\operatorname{diam}(E')\eta^{-\alpha}\|p-q\|^\alpha$. The same argument applies to $h^{-1}$ because $E'=h(E)$ is compact.
 
(iii) Denote by $D$ the diameter of $E$. Then $\|p-q\|^\alpha\leq D^{\alpha-\beta}\|p-q\|^\beta$.
\end{proof}
 
By Corollary~\ref{Cor_HolderSatMapIsHolder} the natural map $\pi_\alpha : X^\alpha\to X$ is a locally bi-$\alpha$-Hölder homeomorphism. Note that $\pi_\alpha$ is a morphism of analytic spaces, because the coordinate functions of $X$ lie in $\Ocal_X\subseteq\Ocal_{X^\alpha}$. Thus by Lemma~\ref{Lem_HolderElementary}~(i) $\pi_\alpha$ is locally Lipschitz; only its inverse is merely $\alpha$-Hölder. We will make use of this asymmetry.
 
\subsection{Homeomorphisms between curves are Hölder}
 
Let $V$ be a reduced curve and let $p\in V$. Denote by $r_p(V)$ the number of analytic branches of $(V,p)$. Let $\pi : \overline V\to V$ be the normalization and let $s : V^+\to V$ be the seminormalization. Recall that $\pi^{-1}(p)$ consists of exactly $r_p(V)$ points, and that if $\pi^{-1}(p)=\{\overline p_1,\dots,\overline p_r\}$, then
\begin{equation}\label{Eq_LocalRingSeminormalization}
\Ocal_{V^+,p}\simeq\Big\{(f_1,\dots,f_r)\in\bigoplus_{j=1}^r\Ocal_{\overline V,\overline p_j}\ \Big|\ f_1(\overline p_1)=\dots=f_r(\overline p_r)\Big\}.
\end{equation}
Thus $s$ is a finite bijective morphism, hence a homeomorphism, and $V^+$ is obtained from $\overline V$ by identifying the points lying over the same point of $V$. Choose a local coordinate $t_j$ at $\overline p_j$. Then the ring in \eqref{Eq_LocalRingSeminormalization} is $\C\{t_1,\dots,t_r\}/(t_it_j)_{i\neq j}$. Therefore $(V^+,p)$ is analytically isomorphic to the germ at $0$ of
\[
A_r:=\bigcup_{j=1}^r\C e_j\subset\C^r,
\]
the union of the $r$ coordinate axes. In particular, the germ of $V^+$ at a point with $r$ branches depends only on $r$.
 
\begin{theorem}\label{Th_HomeoImpliesHolderHomeo}
\begin{enumerate}[label=\textup{(\roman*)}]
\item \emph{Local case.} Suppose $(C_1,0)$ and $(C_2,0)$ are reduced complex analytic curve germs in $(\C^N,0)$ with the same number of branches. Then there exist $\alpha\in(0,1]\cap\Q$ and a bi-$\alpha$-Hölder homeomorphism $(C_1,0)\to(C_2,0)$.
\item \emph{Global case.} Suppose $C_1$ and $C_2$ are reduced projective curves in $\mathbb P^N_\C$ and $h : C_1\to C_2$ is a homeomorphism for the Euclidean topology. Then there exist $\alpha\in(0,1]\cap\Q$ and a bi-$\alpha$-Hölder homeomorphism $\tilde h : C_1\to C_2$, where the Hölder condition refers to the metrics induced by the Fubini--Study metric on $\mathbb P^N_\C$.
\end{enumerate}
In both cases one can take $\alpha=\min(\varepsilon_1,\varepsilon_2)$, where $\varepsilon_i\in(0,1]\cap\Q$ is any number with $C_i^{\varepsilon_i}\simeq C_i^+$.
\end{theorem}
 
\begin{proof}
We first prove (i). Let $r$ be the number of branches of $(C_1,0)$ and $(C_2,0)$, and let $s_i : (C_i^+,0)\to(C_i,0)$ be the seminormalization. By \eqref{Eq_LocalRingSeminormalization} both $(C_1^+,0)$ and $(C_2^+,0)$ are analytically isomorphic to $(A_r,0)$. Hence there exists an analytic isomorphism
\[
d : (C_1^+,0)\longrightarrow(C_2^+,0),
\]
and $d$ is bi-Lipschitz after shrinking the representatives by Lemma~\ref{Lem_HolderElementary}~(i). By Theorem~\ref{Th_HolderSatIsHolderFcts} and the local form of Theorem~\ref{Theo_FiniteNumberOfAlphas}, namely Proposition~\ref{Prop_FiniteNumberOfAlphasCoherentLOCAL}, for each $i\in\{1,2\}$ there exists $\varepsilon_i\in(0,1]\cap\Q$ with $C_i^{\varepsilon_i}\simeq C_i^+$. Then by Corollary~\ref{Cor_HolderSatMapIsHolder} $s_i : C_i^+\to C_i$ is locally bi-$\varepsilon_i$-Hölder. Set
\[
\tilde h:=s_2\circ d\circ s_1^{-1} : (C_1,0)\longrightarrow(C_2,0).
\]
Then $\tilde h$ is a homeomorphism. Because $s_2$ and $d$ are Lipschitz and $s_1^{-1}$ is $\varepsilon_1$-Hölder, $\tilde h$ is $\varepsilon_1$-Hölder by Lemma~\ref{Lem_HolderElementary}~(i). Likewise $\tilde h^{-1}=s_1\circ d^{-1}\circ s_2^{-1}$ is $\varepsilon_2$-Hölder. Therefore $\tilde h$ is bi-$\alpha$-Hölder with $\alpha:=\min(\varepsilon_1,\varepsilon_2)$ by Lemma~\ref{Lem_HolderElementary}~(iii). This proves (i).
 
Next we prove (ii). Let us first record the topological information carried by $h$. Let $V$ be a reduced complex curve and let $p\in V$. Each branch of $(V,p)$ is homeomorphic to a cone over a circle, and distinct branches meet only at $p$. Hence, for $U$ running through a fundamental system of neighborhoods of $p$ in $V$, the punctured neighborhood $U\setminus\{p\}$ has exactly $r_p(V)$ connected components, each of which is the trace of one branch. Thus the number of branches at a point, and the correspondence between branches, are preserved by homeomorphisms.
 
Let $V=\bigcup_{\lambda\in\Lambda}V_\lambda$ be the decomposition of a reduced projective curve into its finitely many irreducible components. Set $B(V):=\{p\in V\mid r_p(V)\geq2\}$; this is a finite set. We show that the connected components of $V\setminus B(V)$ are the sets $V_\lambda\setminus B(V)$. Distinct irreducible components meet only in $B(V)$, so these sets are pairwise disjoint and their union is $V\setminus B(V)$. Each of them is connected: the normalization $\overline V_\lambda$ of $V_\lambda$ is a connected compact Riemann surface, it remains connected after removing the finitely many points of $\pi_\lambda^{-1}(B(V))$, and $V_\lambda\setminus B(V)$ is the continuous image of this set. Finally, each $V_\lambda\setminus B(V)$ is closed in $V\setminus B(V)$, being the trace of the closed set $V_\lambda$, and finitely many pairwise disjoint closed sets covering a space are also open in it. This proves the assertion. Taking closures recovers the irreducible components. Because $h(B(C_1))=B(C_2)$ by the preceding paragraph, $h$ induces a bijection between the connected components of $C_1\setminus B(C_1)$ and those of $C_2\setminus B(C_2)$. After reindexing we can therefore write
\[
C_i=\bigcup_{\lambda\in\Lambda}C_{i,\lambda}
\]
with $h(C_{1,\lambda})=C_{2,\lambda}$ for every $\lambda\in\Lambda$. Note that $h$ also preserves the incidence between branches and irreducible components, because a branch at $p$ is contained in a unique irreducible component.
 
Let $\pi_{i,\lambda} : \overline C_{i,\lambda}\to C_{i,\lambda}$ be the normalization. We show that $\overline C_{1,\lambda}$ and $\overline C_{2,\lambda}$ have the same genus. For $p\in C_{i,\lambda}$ denote by $r_{i,\lambda}(p)$ the number of branches of $C_{i,\lambda}$ at $p$. Because $h|_{C_{1,\lambda}} : C_{1,\lambda}\to C_{2,\lambda}$ is a homeomorphism preserving branches, we have
\[
\chi(C_{1,\lambda})=\chi(C_{2,\lambda})
\qquad\text{and}\qquad
r_{1,\lambda}(p)=r_{2,\lambda}(h(p))\quad\text{for every }p\in C_{1,\lambda}.
\]
Recall that for an irreducible projective curve $C$ with normalization $\pi : \overline C\to C$ one has
\begin{equation}\label{Eq_EulerCharNormalization}
\chi(\overline C)=\chi(C)+\sum_{p\in C}\bigl(r_p(C)-1\bigr).
\end{equation}
Indeed, let $B\subset C$ be the finite set of points with at least two branches. Then $\pi$ restricts to a homeomorphism $\overline C\setminus\pi^{-1}(B)\to C\setminus B$. Triangulate $C$ with the points of $B$ among the vertices and lift the triangulation to $\overline C$. Then $\pi$ is bijective on simplices, except that each vertex $p\in B$ has $r_p(C)$ preimages, which gives \eqref{Eq_EulerCharNormalization}. It follows that $\chi(\overline C_{1,\lambda})=\chi(\overline C_{2,\lambda})$. Because $\overline C_{i,\lambda}$ is a compact connected Riemann surface and $\chi(\overline C_{i,\lambda})=2-2g(\overline C_{i,\lambda})$, the two genera are equal.
 
We now use the correspondence between branches to mark the normalizations. Let
\[
\pi_i : \overline C_i=\coprod_{\lambda\in\Lambda}\overline C_{i,\lambda}\longrightarrow C_i
\]
be the normalization of $C_i$ and set $B_i:=\{p\in C_i\mid\#\pi_i^{-1}(p)\geq2\}=B(C_i)$. For every $p\in B_1$ the homeomorphism $h$ induces a bijection between the branches of $(C_1,p)$ and those of $(C_2,h(p))$, hence a bijection
\[
\pi_1^{-1}(p)\longrightarrow\pi_2^{-1}(h(p)),
\]
which we call the prescribed bijection. Because $h$ preserves the incidence between branches and irreducible components, the prescribed bijection sends the points of $\overline C_{1,\lambda}$ to points of $\overline C_{2,\lambda}$. Set
\[
P_{i,\lambda}:=\overline C_{i,\lambda}\cap\pi_i^{-1}(B_i).
\]
The prescribed bijections give a bijection $P_{1,\lambda}\to P_{2,\lambda}$ between these finite sets. Because $\overline C_{1,\lambda}$ and $\overline C_{2,\lambda}$ have the same genus, the classification of compact oriented surfaces, together with the standard point-pushing argument (the diffeomorphism group of a connected surface acts transitively on ordered tuples of distinct points), gives an orientation-preserving diffeomorphism
\[
d_\lambda : \overline C_{1,\lambda}\longrightarrow\overline C_{2,\lambda}
\]
which restricts to the bijection $P_{1,\lambda}\to P_{2,\lambda}$. Set
\[
d_{\mathrm{nor}}:=\coprod_{\lambda\in\Lambda}d_\lambda : \overline C_1\longrightarrow\overline C_2 .
\]
By construction $d_{\mathrm{nor}}$ is a bijection, $d_{\mathrm{nor}}(\pi_1^{-1}(B_1))=\pi_2^{-1}(B_2)$, and $d_{\mathrm{nor}}(\pi_1^{-1}(p))=\pi_2^{-1}(h(p))$ for every $p\in B_1$.
 
Let $\rho_i : \overline C_i\to C_i^+$ be the normalization map factored through the seminormalization. By \eqref{Eq_LocalRingSeminormalization}, as a topological space $C_i^+$ is the quotient of $\overline C_i$ by the equivalence relation $x\sim_i y$ if and only if $\pi_i(x)=\pi_i(y)$. The map $d_{\mathrm{nor}}$ sends each fiber of $\pi_1$ with at least two points bijectively onto the corresponding fiber of $\pi_2$, and it sends the complement of $\pi_1^{-1}(B_1)$ onto the complement of $\pi_2^{-1}(B_2)$, where the fibers are singletons. Therefore
\[
x\sim_1y\quad\Longleftrightarrow\quad d_{\mathrm{nor}}(x)\sim_2d_{\mathrm{nor}}(y),
\]
so $d_{\mathrm{nor}}$ induces a bijection $d : C_1^+\to C_2^+$ with $d\circ\rho_1=\rho_2\circ d_{\mathrm{nor}}$. Because $\rho_1$ is a quotient map and $\rho_2\circ d_{\mathrm{nor}}$ is continuous, $d$ is continuous. Because $C_1^+$ is compact and $C_2^+$ is Hausdorff, $d$ is a homeomorphism.
 
We show that $d$ is bi-Lipschitz. The curves $C_i^+$ are projective, because finite morphisms are projective, and we use on them the metrics induced by projective embeddings. By Lemma~\ref{Lem_HolderElementary}~(i) and (ii) two such embeddings give bi-Lipschitz equivalent metrics on the compact curve $C_i^+$, so the assertion does not depend on this choice. Away from the points with at least two branches, $\rho_i$ is a local analytic isomorphism by \eqref{Eq_LocalRingSeminormalization}, so $d$ is locally given by the diffeomorphisms $d_\lambda$ read in analytic coordinates. Because a holomorphic chart of a smooth point of an embedded curve is bi-Lipschitz onto its image for the ambient metric, $d$ is locally bi-Lipschitz at these points. Now let $p\in C_1^+$ be a point with $r\geq2$ branches and set $q:=d(p)$. By the local description of the seminormalization there are analytic identifications
\[
(C_1^+,p)\simeq(A_r,0),\qquad(C_2^+,q)\simeq(A_r,0),
\]
in which the $j$-th axis corresponds to the branch through the $j$-th point of $\pi_1^{-1}(s_1(p))$, respectively of $\pi_2^{-1}(s_2(q))$. Because $d_{\mathrm{nor}}$ maps the points over $s_1(p)$ bijectively onto the points over $s_2(q)$, after reordering the target axes $d$ takes the form
\[
d(te_j)=\varphi_j(t)e_j ,
\]
where $\varphi_j : (\C,0)\to(\C,0)$ is a smooth real diffeomorphism germ fixing the origin, namely $d_\lambda$ read in the local coordinates of $\overline C_1$ and $\overline C_2$ at the corresponding marked points. After shrinking, there are constants $0<c\leq C$ such that for all $j$
\[
c|t-s|\leq|\varphi_j(t)-\varphi_j(s)|\leq C|t-s| .
\]
If $x=te_j$ and $y=se_j$ lie on the same branch, then $c\|x-y\|\leq\|d(x)-d(y)\|\leq C\|x-y\|$. If $x=te_j$ and $y=se_k$ with $j\neq k$, then
\[
\|x-y\|^2=|t|^2+|s|^2
\qquad\text{and}\qquad
\|d(x)-d(y)\|^2=|\varphi_j(t)|^2+|\varphi_k(s)|^2 ,
\]
and because $\varphi_j(0)=\varphi_k(0)=0$ we get
\[
c^2\bigl(|t|^2+|s|^2\bigr)\leq\|d(x)-d(y)\|^2\leq C^2\bigl(|t|^2+|s|^2\bigr).
\]
Thus $d$ is locally bi-Lipschitz at $p$ as well. Because the analytic identifications above and their inverses are locally Lipschitz, $d$ is locally bi-Lipschitz everywhere. Hence $d$ is bi-Lipschitz by Lemma~\ref{Lem_HolderElementary}~(ii).
 
Finally, by Theorems~\ref{Th_HolderSatIsHolderFcts} and~\ref{Theo_FiniteNumberOfAlphas}, there exists $\varepsilon_i\in(0,1]\cap\Q$ with $C_i^{\varepsilon_i}\simeq C_i^+$. By Corollary~\ref{Cor_HolderSatMapIsHolder} $s_i : C_i^+\to C_i$ is locally bi-$\varepsilon_i$-Hölder, hence bi-$\varepsilon_i$-Hölder by Lemma~\ref{Lem_HolderElementary}~(ii). Set
\[
\tilde h:=s_2\circ d\circ s_1^{-1} : C_1\longrightarrow C_2 .
\]
Then $\tilde h$ is a homeomorphism. Because $s_2$ and $d$ are Lipschitz and $s_1^{-1}$ is $\varepsilon_1$-Hölder, $\tilde h$ is $\varepsilon_1$-Hölder, and likewise $\tilde h^{-1}=s_1\circ d^{-1}\circ s_2^{-1}$ is $\varepsilon_2$-Hölder, by Lemma~\ref{Lem_HolderElementary}~(i). Therefore $\tilde h$ is bi-$\alpha$-Hölder with $\alpha:=\min(\varepsilon_1,\varepsilon_2)$ by Lemma~\ref{Lem_HolderElementary}~(iii). This proves (ii).
\end{proof}
 
\begin{rmq}\label{Rmq_HomeoHolderExponent}
(i) The exponent in Theorem~\ref{Th_HomeoImpliesHolderHomeo} is the exponent of the seminormalizations. For $\varepsilon_i$ one can take the largest $\alpha\in(0,1]$ with $C_i^\alpha\simeq C_i^+$, which is the first threshold $\alpha_1$ of Proposition~\ref{Prop_FiniteNumberOfAlphasCoherentLOCAL}, respectively of Theorem~\ref{Theo_FiniteNumberOfAlphas}, for $C_i$. In particular, the exponent obtained depends only on $C_1$ and $C_2$ and not on the homeomorphism $h$. For a branch one has $\varepsilon_i=1/\beta_g$ by Corollary~\ref{Cor_HolderEqNormaIffBeta_g}.
 
(ii) The bi-$\alpha$-Hölder homeomorphism $\tilde h$ is not the given homeomorphism $h$. It is built from the topological data determined by $h$, namely the irreducible components, the branches, and their incidences. In the local case, when both germs are branches, $\tilde h$ is the composition of the normalizations with a biholomorphism of discs. In the global case, when every singular point of every irreducible component of $C_1$ is unibranch, each $C_{1,\lambda}$ is a topological surface and $d_{\mathrm{nor}}$ is simply a diffeomorphism realizing the prescribed bijections between marked points.
 
(iii) Theorem~\ref{Th_HomeoImpliesHolderHomeo} says that for curves the existence of a bi-$\alpha$-Hölder homeomorphism for some unspecified exponent $\alpha>0$ is a purely topological condition: it is equivalent to the existence of a homeomorphism, and for germs to the equality of the numbers of branches. In other words, the Hölder geometry of curves carries no information unless one keeps track of the exponent. The proof also shows where the exponent is lost. The seminormalizations $C_1^+$ and $C_2^+$ are bi-Lipschitz homeomorphic, because their singularities are ordinary multiple points, which are determined by their numbers of branches. The only non-Lipschitz maps in the construction are the inverses of the seminormalization maps $s_i$. Once the exponent is fixed, rigidity appears. Theorem~\ref{Th_HolderSatSmoothness} shows that an $\alpha$-Hölder saturation with $\alpha>1/3$ already detects smoothness, and the next two subsections show that bi-$\alpha$-Hölder homeomorphisms with $\alpha$ large detect the multiplicity.
\end{rmq}
 
\subsection{Clusters of points on a branch}\label{Subsec_Clusters}
 
\textbf{Setup.} Throughout this subsection $B\subset(\C^N,0)$ is a branch, that is, an irreducible reduced curve germ. Denote by $n:=m(B)$ its multiplicity and set $A:=\Ocal_{B,0}\subset\overline A=\C\{t\}$. As in Section~3, we fix a good uniformizer $t$, so that there is a transverse parameter $x\in A$ with $x=t^n$, and we denote by
\[
\beta_0=n<\beta_1<\dots<\beta_g,\qquad e_i:=\gcd(\beta_0,\dots,\beta_i),
\]
the characteristic sequence of the set of exponents $E(A)$. Recall that $\beta_0=\min(E(A)\setminus\{0\})$, that $\beta_i=\min\{p\in E(A)\mid e_{i-1}\nmid p\}$ for $1\leq i\leq g$, that $e_0=n>e_1>\dots>e_g=1$ with $e_i\mid e_{i-1}$, and that $g=0$ if and only if $B$ is smooth. By \cite[Lemma VI.3.3]{BrianconGallicoGranger1980DeformationsEqu} the $\beta_i$ are the characteristic exponents of a generic plane projection of $B$. Let
\[
\pi : (\C,0)\longrightarrow(B,0),\qquad t\longmapsto\pi(t)=(\pi_1(t),\dots,\pi_N(t)),
\]
be the normalization. After shrinking, $\pi$ is a homeomorphism from a disc onto a representative of $B$; in particular $\pi(t)\neq\pi(t')$ for $t\neq t'$. Each element of $A$ is of the form $f\circ\pi$ for a holomorphic function germ $f$ on $(\C^N,0)$, and conversely. In particular, the components $\pi_k=z_k\circ\pi$ lie in $A$, so $\operatorname{Ex}(\pi_k)\subseteq E(A)$.
 
Two consequences of the setup will be used constantly. First, all exponents occurring in the $\pi_k$ are at least $n$, and $|t|^n=|x(\pi(t))-x(0)|$ with $x$ Lipschitz on $(\C^N,0)$. Hence there are constants $0<c\leq C$ with
\begin{equation}\label{Eq_NormOnBranch}
c|t|^n\leq\|\pi(t)\|\leq C|t|^n .
\end{equation}
Second, for each $t$ the $n$ points $\pi(\xi t)$, $\xi\in\mu_n:=\{\xi\in\C\mid\xi^n=1\}$, are exactly the points of $B$ lying over the value $x=t^n$ of the transverse parameter. We call them the \textit{sheets} of $B$ over $t^n$. For a plane branch they are the values at $x=t^n$ of the $n$ conjugate Puiseux series $y(\xi x^{1/n})$. The fact that their mutual distances are governed by the characteristic exponents goes back to Zariski's theory of saturation and to Pham and Teissier \cite{PhamTeissier1969LipFractions}; see \cite[Chapter~VI]{BrianconGallicoGranger1980DeformationsEqu} and \cite{FGST2020LipschitzAnAlgebraicApproch} for space branches, and \cite{Fernandes2003TopologicalEquivalence,NeumannPichon2014} for the use of these contact orders in Lipschitz geometry. In our setting they are encoded by the following numbers.
 
\begin{definition}\label{Def_BetaXi}
For $\xi\in\mu_n\setminus\{1\}$ set
\[
i(\xi):=\min\{i\geq1\mid\xi^{e_i}\neq1\}\in\{1,\dots,g\}\qquad\text{and}\qquad\beta_B(\xi):=\beta_{i(\xi)},
\]
and set $\beta_B(1):=+\infty$. For $\gamma>0$ define
\[
N_B(\gamma):=\#\{\xi\in\mu_n\mid\beta_B(\xi)/n\geq\gamma\}.
\]
\end{definition}
 
The index $i(\xi)$ is well defined because $\xi^{e_0}=\xi^n=1$ and $\xi^{e_g}=\xi\neq1$. Note that $\xi^{e_i}=1$ if and only if $\xi\in\mu_{e_i}$. Thus $\{\xi\in\mu_n\mid i(\xi)=i\}=\mu_{e_{i-1}}\setminus\mu_{e_i}$ for $1\leq i\leq g$, and
\[
\#\{\xi\in\mu_n\mid\beta_B(\xi)\geq\beta_i\}=\#\mu_{e_{i-1}}=e_{i-1}.
\]
Hence $N_B(\gamma)=e_{i-1}$ for $\beta_{i-1}/n<\gamma\leq\beta_i/n$ and $1\leq i\leq g$, $N_B(\gamma)=n$ for $\gamma\leq1$, and $N_B(\gamma)=1$ for $\gamma>\beta_g/n$. For a smooth branch $N_B\equiv1$. In the notation of Section~3, $\beta_B(\xi)$ is the order of the diagonal ideal $I_\Delta$ of $B$ along the arc $u\mapsto(u,\xi u)$ of $\overline B\times\overline B$, that is, the number $m_{\varphi_\xi}=\min_{f\in A}\ord_u\bigl(f(\xi u)-f(u)\bigr)$ attached there to the arcs $\varphi_{i,j,\xi}$ in the case of a single branch. This is what the proof of the next lemma computes. The rational numbers $\beta_B(\xi)/n$ are the contact orders between the sheets of $B$.
 
\begin{lemma}\label{Lem_SheetDistances}
Let $\xi\in\mu_n\setminus\{1\}$. There exist constants $0<c\leq C$ such that for all $t$ with $|t|$ small enough
\[
c|t|^{\beta_B(\xi)}\leq\|\pi(\xi t)-\pi(t)\|\leq C|t|^{\beta_B(\xi)} .
\]
\end{lemma}
 
\begin{proof}
Let $f\circ\pi=\sum_{p\in\operatorname{Ex}(f\circ\pi)}a_pt^p$ be an element of $A$. Then
\[
f(\pi(\xi t))-f(\pi(t))=\sum_pa_p(\xi^p-1)t^p .
\]
Let $p\in E(A)$ with $p<\beta_{i(\xi)}$. Then $e_{i(\xi)-1}\mid p$ by the definition of $\beta_{i(\xi)}$, and $\xi^{e_{i(\xi)-1}}=1$ by the definition of $i(\xi)$. Hence $\xi^p=1$, and all terms of order less than $\beta_B(\xi)$ cancel. Applied to the coordinate functions $f=z_k$ this gives the upper bound. On the other hand, $\xi^{\beta_{i(\xi)}}\neq1$, for otherwise $\xi$ would be killed by $\gcd(e_{i(\xi)-1},\beta_{i(\xi)})=e_{i(\xi)}$. Because $E(A)$ is the union of the sets $\operatorname{Ex}(f\circ\pi)$, there is an $f$ with $\beta_{i(\xi)}\in\operatorname{Ex}(f\circ\pi)$. For this $f$ the difference $f(\pi(\xi t))-f(\pi(t))$ has order exactly $\beta_B(\xi)$. Let $L$ be a Lipschitz constant of $f$ near $0$. Then $|f(\pi(\xi t))-f(\pi(t))|\leq L\|\pi(\xi t)-\pi(t)\|$, which gives the lower bound.
\end{proof}
 
\begin{definition}\label{Def_Cluster}
Let $1\leq j\leq g$ and let $t\neq0$ be small. The \textit{$j$-th cluster of $B$ at level $t$} is the set
\[
\mathrm{Cl}_j(t):=\{\pi(\xi t)\mid\xi\in\mu_{e_{j-1}}\}\subset B
\]
of the sheets of $B$ over $t^n$ indexed by the $e_{j-1}$-th roots of unity.
\end{definition}
 
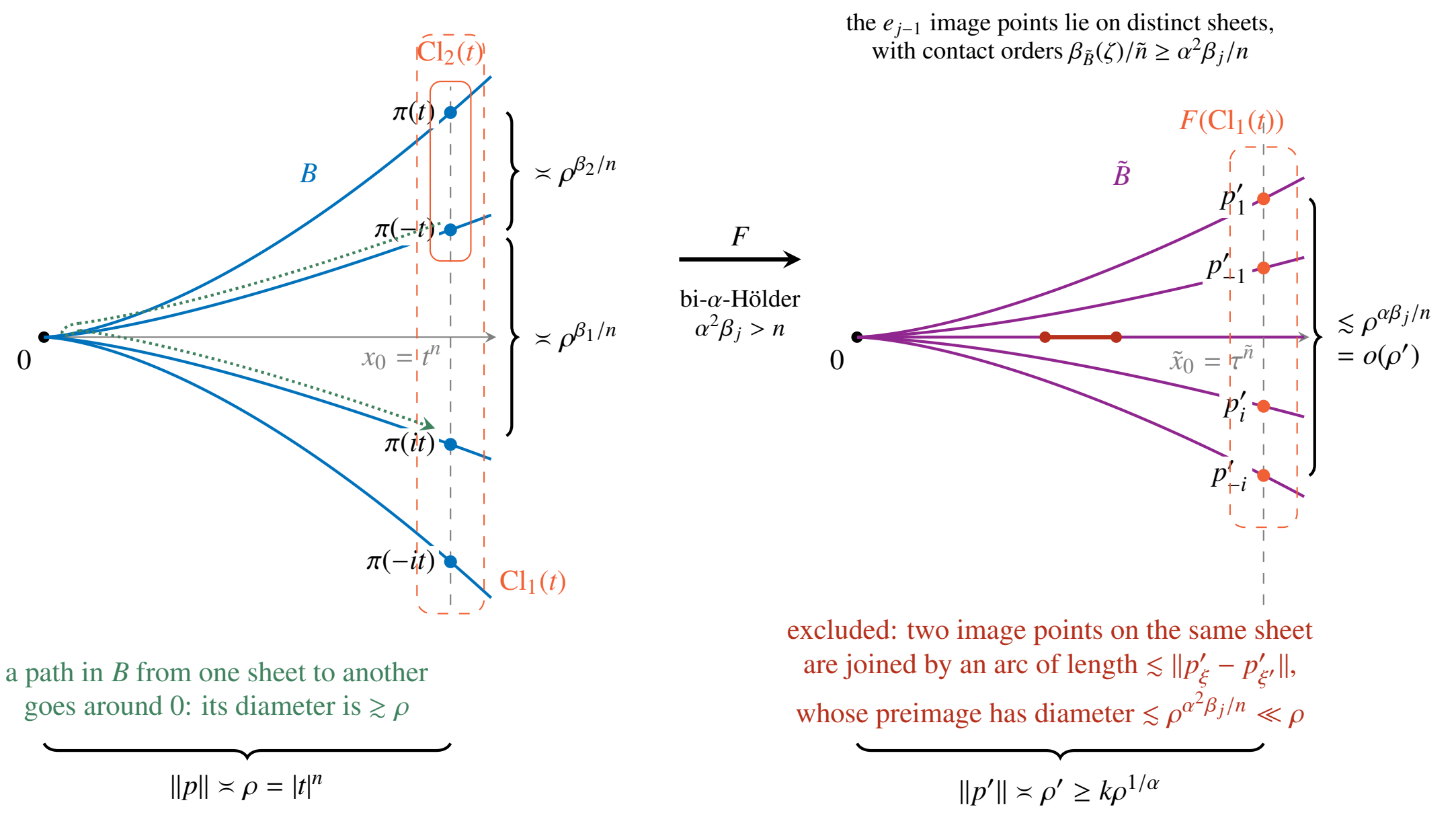
\begin{figure}[htbp]
\centering
\begin{tikzpicture}[x=1cm,y=1.7cm,>=stealth,font=\footnotesize]
\begin{scope}
  \def\xo{4.0}
  \draw[->,gray] (0,0) -- (4.45,0);
  \fill (0,0) circle (1.6pt) node[below left] {$0$};
  \foreach \s/\u in {1/1,1/-1,-1/1,-1/-1}
    \draw[thick,RoyalBlue] plot[domain=0:4.4,samples=80] (\x,{\s*0.12*pow(\x,1.5)+\u*0.03*pow(\x,1.75)});
  \node[RoyalBlue] at (2.6,0.95) {$B$};
  \draw[dashed,gray] (\xo,-1.55) -- (\xo,1.45);
  \node[below left,gray,inner sep=2pt] at (\xo,0) {$x_0=t^n$};
  \coordinate (p1) at (\xo,{ 0.12*pow(\xo,1.5)+0.03*pow(\xo,1.75)});
  \coordinate (p2) at (\xo,{ 0.12*pow(\xo,1.5)-0.03*pow(\xo,1.75)});
  \coordinate (p3) at (\xo,{-0.12*pow(\xo,1.5)+0.03*pow(\xo,1.75)});
  \coordinate (p4) at (\xo,{-0.12*pow(\xo,1.5)-0.03*pow(\xo,1.75)});
  \foreach \p in {p1,p2,p3,p4} \fill[RoyalBlue] (\p) circle (1.8pt);
  \node[left=3pt,fill=white,inner sep=1pt] at (p1) {$\pi(t)$};
  \node[left=3pt,fill=white,inner sep=1pt] at (p2) {$\pi(-t)$};
  \node[left=3pt,fill=white,inner sep=1pt] at (p3) {$\pi(it)$};
  \node[left=3pt,fill=white,inner sep=1pt] at (p4) {$\pi(-it)$};
  \draw[decorate,decoration={brace,amplitude=3pt},thick] ($(p1)+(0.55,0)$) -- ($(p2)+(0.55,0)$)
     node[midway,right=4pt] {$\asymp\rho^{\beta_2/n}$};
  \draw[decorate,decoration={brace,amplitude=3pt},thick] ($(p2)+(0.55,-0.05)$) -- ($(p3)+(0.55,0.05)$)
     node[midway,right=4pt] {$\asymp\rho^{\beta_1/n}$};
  \draw[decorate,decoration={brace,amplitude=4pt,mirror},thick] (0,-2.35) -- (\xo,-2.35)
     node[midway,below=5pt] {$\|p\|\asymp\rho=|t|^n$};
  \draw[rounded corners=3pt,RedOrange] ($(p1)+(-0.2,0.18)$) rectangle ($(p2)+(0.2,-0.18)$);
  \node[RedOrange,above] at ($(p1)+(0,0.2)$) {$\mathrm{Cl}_2(t)$};
  \draw[rounded corners=4pt,RedOrange,dashed] ($(p1)+(-0.33,0.45)$) rectangle ($(p4)+(0.33,-0.3)$);
  \node[RedOrange,right] at ($(p4)+(0.36,-0.12)$) {$\mathrm{Cl}_1(t)$};
  \draw[thick,Green!60!black,densely dotted]
     plot[domain=3.9:0.35,samples=60] (\x,{ 0.12*pow(\x,1.5)-0.03*pow(\x,1.75)+0.06});
  \draw[->,thick,Green!60!black,densely dotted]
     (0.35,{0.12*pow(0.35,1.5)-0.03*pow(0.35,1.75)+0.06}) .. controls (0.12,0.07) and (0.12,-0.04) ..
     (0.35,{-0.12*pow(0.35,1.5)+0.03*pow(0.35,1.75)+0.06})
     plot[domain=0.35:3.85,samples=60] (\x,{-0.12*pow(\x,1.5)+0.03*pow(\x,1.75)+0.06});
  \node[Green!60!black,align=center] at (1.7,-2.05) {a path in $B$ from one sheet to another\\ goes around $0$: its diameter is $\gtrsim\rho$};
\end{scope}
\draw[->,very thick] (6.25,0.45) -- (7.45,0.45) node[midway,above] {$F$};
\node[align=center,font=\scriptsize] at (6.85,0.12) {bi-$\alpha$-H\"older\\ $\alpha^2\beta_j>n$};
\begin{scope}[xshift=8.0cm]
  \def\xt{4.0}
  \draw[->,gray] (0,0) -- (4.45,0);
  \fill (0,0) circle (1.6pt) node[below left] {$0$};
  \foreach \a in {-0.1,-0.05,0,0.05,0.1}
     \draw[thick,Plum] plot[domain=0:4.4,samples=60] (\x,{\a*pow(\x,1.5)});
  \node[Plum] at (2.6,0.95) {$\tilde B$};
  \draw[dashed,gray] (\xt,-1.55) -- (\xt,1.25);
  \node[below left,gray,inner sep=2pt] at (\xt,0) {$\tilde x_0=\tau^{\tilde n}$};
  \draw[decorate,decoration={brace,amplitude=4pt,mirror},thick] (0,-2.35) -- (\xt,-2.35)
     node[midway,below=5pt] {$\|p'\|\asymp\rho'\geq k\rho^{1/\alpha}$};
  \foreach \a/\lab in {0.1/{p'_{1}},0.05/{p'_{-1}},-0.05/{p'_{i}},-0.1/{p'_{-i}}} {
     \fill[RedOrange] (\xt,{\a*pow(\xt,1.5)}) circle (1.8pt) node[left=3pt,black,fill=white,inner sep=1pt] {$\lab$};
  }
  \draw[decorate,decoration={brace,amplitude=3pt},thick] (\xt+0.45,{0.1*pow(\xt,1.5)}) -- (\xt+0.45,{-0.1*pow(\xt,1.5)})
     node[midway,right=4pt,align=left] {$\lesssim\rho^{\alpha\beta_j/n}$\\ $=o(\rho')$};
  \draw[rounded corners=4pt,RedOrange,dashed] (\xt-0.33,{0.1*pow(\xt,1.5)+0.3}) rectangle (\xt+0.33,{-0.1*pow(\xt,1.5)-0.3});
  \node[RedOrange,above left] at (\xt+0.33,{0.1*pow(\xt,1.5)+0.3}) {$F(\mathrm{Cl}_1(t))$};
  \node[align=center,font=\scriptsize] at (2.0,1.72) {the $e_{j-1}$ image points lie on distinct sheets,\\ with contact orders $\beta_{\tilde B}(\zeta)/\tilde n\geq\alpha^2\beta_j/n$};
  \coordinate (q1) at (1.85,0);
  \coordinate (q2) at (2.55,0);
  \draw[very thick,BrickRed] (q1) -- (q2);
  \fill[BrickRed] (q1) circle (1.6pt);
  \fill[BrickRed] (q2) circle (1.6pt);
  \node[BrickRed,align=center] at (1.9,-1.95) {excluded: two image points on the same sheet\\ are joined by an arc of length $\lesssim\|p'_\xi-p'_{\xi'}\|$,\\ whose preimage has diameter $\lesssim\rho^{\alpha^2\beta_j/n}\ll\rho$};
\end{scope}
\end{tikzpicture}
\caption{The cluster argument for $n=4$, $\beta_1=6$, $\beta_2=7$, so that $e_0=4$ and $e_1=2$. Left: the four sheets of $B$ over $x_0=t^n$ form the nested clusters $\mathrm{Cl}_1(t)\supset\mathrm{Cl}_2(t)$, at mutual distances $\rho^{\beta_1/n}$ and $\rho^{\beta_2/n}$, much smaller than their distance $\rho$ to the origin, while any path in $B$ from one sheet to another has to go around the origin. Right: the image of $\mathrm{Cl}_1(t)$ under $F$ is still a cluster, of size $o(\rho')$ because $\alpha\beta_j/n>1/\alpha$. Two image points on the same sheet of $\tilde B$ could be joined by a short arc along that sheet, whose preimage would be a short path between two sheets of $B$. Thus the image points occupy $e_{j-1}$ distinct sheets, whose contact orders are bounded below by $\alpha^2\beta_j/n$.}
\label{Fig_Clusters}
\end{figure}

Thus $\mathrm{Cl}_1(t)$ consists of all $n$ sheets over $t^n$, and the clusters are nested: $\mathrm{Cl}_1(t)\supset\mathrm{Cl}_2(t)\supset\dots\supset\mathrm{Cl}_g(t)\ni\pi(t)$. The sheets in $\mathrm{Cl}_j(t)$ are the sheets whose contact order with $\pi(t)$ is at least $\beta_j/n$. The cluster $\mathrm{Cl}_j(t)$ splits into the $e_{j-1}/e_j$ groups $\{\pi(\xi\zeta t)\mid\zeta\in\mu_{e_j}\}$, $\xi\in\mu_{e_{j-1}}/\mu_{e_j}$, of $e_j$ sheets each; these groups are the translates of $\mathrm{Cl}_{j+1}(t)$ when $j<g$ and single points when $j=g$, and two sheets in different groups are at distance $\asymp|t|^{\beta_j}$. This nested grouping of the sheets of a branch according to their contact orders is the structure described by L\^e's carrousel \cite{Le1978Carrousel,Le2003Carousels}, which has become a basic tool in the Lipschitz geometry of curves and surfaces (see \cite{NeumannPichon2014,BirbrairNeumannPichon2014}). We only need the following three properties of the clusters, all consequences of Lemma~\ref{Lem_SheetDistances}. Set $\rho:=|t|^n$.
\begin{enumerate}[label=\textup{(C\arabic*)}]
\item\label{C1} $\mathrm{Cl}_j(t)$ consists of $e_{j-1}$ distinct points, all at distance $\asymp\rho$ from the origin by \eqref{Eq_NormOnBranch}.
\item\label{C2} The points of $\mathrm{Cl}_j(t)$ are pairwise at distance at most $C\rho^{\beta_j/n}$. Indeed, let $\xi\neq\xi'$ in $\mu_{e_{j-1}}$. Then $\xi'\xi^{-1}\in\mu_{e_{j-1}}\setminus\{1\}$, so $i(\xi'\xi^{-1})\geq j$, and $\|\pi(\xi't)-\pi(\xi t)\|=\|\pi(\xi'\xi^{-1}\cdot\xi t)-\pi(\xi t)\|\asymp|t|^{\beta_B(\xi'\xi^{-1})}\leq|t|^{\beta_j}$.
\item\label{C3} The points of $\mathrm{Cl}_j(t)$ cannot be joined by short paths. This is the content of Lemma~\ref{Lem_ClusterSeparation} below.
\end{enumerate}
Property \ref{C3} is a form of the well-known fact that the inner metric of a branch is bi-Lipschitz equivalent to the metric of a cone over a circle, so that the inner distance between two distinct sheets over $t^n$ is comparable to their distance to the origin, whatever their outer distance is (see for instance \cite{NeumannPichon2014}). We give the elementary argument, which is the monodromy of the Puiseux parametrization. The clusters of a branch and the properties \ref{C1}--\ref{C3} are illustrated on the left of Figure~\ref{Fig_Clusters}; the right-hand side of the figure illustrates the proof of Proposition~\ref{Prop_ClusterLemma} below.
 
\begin{lemma}[separation of the sheets]\label{Lem_ClusterSeparation}
Let $\mathcal C$ be a reduced curve germ containing $B$ as a branch and let $\xi\neq\xi'$ be in $\mu_n$. There exists a constant $c_0>0$ such that for all $t\neq0$ with $|t|$ small enough every path $\gamma : [0,1]\to\mathcal C$ with $\gamma(0)=\pi(\xi t)$ and $\gamma(1)=\pi(\xi't)$ has diameter at least $c_0|t|^n$.
\end{lemma}
 
\begin{proof}
Suppose $\gamma$ passes through the origin. Then its diameter is at least $\|\pi(\xi t)\|\geq c|t|^n$ by \eqref{Eq_NormOnBranch}. Suppose $\gamma$ does not pass through the origin. Then $\gamma([0,1])$ is a connected subset of $\mathcal C\setminus\{0\}$, so it is contained in a single connected component of $\mathcal C\setminus\{0\}$. Because the branches of $\mathcal C$ are the closures of these components, this component is $B\setminus\{0\}$. Because $\pi$ restricts to a homeomorphism from the punctured disc onto $B\setminus\{0\}$, the path $\gamma$ lifts to a path $\tilde\gamma : [0,1]\to\C\setminus\{0\}$ with $\pi\circ\tilde\gamma=\gamma$, $\tilde\gamma(0)=\xi t$ and $\tilde\gamma(1)=\xi't$. Then $x\circ\gamma=\tilde\gamma^{\,n}$ is a closed path in $\C\setminus\{0\}$ based at $x_0:=t^n$. Its winding number around $0$ is $n$ times the variation of the argument of $\tilde\gamma$ divided by $2\pi$. Write $\xi'\xi^{-1}=e^{2\pi ia/n}$ with $0<a<n$. Then this winding number is $a+nk$ for some integer $k$, hence it is nonzero. A closed path contained in the disc $\{w\in\C\mid|w-x_0|<|x_0|\}$, which does not contain $0$, has winding number zero around $0$. Hence $x\circ\gamma$ contains a point $w$ with $|w-x_0|\geq|x_0|=|t|^n$. Let $L$ be a Lipschitz constant of $x$ on $(\C^N,0)$. Then the point $\gamma(s)$ with $x(\gamma(s))=w$ satisfies $\|\gamma(s)-\gamma(0)\|\geq|w-x_0|/L\geq|t|^n/L$. Set $c_0:=\min(c,1/L)$.
\end{proof}
 
\subsection{The cluster proposition}
 
We can now state the main result of this section. Fix a bi-$\alpha$-Hölder homeomorphism $F : \mathcal C\to\tilde{\mathcal C}$ of reduced curve germs with $F(0)=0$, and fix constants $K_1$ and $K_2$ such that
\begin{equation}\label{Eq_HolderConstantsF}
\|F(p)-F(q)\|\leq K_1\|p-q\|^\alpha\qquad\text{and}\qquad\|F^{-1}(p')-F^{-1}(q')\|\leq K_2\|p'-q'\|^\alpha
\end{equation}
on representatives. Because the branches of a reduced curve germ are the closures of the connected components of the punctured germ, $F$ maps each branch of $\mathcal C$ homeomorphically onto a branch of $\tilde{\mathcal C}$. For a branch $\tilde B$ of $\tilde{\mathcal C}$ we use the notation of Subsection~\ref{Subsec_Clusters} decorated with a tilde: $\tilde n$, $\tilde\beta_0<\dots<\tilde\beta_{\tilde g}$, $\tilde e_i$, $\tilde\pi$, $\tilde x$, $\tau$, $\beta_{\tilde B}(\zeta)$ for $\zeta\in\mu_{\tilde n}$, and $N_{\tilde B}(\gamma)$.
 
\begin{prop}[cluster proposition]\label{Prop_ClusterLemma}\label{Lem_MetricClusterContinuation}
Let $F : \mathcal C\to\tilde{\mathcal C}$ be a bi-$\alpha$-Hölder homeomorphism of reduced curve germs. Let $B$ be a singular branch of $\mathcal C$ of multiplicity $n$, with characteristic sequence $\beta_0=n<\beta_1<\dots<\beta_g$ and $e_i=\gcd(\beta_0,\dots,\beta_i)$, and set $\tilde B:=F(B)$. Suppose $\alpha^2\beta_j>n$ for some $1\leq j\leq g$. Then
\[
e_{j-1}\leq N_{\tilde B}\bigl(\alpha^2\beta_j/n\bigr).
\]
In particular, the following holds.
\begin{enumerate}[label=\textup{(\roman*)},leftmargin=*]
\item If $\alpha^2>n/\beta_1$, then $m(\tilde B)\geq n$.
\item If $\alpha^2>n/\beta_g$, then $\tilde B$ is singular and
\[
\alpha^2\frac{\beta_g}{n}\leq\frac{\tilde\beta_{\tilde g}}{\tilde n}.
\]
\end{enumerate}
\end{prop}
 
\begin{proof}
The argument is illustrated in Figure~\ref{Fig_Clusters}. Fix $j$ with $\alpha^2\beta_j>n$. For $t\neq0$ with $|t|$ small consider the cluster $\mathrm{Cl}_j(t)=\{p_\xi:=\pi(\xi t)\mid\xi\in\mu_{e_{j-1}}\}$ and its image $\{p'_\xi:=F(p_\xi)\}\subset\tilde B$. Set $\rho:=|t|^n$. All constants below are independent of $t$. We write $u\asymp v$ if $u/v$ is bounded above and below by positive constants, and $u=o(v)$ if $u/v\to0$ as $t\to0$. The proof proceeds in five steps.
 
\medskip
\noindent\emph{Step 1: the image of the cluster is still a cluster of comparable points.} Let $\xi\neq\xi'$ be in $\mu_{e_{j-1}}$. By \ref{C2} and \eqref{Eq_HolderConstantsF},
\begin{equation}\label{Eq_ImageDistance}
\|p'_\xi-p'_{\xi'}\|\leq K_1\|p_\xi-p_{\xi'}\|^\alpha\leq K\rho^{\alpha\beta_j/n}.
\end{equation}
By \ref{C1} and the Hölder estimate for $F^{-1}$ applied to the pair $(p'_\xi,0)$, we have $c\rho\leq\|p_\xi\|=\|F^{-1}(p'_\xi)-F^{-1}(0)\|\leq K_2\|p'_\xi\|^\alpha$. Hence
\begin{equation}\label{Eq_ImageNorm}
\|p'_\xi\|\geq k\rho^{1/\alpha}\qquad\text{where }k:=(c/K_2)^{1/\alpha}.
\end{equation}
The hypothesis $\alpha^2\beta_j>n$ says precisely that $\alpha\beta_j/n>1/\alpha$. Thus $\rho^{\alpha\beta_j/n}=o(\rho^{1/\alpha})$: the mutual distances of the points $p'_\xi$ are $o$ of their norms. Set $\rho':=\|p'_1\|$. Then $\|p'_\xi\|=\rho'(1+o(1))$ for every $\xi$, and $\rho'\geq k\rho^{1/\alpha}$.
 
\medskip
\noindent\emph{Step 2: the sheets of $\tilde B$ carrying the image points.} Let $\tilde\pi : (\C,0)\to\tilde B$ be the normalization, with $\tilde x\circ\tilde\pi=\tau^{\tilde n}$, and write $p'_\xi=\tilde\pi(\tau_\xi)$. By \eqref{Eq_NormOnBranch} for $\tilde B$ we have $|\tau_\xi|^{\tilde n}\asymp\rho'$. Because $\tilde x$ is Lipschitz,
\[
\bigl|\tau_{\xi'}^{\tilde n}-\tau_\xi^{\tilde n}\bigr|=|\tilde x(p'_{\xi'})-\tilde x(p'_\xi)|\leq L\|p'_{\xi'}-p'_\xi\|=o(\rho').
\]
Hence $(\tau_{\xi'}/\tau_\xi)^{\tilde n}\to1$ as $t\to0$. Because the roots of $w^{\tilde n}=1$ are simple, for $|t|$ small there is a unique root of unity $\zeta=\zeta(\xi,\xi';t)\in\mu_{\tilde n}$ nearest to $\tau_{\xi'}/\tau_\xi$. Set $\eta=\eta(\xi,\xi';t):=\tau_{\xi'}/(\zeta\tau_\xi)-1$. Then
\begin{equation}\label{Eq_SheetDecomposition}
\tau_{\xi'}=\zeta\,\tau_\xi\,(1+\eta)\qquad\text{and}\qquad|\eta|\leq C_1\bigl|(\tau_{\xi'}/\tau_\xi)^{\tilde n}-1\bigr|\leq C_2\,\frac{\|p'_{\xi'}-p'_\xi\|}{\rho'}=o(1).
\end{equation}
In words, $p'_{\xi'}$ lies on the sheet of $\tilde B$ labelled by $\zeta$ relative to the sheet of $p'_\xi$, up to a small relative error $\eta$. We also record the following estimate. The components of $\tilde\pi$ have order at least $\tilde n$, so their derivatives satisfy $|\tilde\pi_k'(\tau)|\leq C|\tau|^{\tilde n-1}$. Therefore, for $\tau$ with $|\tau|\asymp|\tau_\xi|$ and $|\eta'|\leq|\eta|$,
\begin{equation}\label{Eq_MoveAlongSheet}
\|\tilde\pi(\tau(1+\eta'))-\tilde\pi(\tau)\|\leq C|\tau|^{\tilde n}|\eta|\leq C'\rho'|\eta|\leq C''\|p'_{\xi'}-p'_\xi\| ,
\end{equation}
where the last inequality follows from \eqref{Eq_SheetDecomposition}.
 
\medskip
\noindent\emph{Step 3: distinct points of the cluster go to distinct sheets.} We show that $\zeta\neq1$ in \eqref{Eq_SheetDecomposition} for $|t|$ small. Suppose $\zeta=1$, that is, $\tau_{\xi'}=\tau_\xi(1+\eta)$. Define
\[
\sigma : [0,1]\longrightarrow\tilde B,\qquad\sigma(s):=\tilde\pi\bigl(\tau_\xi(1+s\eta)\bigr).
\]
This is a path from $p'_\xi$ to $p'_{\xi'}$ inside $\tilde B$. By \eqref{Eq_MoveAlongSheet} and \eqref{Eq_ImageDistance} its diameter is at most $C''\|p'_{\xi'}-p'_\xi\|\leq C''K\rho^{\alpha\beta_j/n}$. Its preimage $F^{-1}\circ\sigma$ is a path in $\mathcal C$ from $p_\xi$ to $p_{\xi'}$, and by \eqref{Eq_HolderConstantsF} its diameter is at most $K_2(\operatorname{diam}\sigma)^\alpha\leq K'\rho^{\alpha^2\beta_j/n}$. On the other hand, by Lemma~\ref{Lem_ClusterSeparation} its diameter is at least $c_0\rho$. Thus $c_0\rho\leq K'\rho^{\alpha^2\beta_j/n}$, which is impossible for $\rho$ small because $\alpha^2\beta_j/n>1$. This proves the claim.
 
\medskip
\noindent\emph{Step 4: the contact orders of the sheets involved are large.} Let $\xi\neq\xi'$ and $\zeta\neq1$ be as in \eqref{Eq_SheetDecomposition}. By Lemma~\ref{Lem_SheetDistances} applied to $\tilde B$ and by \eqref{Eq_MoveAlongSheet},
\[
\|p'_{\xi'}-p'_\xi\|\geq\|\tilde\pi(\zeta\tau_\xi)-\tilde\pi(\tau_\xi)\|-\|\tilde\pi(\zeta\tau_\xi(1+\eta))-\tilde\pi(\zeta\tau_\xi)\|
\geq c|\tau_\xi|^{\beta_{\tilde B}(\zeta)}-C''\|p'_{\xi'}-p'_\xi\| .
\]
Because $|\tau_\xi|^{\tilde n}\asymp\rho'$, together with \eqref{Eq_ImageDistance} this gives
\[
c'\rho'^{\,\beta_{\tilde B}(\zeta)/\tilde n}\leq(1+C'')\|p'_{\xi'}-p'_\xi\|\leq K''\rho^{\alpha\beta_j/n}.
\]
Combined with $\rho'\geq k\rho^{1/\alpha}$ from Step~1 we obtain
\[
\rho^{\beta_{\tilde B}(\zeta)/(\alpha\tilde n)}\leq K'''\rho^{\alpha\beta_j/n}
\]
for $|t|$ small. The root of unity $\zeta$ depends on $t$, but it takes only finitely many values. Choose a sequence $t_k\to0$ along which $\zeta(\xi,\xi';t_k)$ is constant for every pair $(\xi,\xi')$. Letting $\rho\to0$ along this sequence we conclude that
\begin{equation}\label{Eq_RateInequality}
\frac{\beta_{\tilde B}(\zeta)}{\tilde n}\geq\alpha^2\frac{\beta_j}{n}
\qquad\text{for every }\zeta=\zeta(\xi,\xi';t_k)\text{ with }\xi\neq\xi'.
\end{equation}
 
\medskip
\noindent\emph{Step 5: counting.} Fix the point $p'_1$ of the image cluster. For $\xi\in\mu_{e_{j-1}}\setminus\{1\}$ denote by $\zeta_\xi:=\zeta(1,\xi;t_k)\in\mu_{\tilde n}\setminus\{1\}$ the label of the sheet of $p'_\xi$ relative to $p'_1$, as in \eqref{Eq_SheetDecomposition}. Suppose $\zeta_\xi=\zeta_{\xi'}$ for some $\xi\neq\xi'$. Then $\tau_{\xi'}/\tau_\xi=(1+\eta_{\xi'})/(1+\eta_\xi)\to1$ along the sequence, so the root of unity nearest to $\tau_{\xi'}/\tau_\xi$ is $1$ for $k$ large, that is, $\zeta(\xi,\xi';t_k)=1$. This contradicts Step~3. Hence $\xi\mapsto\zeta_\xi$ is an injective map from $\mu_{e_{j-1}}\setminus\{1\}$ into the set $\{\zeta\in\mu_{\tilde n}\setminus\{1\}\mid\beta_{\tilde B}(\zeta)/\tilde n\geq\alpha^2\beta_j/n\}$, by \eqref{Eq_RateInequality}. Because $\beta_{\tilde B}(1)=+\infty$, this set has $N_{\tilde B}(\alpha^2\beta_j/n)-1$ elements. Therefore $e_{j-1}-1\leq N_{\tilde B}(\alpha^2\beta_j/n)-1$, as claimed. (If $\tilde B$ is smooth, then $\mu_{\tilde n}=\{1\}$, so Step~2 forces $\zeta=1$ for every pair, and Step~3 yields a contradiction as soon as the cluster has two points. Because $e_{j-1}\geq2$ for $j\leq g$, a smooth $\tilde B$ is impossible. This is in accordance with the inequality, which would read $e_{j-1}\leq1$.)
 
\medskip
For (i) take $j=1$. Then $n=e_0\leq N_{\tilde B}(\alpha^2\beta_1/n)\leq\#\mu_{\tilde n}=\tilde n$. For (ii) take $j=g$. Because $e_{g-1}\geq2$, we get $N_{\tilde B}(\alpha^2\beta_g/n)\geq2$. Thus $\tilde n\geq2$, and some $\zeta\in\mu_{\tilde n}\setminus\{1\}$ satisfies $\beta_{\tilde B}(\zeta)/\tilde n\geq\alpha^2\beta_g/n$. Because $\beta_{\tilde B}(\zeta)\leq\tilde\beta_{\tilde g}$, the inequality follows.
\end{proof}
 
\begin{cor}\label{Cor_MetricMultiplicityContinuation}
Let $F : \mathcal C\to\tilde{\mathcal C}$ be a bi-$\alpha$-Hölder homeomorphism of reduced curve germs and set $M:=\max\{m(\mathcal C),m(\tilde{\mathcal C})\}$. Suppose $\alpha^2>M/(M+1)$. Then $m(F(B))=m(B)$ for every branch $B$ of $\mathcal C$; in particular $m(\mathcal C)=m(\tilde{\mathcal C})$. Moreover, a singular branch of multiplicity $n$ and last characteristic exponent $\beta_g$ is not bi-$\alpha$-Hölder homeomorphic to a smooth germ if $\alpha^2>n/\beta_g$.
\end{cor}
 
\begin{proof}
Let $B$ be a branch of $\mathcal C$ and set $\tilde B:=F(B)$. Suppose $B$ is singular, of multiplicity $n\geq2$. Then $\beta_1\geq n+1$. Because $n\leq m(\mathcal C)\leq M$ and $u\mapsto u/(u+1)$ is increasing, $n/\beta_1\leq n/(n+1)\leq M/(M+1)<\alpha^2$. Thus Proposition~\ref{Prop_ClusterLemma}~(i) gives $m(\tilde B)\geq m(B)$. If $\tilde B$ is singular, the same argument applied to $F^{-1}$ gives $m(B)\geq m(\tilde B)$. If both branches are smooth there is nothing to prove, and if exactly one of them is singular the inequality just obtained is absurd. Because the multiplicity of a reduced curve germ is the sum of the multiplicities of its branches, $m(\mathcal C)=m(\tilde{\mathcal C})$. The last assertion is Proposition~\ref{Prop_ClusterLemma}~(ii).
\end{proof}
 
\begin{rmq}\label{Rmq_ClusterThresholds}
(i) The square of $\alpha$ in Proposition~\ref{Prop_ClusterLemma} reflects the use of the Hölder estimates for both $F$ and $F^{-1}$: once in Step~1, to compare the size of the image cluster with the distance of its points to the origin, and once in Step~3, to transport a short path back to $\mathcal C$. The exponent of the proposition is at most a square root away from the truth. Consider the monomial branch $C=\{x^m=y^{\beta_1}\}$ with $\gcd(m,\beta_1)=1$, of multiplicity $m$ and single characteristic exponent $\beta_1$. Let $f=(f_1,f_2) : S^1\to C\cap S^3$ be a diffeomorphism onto the link, and consider the radial map
\[
F(tv):=\bigl(t\,f_1(v),\,t^{m/\beta_1}f_2(v)\bigr),\qquad v\in S^1,\ t\geq0 .
\]
The computation of \cite[Proposition~3.2]{Sampaio2025LogLipschitz}, which treats the case $\beta_1=m+1$ and applies verbatim, shows that $F$ is $(m/\beta_1)$-Hölder and that $F^{-1}$ is Lipschitz. Thus $F$ is a bi-$(m/\beta_1)$-Hölder homeomorphism from a disc onto $C$. By Pham and Teissier \cite{PhamTeissier1969LipFractions} (see also \cite{NeumannPichon2014}), two plane branches with the same characteristic exponents are bi-Lipschitz homeomorphic, and a generic plane projection in the sense of Definition~\ref{Def_GenericProjection} is bi-Lipschitz onto its image, because its kernel avoids the cone $C_5$ of limits of secants. Hence every branch with a single characteristic pair $(m,\beta_1)$ is bi-Lipschitz homeomorphic to the monomial branch, and the same conclusion holds for it. For these branches the threshold for being bi-$\alpha$-Hölder homeomorphic to a smooth germ therefore lies in the interval $[m/\beta_1,\sqrt{m/\beta_1}]$, and the same holds for the invariance of the multiplicity. In particular, the cusp is bi-$\frac23$-Hölder homeomorphic to a disc, although its $\alpha$-Hölder saturation is its normalization only for $\alpha\leq\frac13$ (Example~\ref{Ex_Cusp} and Theorem~\ref{Th_HolderSatSmoothness}). The radial map is not holomorphic, and non-analytic Hölder homeomorphisms can be far better than the ones produced by the saturation. It is worth noting that by Sampaio's theorem \cite[Theorem~3.1]{Sampaio2025LogLipschitz} a complex analytic set which is bi-$\alpha$-Hölder homeomorphic to a ball for \emph{all} $\alpha<1$ is smooth, and \cite[Proposition~3.2]{Sampaio2025LogLipschitz} shows that no exponent depending only on the dimension can replace ``all $\alpha$''. Theorem~\ref{Th_HomeoImpliesHolderHomeo} shows that for curves no exponent at all is forced by the topology.
 
(ii) The comparison of characteristic data under Hölder homeomorphisms of plane curves is studied by Fernandes, Sampaio and Silva \cite{FernandesSampaioSilva2018}. They prove that a bi-$\alpha$-Hölder homeomorphism distorts the contact order of two arcs by a factor in $[\alpha^2,\alpha^{-2}]$ \cite[Proposition~2.2]{FernandesSampaioSilva2018}. For two plane branches with distinct characteristic rates $a=\beta_j/n$ and $b=\tilde\beta_i/\tilde n$ they derive the obstruction
\[
\alpha^4\leq\frac{a+b}{2\max(a,b)}
\]
\cite[Lemmas~4.2 and~4.3]{FernandesSampaioSilva2018}; the paper states this with the exponent $4$ in its introduction and in the proofs, and with the exponent $2$ in the statement of Lemma~4.3. This obstruction compares the rates of the two branches pairwise at all levels, whereas Proposition~\ref{Prop_ClusterLemma} compares the sizes of the clusters of the source with the rates of the target and loses only a square root. For the multiplicity the latter is both simpler and stronger. It gives $\alpha^2\leq n/\beta_1$ when the multiplicity drops. The obstruction of \cite{FernandesSampaioSilva2018} is $\alpha^4\leq\frac{a+b}{2\max(a,b)}$ for some pair of distinct rates $a=\beta_i/n$ and $b=\tilde\beta_k/\tilde n$. Because $|a-b|\geq1/(n\tilde n)$ and $\max(a,b)\leq\max(\beta_g/n,\tilde\beta_{\tilde g}/\tilde n)$, it is at best $\alpha^4\leq1-\frac{1}{2\max(\tilde n\beta_g,\,n\tilde\beta_{\tilde g})}$, which involves the last characteristic exponents of both branches.
 
(iii) Proposition~\ref{Prop_ClusterLemma} and Corollary~\ref{Cor_MetricMultiplicityContinuation} are used in Section~\ref{Sec_HolderEquisaturation} to show that a bi-$\alpha$-Hölder trivialization of a family of curves with $\alpha^2$ above the multiplicity threshold forces the multiplicity to be constant, hence, together with the constancy of the Milnor number, Whitney equisingularity.
\end{rmq}

\section{Relative Hölder saturation and generic specialization}\label{Sec_RelativeSpecialization}

In this section we study the Hölder saturation of a family of reduced curves. Let $p : X\to Y$ be such a family and let $\alpha\in(0,1]$. We define the relative $\alpha$-Hölder saturation $X^\alpha_Y$ of $X$ over $Y$ and ask whether its fiber over $y$ is the $\alpha$-Hölder saturation of $X_y$. We prove in Theorem~\ref{Th_GenericHolderSpecialization} that this holds for $y$ in a dense Zariski-open subset of $Y$. Example~\ref{Ex_ProductFactorDifferent} shows that it can fail on the special fiber.

\begin{definition}\label{Def_FamilyOfCurves}
A \textit{family of reduced curve germs} is a flat morphism of complex analytic germs $p : (X,0)\to(Y,0)$, where $Y$ is smooth, $X$ is reduced and equidimensional, and the fibers $X_y:=p^{-1}(y)$ are reduced curves of pure dimension one.
\end{definition}

Let $\pi : \Xbar\to X$ be the normalization. Denote by $(\Xbar)_y$ the fiber of $p\circ\pi$ over $y$ and by $\overline{X_y}$ the normalization of $X_y$. The two need not coincide, and the first task is to show that they do over a dense Zariski-open subset of $Y$.

\subsection{The relative Hölder saturation}

In Section~\ref{SectionSaturations} we have defined saturations for a sequence of rings $R\to A\to B$. This relative setup is designed for families.

\begin{definition}\label{Def_RelativeHolderSaturation}
Let $p : (X,0)\to(Y,0)$ be a family of reduced curve germs and let $\alpha\in(0,1]\cap\Q$. The local ring of the \textit{relative $\alpha$-Hölder saturation} of $X$ over $Y$ is the saturation
$$
\Ocal_{X_Y^\alpha,0}:=A^{(.)_\alpha}_{B,R},
\qquad R:=\Ocal_{Y,0},\quad A:=\Ocal_{X,0},\quad B:=(\pi_*\Ocal_{\Xbar})_0 .
$$
\end{definition}

Set $W:=\Xbar\times_Y\Xbar$, let $\mathrm{pr}_1,\mathrm{pr}_2 : W\to\Xbar$ be the two projections, and let $\Ical$ be the ideal of $\Xbar\times_X\Xbar$ in $W$. It is generated by the differences $a\circ\pi\circ\mathrm{pr}_1-a\circ\pi\circ\mathrm{pr}_2$ with $a\in\Ocal_X$. As in the proof of Lemma~\ref{Lem_Coherence}, the local rings of Definition~\ref{Def_RelativeHolderSaturation} are the stalks of the coherent $\Ocal_X$-algebra
\begin{equation}\label{Eq_RelativeKernelContinuation}
\Ocal_{X_Y^\alpha}
=\Ker\Big(\pi_*\Ocal_{\Xbar}\xrightarrow{\ \Delta\ }
(\pi\circ\mathrm{pr}_1)_*(\Ocal_{W}/\Ical_\alpha)\Big),
\qquad\Delta(h)=h\otimes1-1\otimes h .
\end{equation}
Denote by $\pi^\alpha : X_Y^\alpha\to X$ the corresponding finite morphism. The inclusions $A\subseteq A^{(.)_\alpha}_{B,R}\subseteq B$ show that $X_Y^\alpha$ has normalization $\Xbar$. Because $\Ical_\alpha\subseteq\sqrt{\Ical}$, every section of $\Ocal_{X_Y^\alpha}$ takes the same value at any two points of $\Xbar$ with the same image in $X$. Thus $\pi^\alpha$ is bijective, and the sections of $\Ocal_{X_Y^\alpha}$ are continuous meromorphic functions on $X$. For $\alpha=1$ we write $X_Y^s$ for the relative Lipschitz saturation.

The following proposition extends the characterization of the relative Lipschitz saturation given by Pham and Teissier \cite[Theorem~3.1]{PhamTeissier1969LipFractions}. It says that the relative saturation consists of the functions that are $\alpha$-Hölder on the fibers with a constant independent of the fiber.

\begin{prop}\label{Prop_RelativeHolderEstimate}
Let $p : (X,0)\to(Y,0)$ be a family of reduced curve germs. Choose an embedding $X\subset\C^N\times Y$ with $p(x,y)=y$. Let $\alpha\in(0,1]\cap\Q$ and let $f\in\overline{\Ocal_{X,0}}$. Then $f\in\Ocal_{X_Y^\alpha,0}$ if and only if there exist a neighborhood $U$ of $0$ in $X$ and a constant $C>0$ such that
$$
|f(x,y)-f(x',y)|\leq C\|x-x'\|^\alpha
$$
for all $(x,y),(x',y)\in U$ lying in the same fiber of $p$.
\end{prop}

\begin{proof}
The pullback $f\circ\pi$ is holomorphic on $\Xbar$. By \eqref{Eq_RelativeKernelContinuation}, $f\in\Ocal_{X_Y^\alpha,0}$ if and only if the difference of the two pullbacks of $f\circ\pi$ to $W$ belongs to $\Ical_\alpha$ at every point of $W$ above $(0,0)$. Because the $Y$-coordinates agree on $W$, the ideal $\Ical$ is generated by the pullbacks of $x_1-x'_1,\dots,x_N-x'_N$. Lemma~\ref{Lem_InequalityPowerIdeal} therefore gives the stated inequality near each point above $(0,0)$. There are finitely many such points, so after shrinking $U$ we can take the same constant $C$ for all of them.
\end{proof}

\subsection{Generic specialization of the relative Hölder saturation}

Let $p : (X,x_0)\to(Y,0)$ be a family of reduced curve germs. We prove that, for $y$ in a dense Zariski-open subset of $Y$, the restriction of $X_Y^\alpha$ to $X_y$ is the $\alpha$-Hölder saturation of $X_y$. We need three preliminary results, on generic flatness, generic simultaneous normalization and rational powers.

\begin{lemma}[Generic flatness]\label{Lem_GenericFlatness}
Let $f : X\to Y$ be a morphism of complex analytic spaces, with $X$ reduced and $Y$ smooth and connected. Let $\Fcal$ be a coherent $\Ocal_X$-module and set
$$
\Sigma_{\Fcal}:=\{x\in X\mid\Fcal_x\text{ is not flat over }\Ocal_{Y,f(x)}\}.
$$
Suppose that $f|_{\Sigma_{\Fcal}} : \Sigma_{\Fcal}\to Y$ is proper. Then there exists a dense Zariski-open subset $U\subseteq Y$ such that $\Fcal|_{f^{-1}(U)}$ is flat over $U$.
\end{lemma}

\begin{proof}
By Frisch's generic flatness theorem \cite{Frisch1967PlatitudeMorphismes}, $\Sigma_{\Fcal}$ is a closed analytic subset of $X$ and $f(\Sigma_{\Fcal})$ has empty interior in $Y$; see \cite[Remark~2.5]{AdamusBierstoneMilman2012}. Because $f|_{\Sigma_{\Fcal}}$ is proper, Remmert's proper mapping theorem shows that $E:=f(\Sigma_{\Fcal})$ is a closed analytic subset of $Y$. Set $U:=Y\setminus E$. Then $U$ is dense and Zariski-open, and $\Fcal_x$ is flat over $\Ocal_{Y,f(x)}$ for every $x\in f^{-1}(U)$.
\end{proof}

\begin{definition}\label{Def_RelativeSingularLocus}
Let $f : X\to Y$ be a morphism of complex analytic spaces. The \textit{relative singular locus} of $f$ is
$$
\sing(f):=\{x\in X\mid X_{f(x)}\text{ is singular at }x\}=\bigcup_{y\in Y}\sing(X_y).
$$
If $f$ is flat, then $\sing(f)$ is the locus where $f$ is not smooth.
\end{definition}

\begin{lemma}[Generic simultaneous normalization]\label{Lem_GenericNormalizationProperSingularLocus}
Let $f : X\to Y$ be a flat morphism of complex analytic spaces, with $Y$ smooth and connected and $X$ reduced and equidimensional. Suppose that the fibers are reduced and pure $d$-dimensional and that $f|_{\sing(f)}$ is proper. Let $\pi : \Xbar\to X$ be the normalization. Then there exists a dense Zariski-open subset $U\subseteq Y$ such that $f\circ\pi$ is flat over $U$ and, for every $y\in U$,
$$
\pi_y : (\Xbar)_y\longrightarrow X_y
$$
is the normalization of $X_y$. In particular $(\Xbar)_y\simeq\overline{X_y}$.
\end{lemma}

\begin{proof}
We say that $f\circ\pi$ is flat with normal fiber at a point $\overline x\in\Xbar$ if $\Ocal_{\Xbar,\overline x}$ is flat over $\Ocal_{Y,y}$, where $y=f\circ\pi(\overline x)$, and the local ring $\Ocal_{(\Xbar)_y,\overline x}=\Ocal_{\Xbar,\overline x}/\m_y\Ocal_{\Xbar,\overline x}$ is normal. Let $\Sigma\subseteq\Xbar$ be the locus where $f\circ\pi$ fails to be flat with normal fiber. By \cite[Theorem~I.1.100]{GreuelLossenShustin2007}, $\Sigma$ is a closed analytic subset of $\Xbar$. Set $E:=(f\circ\pi)(\Sigma)$. We show that $E$ is a nowhere dense closed analytic subset of $Y$.

We first show that $(f\circ\pi)|_\Sigma$ is proper. Outside $\sing(f)$ the morphism $f$ is smooth, so for every $x\in X\setminus\sing(f)$ the germ $(X,x)$ is isomorphic to $(Y,f(x))\times(\C^d,0)$. Because $Y$ is smooth, $X$ is smooth at $x$. Hence $\pi$ is an isomorphism over $X\setminus\sing(f)$ and $f\circ\pi$ is smooth there. Consequently $\Sigma\subseteq\pi^{-1}(\sing(f))$. Because $\pi$ is finite and $f|_{\sing(f)}$ is proper, $(f\circ\pi)|_\Sigma$ is proper.

Now \cite[Theorem~6.2(2)(ii)]{Greuel2017Equisingular}, applied to $f\circ\pi$, whose source is normal and whose base is smooth, shows that $E$ has empty interior. By Remmert's proper mapping theorem $E$ is closed analytic. Set $U:=Y\setminus E$. Then $U$ is dense and Zariski-open, and $f\circ\pi$ is flat with normal fibers over $U$.

It remains to show that these fibers normalize the $X_y$. Fix $y\in U$. The morphism $\pi_y$ is finite and is an isomorphism over the smooth locus of $X_y$, since at such points $X$ is smooth. Because normalization preserves dimension, $\Xbar$ is pure of dimension $d+\dim Y$, so by flatness $(\Xbar)_y$ is pure $d$-dimensional. Because $\pi_y$ is finite, each irreducible component of $(\Xbar)_y$ has a $d$-dimensional image, which cannot lie in $\sing(X_y)$. Thus $\pi_y$ is finite and bimeromorphic with normal source, hence it is the normalization of $X_y$.
\end{proof}

\begin{rmq}\label{Rem_GenericNormalityAlongProperSubset}
The same argument gives the following statement, which we use below. Let $h : V\to Y$ be a morphism from a normal complex space to a smooth connected base and let $P\subseteq V$ be a closed analytic subset proper over $Y$. Then there exists a dense Zariski-open subset $U\subseteq Y$ such that, for every $y\in U$, $h$ is flat with normal fiber at every point of $P\cap V_y$. Indeed, let $\Sigma$ be the locus where $h$ fails to be flat with normal fiber. By \cite[Theorem~I.1.100]{GreuelLossenShustin2007} and \cite[Theorem~6.2(2)(ii)]{Greuel2017Equisingular}, $\Sigma$ is closed analytic and $h(\Sigma)$ has empty interior. Because $P$ is proper over $Y$, $h(P\cap\Sigma)$ is closed analytic by Remmert's theorem, and its complement is the required subset.
\end{rmq}

\begin{lemma}[Generic normalization for families of curves]\label{Lem_GenericNormalization}
Let $p : (X,x_0)\to(Y,0)$ be a family of reduced curve germs. After choosing sufficiently small representatives, $p|_{\sing(p)}$ is finite, and there exists a dense Zariski-open subset $U_0\subseteq Y$ such that $p\circ\pi$ is smooth over $U_0$ and
$$
\pi_y : (\Xbar)_y\longrightarrow X_y
$$
is the normalization of $X_y$ for every $y\in U_0$. Consequently $W$ is smooth over $U_0$ and $W_y=\overline{X_y}\times\overline{X_y}$ for every $y\in U_0$.
\end{lemma}

\begin{proof}
After choosing sufficiently small representatives, every fiber is a reduced curve \cite[Remark~I.1.100.1]{GreuelLossenShustin2007} and $X_0\setminus\{x_0\}$ is smooth. Because $p$ is flat, $\sing(p)\cap X_0=\sing(X_0)\subseteq\{x_0\}$. Thus the central fiber of $p|_{\sing(p)} : \sing(p)\to Y$ is at most one point, and by the local finiteness theorem \cite[Theorem~I.1.66]{GreuelLossenShustin2007} we can choose representatives for which this morphism is finite, hence proper. Lemma~\ref{Lem_GenericNormalizationProperSingularLocus} gives $U_0$. Over $U_0$ the morphism $p\circ\pi$ is flat with normal, hence smooth, curve fibers, so it is smooth. The last assertion follows, since $W=\Xbar\times_Y\Xbar$.
\end{proof}

From now on we fix representatives as in Lemma~\ref{Lem_GenericNormalization}. Set $\Ical_y:=\Ical\Ocal_{W_y}$ for $y\in U_0$. We will compare the restriction of $\Ical_{a/b}$ to $W_y$ with $(\Ical_y)_{a/b}$. The following description of rational powers allows us to treat all numerators $a$ with the same denominator $b$ at once.

\begin{lemma}[Graded description of rational powers]\label{Lem_GradedRationalPowers}
Let $R$ be a normal analytic local domain, let $I\subseteq R$ be a nonzero ideal and let $b\in\N^*$. Set $\Rcal_b(I):=R[It,It^b]\subseteq R[t]$. Then the normalization of $\Rcal_b(I)$ in its fraction field is
$$
\overline{\Rcal_b(I)}=R\oplus\bigoplus_{a\geq1}I_{a/b}\,t^a .
$$
\end{lemma}

\begin{proof}
Choose $0\neq c\in I$. Because $t=(ct)/c$, the rings $\Rcal_b(I)$ and $R[t]$ have the same fraction field. Hence $\overline{\Rcal_b(I)}\subseteq R[t]$ by normality of $R[t]$, and this normalization is graded by \cite[Theorem~2.3.2]{SwansonHuneke2006IntegralClosure}. Moreover $\Rcal_b(I)$ is integral over $R[It^b]$, because $(ct)^b=c^{b-1}(ct^b)$ for every $c\in I$. Let $h\in R$ and $a\in\N^*$. Set $u:=t^b$. Using $(ht^a)^b=h^bu^a$ we get
$$
ht^a\in\overline{\Rcal_b(I)}
\iff h^bu^a\text{ is integral over }R[Iu]
\iff h^b\in\overline{I^a}
\iff h\in I_{a/b},
$$
where the second equivalence is \cite[Proposition~5.2.1]{SwansonHuneke2006IntegralClosure}. Because the degree zero part of $\overline{\Rcal_b(I)}$ is $R$, the lemma follows.
\end{proof}

\begin{lemma}[Generic specialization of rational powers]\label{Lem_RationalPowersGeneralFibres}
Let $b\in\N^*$. There exists a dense Zariski-open subset $U_b\subseteq U_0$ such that, for every $y\in U_b$ and every $a\in\N^*$, the natural map
$$
\Ical_{a/b}\otimes_{\Ocal_Y}k(y)\longrightarrow\Ocal_{W_y}
$$
is injective with image $(\Ical_y)_{a/b}$.
\end{lemma}

\begin{proof}
Fix $b\in\N^*$, let $y\in U_0$ and let $w=(u,v)\in W_y$. We first treat the points not lying over a singular point of $X_y$. Suppose $\pi(u)\neq\pi(v)$. Then $\Ical_w=\Ocal_{W,w}$ and the assertion is clear. Suppose $\pi(u)=\pi(v)$ is a smooth point of $X_y$. Then $u=v$ and $\pi$ is an isomorphism near $u$. Choose a relative parameter $s$ there and denote by $s_1,s_2$ its pullbacks to the two factors of $W$. Then $\Ical_w=(s_1-s_2)$. Because $s_1-s_2$ is prime, a germ $h$ lies in $(\Ical_{a/b})_w$ if and only if $(s_1-s_2)^a$ divides $h^b$, so
$$
(\Ical_{a/b})_w=\big((s_1-s_2)^{\lceil a/b\rceil}\big).
$$
The same formula computes the rational power on $W_y$, because the restriction of $s_1-s_2$ is again a local parameter. Thus restriction takes a generator to a generator, and multiplication by the restricted generator is injective in $\Ocal_{W_y,w}$. This proves the assertion at $w$ for every $y\in U_0$.

It remains to treat the set
$$
Z:=\{(u,v)\in W\mid\pi(u)=\pi(v)\in\sing(p)\},
$$
which is finite over $Y$ because $\Xbar\times_X\Xbar\to X$ and $\sing(p)\to Y$ are finite. To obtain one open subset of $Y$ for all numerators $a$, we use the graded algebra
$$
\Rcal_b(\Ical):=\Ocal_W[\Ical t,\Ical t^b].
$$
For $b=1$ this is the Rees algebra of $\Ical$. Over $U_0$ the stalks of $\Ocal_W$ are regular, and $\Ical$ is nonzero at every point, so Lemma~\ref{Lem_GradedRationalPowers} gives
$$
\overline{\Rcal_b(\Ical)}=\Ocal_W\oplus\bigoplus_{a\geq1}\Ical_{a/b}\,t^a
$$
over $U_0$. The same formula holds on $W_y$ for $\Rcal_b(\Ical_y):=\Ocal_{W_y}[\Ical_yt,\Ical_yt^b]$, because $V(\Ical_y)$ is finite over the curve $X_y$, so $\Ical_y$ is nonzero at every point of the smooth surface $W_y$.

An equation of integral dependence over $\Rcal_b(\Ical)$ restricts to an equation of integral dependence over $\Rcal_b(\Ical_y)$. Restriction of coefficients therefore induces a graded morphism
\begin{equation}\label{Eq_RestrictionNormalizedGradedAlgebras}
\overline{\Rcal_b(\Ical)}\otimes_{\Ocal_Y}k(y)
\longrightarrow
\overline{\Rcal_b(\Ical_y)} .
\end{equation}
On the left we first normalize on $W$ and then restrict to $W_y$; on the right we first restrict $\Ical$ to $W_y$ and then normalize. In degree $a$ the morphism sends $ht^a\otimes1$ to $(h|_{W_y})t^a$, so the assertion of the lemma at $w$ says precisely that \eqref{Eq_RestrictionNormalizedGradedAlgebras} is an isomorphism at $w$ in every degree.

Fix $y\in U_0$ and $w\in Z\cap W_y$. Denote by $S$ and $T$ the stalks at $w$ of the source and of the target of \eqref{Eq_RestrictionNormalizedGradedAlgebras}. The degree zero part of $S$ is the local ring $\Ocal_{W_y,w}$. Let $\n$ be the homogeneous maximal ideal of $S$, generated by the maximal ideal of $\Ocal_{W_y,w}$ and by the elements of positive degree. We first show that it suffices to prove that $S\to T$ becomes an isomorphism after localization at $\n$. Indeed, the kernel and the cokernel of $S\to T$ are graded $S$-modules. The annihilator of a nonzero homogeneous element of either module is a proper homogeneous ideal, hence it is contained in $\n$, and such an element remains nonzero after localization at $\n$.

The morphism $S\to T$ is finite, because its image contains $\Rcal_b(\Ical_y)_w$, whose normalization $T$ is finite over it by excellence of $\Ocal_{W_y,w}$. It is an isomorphism away from $V(\Ical_y)$. To see this, choose $g\in\Ical_w$ whose restriction $g_y$ is nonzero. After inverting $g$ we have $t=(gt)/g$, so $\Rcal_b(\Ical)_w[1/g]=\Ocal_{W,w}[1/g][t]$, which is normal, and the same computation holds on the fiber. Hence both $S[1/g_y]$ and $T[1/g_y]$ are identified with $\Ocal_{W_y,w}[1/g_y][t]$.

Suppose now that $S_\n$ is a normal domain. The localization of $T$ at the image of $\n$ is a domain, because $T$ is a subring of $\Ocal_{W_y,w}[t]$ and every element of $S\setminus\n$ maps to a polynomial with invertible constant term. The element $g_y$ is nonzero in the localized target, hence also in the localized source. After inverting $g_y$ the two localized rings agree. Thus the localized morphism is injective, and the two domains have the same fraction field. Because the morphism is finite and $S_\n$ is integrally closed in this field, it is an isomorphism. It remains to show that $S_\n$ is a normal domain for every $w\in Z\cap W_y$ and every $y$ in a dense Zariski-open subset of $Y$.

Let $C$ be the reduction of the analytic cone $\operatorname{Specan}\Rcal_b(\Ical)$ over $W$, let $\theta : C\to W$ be its projection, and let $\overline C\to C$ be the normalization. The zero section of $C$ is defined by the vanishing of the elements of positive degree, and $\theta$ maps it isomorphically onto $W_{\mathrm{red}}$. Let $w\in W$. Because the degree zero part of $\overline{\Rcal_b(\Ical)}_w$ is the local ring $\Ocal_{W,w}$, the algebra $\overline{\Rcal_b(\Ical)}_w$ has exactly one homogeneous maximal ideal, so there is exactly one point of the zero section of $\overline C$ over $w$. Denote it by $\tilde w$; its analytic local ring is the analytic localization of $\overline{\Rcal_b(\Ical)}_w$ at this homogeneous maximal ideal. The points $\tilde w$ with $w\in Z$ form a closed analytic subset $P$ of $\overline C$ which is finite over $Y$, because $Z$ is. Apply Remark~\ref{Rem_GenericNormalityAlongProperSubset} to $\overline C\to Y$ and to $P$. We obtain a dense Zariski-open subset $U_b\subseteq U_0$ such that, for every $y\in U_b$ and every $w\in Z\cap W_y$, the fiber $(\overline C)_y$ is normal at $\tilde w$.

Let $y\in U_b$ and $w\in Z\cap W_y$. The analytic local ring of $(\overline C)_y$ at $\tilde w$ and the ring $S_\n$ have the same completion. Indeed, before passing to the fiber, the polynomial presentation of the cone over $\Ocal_{W,w}$ and its convergent presentation have the same completion at the point $\tilde w$, and analytification commutes with finite normalization there, as in \cite[proof of Lemma~5.7.4]{ChiangHsiehLipman2006}. On both sides, passing to the fiber is taking the quotient by the maximal ideal of $\Ocal_{Y,y}$. Because the analytic local ring of $(\overline C)_y$ at $\tilde w$ is normal and excellent, its completion is normal \cite[Tag~0C23]{StacksProject}. Normality descends along the faithfully flat completion map, so $S_\n$ is a normal domain.

Thus \eqref{Eq_RestrictionNormalizedGradedAlgebras} is an isomorphism at every $w\in Z\cap W_y$ for $y\in U_b$, and its degree $a$ part is
$$
\Ical_{a/b}\otimes_{\Ocal_Y}k(y)\xrightarrow{\ \sim\ }(\Ical_y)_{a/b}\subseteq\Ocal_{W_y}
$$
for every $a\geq1$. Together with the first part of the proof this proves the lemma.
\end{proof}

\begin{theorem}\label{Th_GenericRelativeHolderSpecialization}\label{Th_GenericHolderSpecialization}
Let $p : (X,x_0)\to(Y,0)$ be a family of reduced curve germs and let $\alpha\in(0,1]$. After shrinking the representatives, there exists a dense Zariski-open subset $U\subseteq Y$ such that, for every $y\in U$, $\pi_y : (\Xbar)_y\to X_y$ is the normalization of $X_y$ and the natural map
$$
\Ocal_{X_Y^\alpha}\otimes_{\Ocal_Y}k(y)
\longrightarrow(\pi_y)_*\Ocal_{\overline{X_y}}
$$
is injective with image $\Ocal_{X_y^\alpha}$. In particular $(X_Y^\alpha)_y\simeq X_y^\alpha$ canonically as finite spaces over $X_y$. For finitely many exponents $\alpha_1,\dots,\alpha_k$, the open subset $U$ can be chosen so that the conclusion holds for all of them simultaneously.
\end{theorem}

\begin{proof}
By \eqref{Eq_RelativeKernelContinuation}, $\Ocal_{X_Y^\alpha}$ is the kernel of the morphism of coherent $\Ocal_X$-modules
$$
\Delta : \pi_*\Ocal_{\Xbar}\longrightarrow(\pi\circ\mathrm{pr}_1)_*(\Ocal_W/\Ical_\alpha),
\qquad h\longmapsto h\otimes1-1\otimes h .
$$
We first show that taking this kernel commutes with restriction to a general fiber. On $X\setminus\sing(p)$ the normalization is an isomorphism. Along the diagonal of $W$ we have $\Ical=(s_1-s_2)$ in relative coordinates, and $\Ical_\alpha=(s_1-s_2)$ because $s_1-s_2$ is prime and $0<\alpha\leq1$; off the diagonal the quotient $\Ocal_W/\Ical_\alpha$ vanishes. Thus on $X\setminus\sing(p)$ the target of $\Delta$ is $\Ocal_X$ and $\Delta=0$. Because $p$ is flat, the nonflat loci over $Y$ of the image and of the cokernel of $\Delta$ are contained in $\sing(p)$, which is finite over $Y$. By Lemma~\ref{Lem_GenericFlatness} there is a dense Zariski-open subset $U\subseteq U_0$ over which the image and the cokernel of $\Delta$ are flat. Then, for every $y\in U$,
$$
\Ocal_{X_Y^\alpha}\otimes_{\Ocal_Y}k(y)
=\Ker\big(\Delta\otimes_{\Ocal_Y}k(y)\big).
$$
Indeed, write $0\to\Ocal_{X^\alpha_Y}\to\pi_*\Ocal_{\Xbar}\to\mathcal E\to0$ and $0\to\mathcal E\to\mathcal T\to\mathcal T/\mathcal E\to0$, where $\mathcal E$ is the image and $\mathcal T$ the target of $\Delta$. Flatness of $\mathcal E$ and of $\mathcal T/\mathcal E$ gives the vanishing of the relevant $\mathrm{Tor}_1$, so both sequences stay exact after tensoring with $k(y)$.

Next we identify the source. Direct images under finite morphisms commute with restriction to a fiber, so
$$
(\pi_*\Ocal_{\Xbar})\otimes_{\Ocal_Y}k(y)\simeq(\pi_y)_*\Ocal_{(\Xbar)_y},
$$
and $\pi_y$ is the normalization of $X_y$ for $y\in U_0$ by Lemma~\ref{Lem_GenericNormalization}.

Next we identify the target. Because $\Ical\subseteq\Ical_\alpha$, the quotient $\Ocal_W/\Ical_\alpha$ is a sheaf on $\Xbar\times_X\Xbar$, which is finite over $X$, so its direct image also commutes with restriction to a fiber. By right exactness of the tensor product,
$$
(\Ocal_W/\Ical_\alpha)\otimes_{\Ocal_Y}k(y)\simeq\Ocal_{W_y}/\Ical_\alpha\Ocal_{W_y}.
$$
Suppose $\alpha=a/b$ is rational. Replacing $U$ by $U\cap U_b$, Lemma~\ref{Lem_RationalPowersGeneralFibres} gives $\Ical_\alpha\Ocal_{W_y}=(\Ical_y)_\alpha$ for every $y\in U$. Suppose $\alpha$ is irrational. By Proposition~\ref{Prop_FinitenessOfRationalIdeals}, at each of the finitely many points of $\pi^{-1}(x_0)\times\pi^{-1}(x_0)$ the exponents in $(0,1]$ at which the rational powers of $\Ical$ change belong to a finite set of rational numbers. Choose rational numbers $0<\beta<\alpha<\gamma\leq1$ with $\Ical_\beta=\Ical_\gamma$ at all these stalks. By coherence this equality holds on a neighborhood of $\pi^{-1}(x_0)\times\pi^{-1}(x_0)$, hence on $W$ after shrinking the representatives, since $W\to X\times_YX$ is finite. Write $\beta$ and $\gamma$ with a common denominator $b$. For $y\in U\cap U_b$ the rational case gives $(\Ical_y)_\gamma=\Ical_\gamma\Ocal_{W_y}=\Ical_\beta\Ocal_{W_y}=(\Ical_y)_\beta$. By monotonicity of rational powers, both $\Ical_\alpha\Ocal_{W_y}$ and $(\Ical_y)_\alpha$ lie between these two equal ideals, so $\Ical_\alpha\Ocal_{W_y}=(\Ical_y)_\alpha$ in this case as well.

Finally, $\Ical_y$ is the diagonal ideal of $X_y$. Indeed, $\Ical$ is locally generated by the differences $a\otimes1-1\otimes a$ with $a\in\Ocal_X$, every local function on $X_y$ lifts locally to $X$, and the restrictions of these generators are exactly the differences defining the diagonal ideal of $X_y$ on $\overline{X_y}\times\overline{X_y}$. Under the identifications above, $\Delta\otimes_{\Ocal_Y}k(y)$ sends $h$ to $h\otimes1-1\otimes h$ modulo $(\Ical_y)_\alpha$, that is, it is the morphism defining the $\alpha$-Hölder saturation of $X_y$. Its kernel is $\Ocal_{X_y^\alpha}$, which proves the assertion for one exponent. For finitely many exponents, intersect the corresponding open subsets.
\end{proof}

\begin{rmq}\label{Rmq_SpecialFiberNotSaturation}
Theorem~\ref{Th_GenericHolderSpecialization} says nothing about the special fiber. The restriction $\Ocal_{X^\alpha_Y}\otimes_{\Ocal_Y}k(0)\to(\pi_0)_*\Ocal_{(\Xbar)_0}$ need not be injective, $(\Xbar)_0$ need not be the normalization of $X_0$, and even when both hold the image can be strictly smaller than $\Ocal_{X_0^\alpha}$. Example~\ref{Ex_ProductFactorDifferent} illustrates the last phenomenon. These failures are the subject of Section~\ref{Sec_HolderEquisaturation}.
\end{rmq}

\section{Hölder equisaturation and equisingularity}\label{Sec_HolderEquisaturation}

In this section we apply the relative Hölder saturation to equisingularity problems for families of curves. We first review the setup of \cite{BengusGaffneyRangachev2026}, the invariants of curve singularities that enter, and the three classical notions of equisingularity: Whitney, topological and Lipschitz. We then define $\alpha$-Hölder equisaturation, prove that it implies $\alpha$-Hölder equisingularity, and compare it with the classical notions through the numerical invariants of the special fiber.

\textbf{Family setup.} As in Section~\ref{Sec_RelativeSpecialization}, a family $p : X\to Y$ is a flat morphism of complex analytic germs $p : (X,0)\to(Y,0)$, where $Y$ is smooth and connected of dimension $k$, $X$ is reduced and equidimensional, and the fibers $X_y$ are reduced curves. We say that $p$ is \textit{centered} if there exists a section $\sigma : Y\to X$ with $\sing(p)=\sigma(Y)$, where $\sing(p)$ is the relative singular locus of Definition~\ref{Def_RelativeSingularLocus}. In this case we identify $Y$ with $\sigma(Y)$ and choose an embedding $(X,0)\subset(\C^n\times\C^k,0)$ with $\sigma(y)=(0,y)$; then $p$ is the projection to the second factor and $\sigma(y)$ is the only singular point of $X_y$ near $0$.

\textbf{Invariants.} Let $(X,x)$ be a reduced curve germ. Denote by $m(X,x)$ its multiplicity, by $r(X,x)$ its number of branches, by $\delta(X,x):=\dim_\C\Ocal_{\overline X,x}/\Ocal_{X,x}$ its delta invariant, and by $\mu(X,x)$ its Milnor number in the sense of Buchweitz and Greuel \cite{BuchweitzGreuel1980ConstancyOfMilnorNbr}. Recall from \cite[Proposition~1.2.1]{BuchweitzGreuel1980ConstancyOfMilnorNbr} that
$$
\mu(X,x)=2\delta(X,x)-r(X,x)+1 .
$$
For a reduced curve $X$ with finitely many singular points, the total delta invariant and the total Milnor number are
$$
\delta(X):=\sum_{x\in\sing(X)}\delta(X,x),\qquad\mu(X):=\sum_{x\in\sing(X)}\mu(X,x).
$$
For a family $p : X\to Y$ we apply these to a small representative of the fiber $X_y$. If $p$ is centered, then $\delta(X_y)=\delta(X_y,\sigma(y))$ and $\mu(X_y)=\mu(X_y,\sigma(y))$. Set
$$
m_0:=m(X_0,0),\qquad r_0:=r(X_0,0),\qquad\delta_0:=\delta(X_0,0),\qquad\mu_0:=\mu(X_0,0).
$$
The function $y\mapsto\mu(X_y)$ is upper semicontinuous \cite[Theorem~6.7]{BengusGaffneyRangachev2026}, and so is $y\mapsto m(X_y,\sigma(y))$ for a centered family, by upper semicontinuity of the multiplicity.

\textbf{Whitney equisingularity.} Suppose $p$ is centered and embedded as above. Let $\{x_i\}$ be a sequence of smooth points of $X$ and $\{y_i\}$ a sequence of points of $\sigma(Y)$, both converging to $0$. Consider the secant lines $\overline{x_iy_i}$ and the tangent spaces $T_{x_i}X$. After passing to subsequences they converge to a line $l$ and to a linear space $T$. We say that $p$ is \textit{Whitney equisingular} if $l\subseteq T$ for every such pair of sequences, that is, if the pair $(X\setminus\sigma(Y),\sigma(Y))$ satisfies the Whitney condition (b) at $0$. By \cite[Theorem~1.3]{BengusGaffneyRangachev2026}, a centered family is Whitney equisingular if and only if $\mu(X_y,\sigma(y))$ and $m(X_y,\sigma(y))$ are constant, if and only if $p$ admits a strong simultaneous resolution, that is, the fibers $(\Xbar)_y$ are smooth and $\pi^{-1}(\sigma(Y))$, with its natural analytic structure, is locally analytically trivial over $Y$. For $k=1$ this is due to Briançon, Galligo and Granger \cite{BrianconGallicoGranger1980DeformationsEqu}.

\textbf{Topological equisingularity.} We say that $p$ is \textit{topologically equisingular} if there exists a homeomorphism of germs $\Phi : (X,0)\to(X_0,0)\times(Y,0)$ with $p=\mathrm{pr}_Y\circ\Phi$. By \cite[Theorem~1.4]{BengusGaffneyRangachev2026}, a centered family is topologically equisingular if and only if $\mu(X_y,\sigma(y))$ is constant, and for a centered family the constancy of $\mu$ is equivalent to the constancy of $\delta$ and $r$. For $k=1$ this is due to Buchweitz and Greuel \cite{BuchweitzGreuel1980ConstancyOfMilnorNbr}. Because Whitney equisingularity implies the constancy of $\mu$, it implies topological equisingularity; the converse fails, since the multiplicity can drop in a $\mu$-constant family.

\textbf{Hölder equisingularity.} Let $\alpha\in(0,1]$. We say that $p$ is \textit{$\alpha$-Hölder equisingular} if there exists a homeomorphism of germs $\Phi : (X,0)\to(X_0,0)\times(Y,0)$ with $p=\mathrm{pr}_Y\circ\Phi$ such that $\Phi$ and $\Phi^{-1}$ are $\alpha$-Hölder for the outer metrics of representatives. An $\alpha$-Hölder equisingular family is topologically equisingular, and it is $\beta$-Hölder equisingular for every $\beta\leq\alpha$ by Lemma~\ref{Lem_HolderElementary}~(iii). For $\alpha=1$ we say that $p$ is \textit{Lipschitz equisingular}; this is the notion studied by Mostowski \cite{Mostowski1985}, Parusi\'nski \cite{Parusinski1994}, Parusi\'nski and P\u{a}unescu \cite{ParusinskiPaunescu2017} and Neumann and Pichon \cite{NeumannPichon2014,NeumannPichon2022}. Hölder equisingularity is a metric condition on the total space. The notion of $\alpha$-Hölder equisaturation introduced below is an algebraic condition, the triviality of the relative saturation $X^\alpha_Y$ of Section~\ref{Sec_RelativeSpecialization}, and the first thing we show is that it implies $\alpha$-Hölder equisingularity (Proposition~\ref{Prop_ProductHolderTriviality}). We then relate it to the three classical notions: it implies the constancy of $\mu$, hence topological equisingularity, for every $\alpha$ (Corollary~\ref{Cor_MuSomeAlpha}); it implies Lipschitz equisaturation for $\alpha$ above a bound depending on $\mu_0$ and $m_0$ (Subsection~\ref{Subsec_LipschitzEquisaturation}); it is implied by the constancy of $\mu$ for $\alpha$ below an explicit bound depending on $\delta_0$ (Theorem~\ref{Th_UniformSmallExponent}); and it implies Whitney equisingularity for $\alpha$ above a bound depending on $m_0$ and $r_0$ (Theorem~\ref{Th_EquisaturationWhitney}).

\subsection{Product structures and their fibers}

We introduce a new notion of equisingularity that we call "$\alpha$-equisaturation".

\begin{definition}\label{Def_AlphaEquisaturation}
A family $p:X\to Y$ is said to be $\alpha$-equisaturated if there exists a germ of curve $(F,0)$ and an isomorphism of germs $\varphi:X_Y^\alpha\longrightarrow F \times Y$  preserving the base, i.e. such that the following diagram commutes :
\[\begin{tikzcd}
	{X_Y^\alpha} & X & Y \\
	{F\times Y}
	\arrow["{\pi^\alpha}", from=1-1, to=1-2]
	\arrow["\varphi", from=1-1, to=2-1]
	\arrow["p", from=1-2, to=1-3]
	\arrow["{\text{pr}_Y}"', from=2-1, to=1-3]
\end{tikzcd}\]
where $\pi^\alpha$ is the relative saturation morphism and pr$_Y$ is the projection to $Y$. For $\alpha=1$, we say that $p$ is Lipschitz-equisaturated.
\end{definition}

A family having fibers with isomorphic $\alpha$-Hölder saturations is not necessarily $\alpha$-equisaturated. See the $1/5$-equisaturation in Example \ref{Ex_RelativeSeminormalExponent} for instance. Note also that, when $\alpha<1$, the factor $F$ in the definition of $\alpha$-equisaturation need not be isomorphic to $X_0^\alpha$. Moreover, if a family is $\alpha$-equisaturated, it does not imply that the $\alpha$-Hölder saturations of $X_\lambda$ and $X_0$ are isomorphic. Let us give an example of such a family. 

\begin{ex}\label{Ex_ProductFactorDifferent}
Let $(Y,0)=(\C,0)$ and consider the centered $\mu$-constant family
$$
X_\lambda:t\longmapsto(t^3+\lambda t^2,t^4,t^5).
$$
One can easily see that the characteristic sequences of $X_0$ is $(3,4)$ and, after a change of variable in order to have a good uniformizer, one can verify that the characteristic sequence of $X_\lambda$ is $(2,5)$ when $\lambda\neq 0$. By Proposition \ref{Prop_ComputeSemigroupOfTheHolderSat} and by using the arc criterion, one can check that 
$$
\Ocal_{X_Y^{1/4}}=\C\{\lambda,t^2,t^3\} \text{, }\quad\Ocal_{X_0^{1/4}}=\C\{t\},\quad \text{and}\quad\Ocal_{X_\lambda^{1/4}}=\C\{t^2,t^3\}\text{ for }\lambda\neq 0
$$
Hence, the relative $1/4$-Hölder saturation is of the form $X^{1/4}_Y\simeq F\times Y$ with $F$ being the cuspidal curve. Thus the family is $1/4$-equisaturated, although $F\not\simeq X_0^{1/4}$ and $X_0^{1/4} \not\simeq X_\lambda^{1/4}$, for $\lambda\neq 0$.
\end{ex}

We now give the following proposition which clarifies the motivation for defining $\alpha$-equisaturation through a product condition: point 1 and 4 show that such a condition implies constancy of the $\delta$ invariant, and of the $\mu$ invariant if the family is centered. Moreover, this proposition will be a useful tool throughout this section. For instance, point 2 shows that, even if the factor $F$ might not be isomorphic to $X_0^\alpha$, it is actually not so far from it.

\begin{prop}\label{Prop_ProductConsequences}
Let $p:X\to Y$ be an $\alpha$-equisaturated family. Let us write $q:F\times Y \simeq X^\alpha_Y\to X$ the relative saturation morphism and let us consider $\sigma:Y\to X$ such that $\sigma(y)=q(0,y)$. Then the following hold.
\begin{enumerate}
\item We have
$$
\Xbar\simeq\overline F\times Y\text{ and } (\Xbar)_y\simeq\overline{X_y} \quad\text{for every }y.
$$
In particular, the family admits a simultaneous normalization, and the total delta invariant of $X_y$ is constant.

\item For all $y$ in a dense Zariski-open subset of $Y$, one has $X_y^\alpha\simeq F$. Restriction to the special fiber gives the following inclusions
$$
\Ocal_{X_0}\subseteq\Ocal_{X_0^s}\subseteq\Ocal_F\subseteq\Ocal_{X_0^\alpha}\subseteq\Ocal_{\overline{X_0}}.
$$

\item The map $q : F\times Y\to X$ is a finite bijective morphism. Every singular point of $X_y$ distinct from $\sigma(y)$ is unibranch, and its unique preimage under $q$ is a smooth point of the fiber $F\times\{y\}$.

\item If the family in centered, then $\mu(X_y,\sigma(y))$ is constant for every $y\in Y$.
\end{enumerate}
\end{prop}

\begin{proof}
Let us prove 1. Since $X_Y^\alpha$ has the same normalization as $X$. Normalizing the product, and using the smoothness of $Y$, gives an isomorphism $\Xbar\simeq\overline F\times Y$ which preserves the base. Under this identification, the normalization map
$\pi:\Xbar\to X$ is the composition
$$
\overline F\times Y
\xrightarrow{\ \pi_F\times\mathrm{id}_Y\ }
F\times Y\xrightarrow{\ q\ }X,
$$
where $\pi_F:\overline F\to F$ is the normalization of $F$. For every $y$, the restriction $q_y:F\times\{y\}\to X_y$ is finite and bijective. Since the curves are reduced and defined over $\C$, it is birational. Consequently, $\pi_y:(\Xbar)_y\to X_y$ is finite and birational with normal source, hence is the normalization of $X_y$. Moreover, $\Xbar\to Y$ is smooth, so $\pi$ is a simultaneous normalization. Finally, by the $\delta$-constant criterion for simultaneous normalization \cite[Theorem~5.6]{ChiangHsiehLipman2006}, the total delta invariant of $X_y$ is constant.

Let us prove 2. By Theorem~\ref{Th_GenericHolderSpecialization}, there is a dense Zariski-open subset $U\subseteq Y$ such that $F\simeq(X_Y^\alpha)_y\simeq X_y^\alpha$, for every $y\in U$. By restricting the natural inclusions $\Ocal_X \subseteq \Ocal_{X^\alpha_Y} \subseteq \Ocal_{\overline{X}}$ to the special fiber, we get
$$
\Ocal_{X_0}\subseteq\Ocal_F
\subseteq\Ocal_{\overline{X_0}}.
$$
Let $h\in\Ocal_F$. Viewed as a function on $F\times Y$ independent of $y$, it is a section of algebra $\Ocal_{X_Y^\alpha}$. Its pullback to $\Xbar$ satisfies the relative $\alpha$-Hölder inequality of Proposition~\ref{Prop_RelativeHolderEstimate}. Restricting this inequality to $(\Xbar)_0=\overline{X_0}$ shows that $h$ belongs to $\Ocal_{X_0^\alpha}$. Finally, since $X_y^\alpha$ is Lipschitz saturated, the generic identification $F\simeq X_y^\alpha$ implies that $F$ is Lipschitz saturated (i.e $F^s=F$) and since we have $\Ocal_{X_0}\subseteq\Ocal_F\subseteq \Ocal_{\overline{X_0}}$, the universal property of Lipschitz saturation implies $\Ocal_{X_0^s}\subset\Ocal_F$.

We now prove 3. After shrinking the curve germ $(F,0)$, we can assume it is smooth away from $0$. Let $x\in X_y$ be a point distinct from $\sigma(y)$. Its unique preimage under $q_y$ is different from $0$ therefore it is a smooth point of $F\times\{y\}$. The normalization of $F$ has exactly one point above that point, so $\pi_y : (\overline{X})_y \to X_y$ has exactly one point above $x$. Thus $(X_y,x)$ is unibranch.

Finally, let us prove 4. Suppose that $\sing(p)=\sigma(Y)$. Since $q$ is bijective and $\sigma(y)=q(0,y)$, we have the set-theoretically equality $\pi^{-1}(\sigma(Y))=\pi_F^{-1}(0)\times Y$. As $\pi_y : (\overline{X})_y \to X_y$ is the normalization of $X_y$, it follows that $r(X_y,\sigma(y))=\#\pi_F^{-1}(0)=r(F,0)$ for every $y$. By 1, we already know that $\delta(X_y,\sigma(y))$ is constant, so the identity
$$
\mu(X_y,\sigma(y))
=2\delta(X_y,\sigma(y))-r(X_y,\sigma(y))+1
$$
gives constancy of $\mu$.
\end{proof}

\begin{rmq}\label{Rmq_LipschitzEquisaturation}
From Proposition \ref{Prop_ProductConsequences}(2) we can see that a Lipschitz equisaturated family satisfies $X_y^s\simeq X_0^s$ for every $y\in Y$.
\end{rmq}


\begin{rmq}\label{rmk_LipEquiIFFisoLipSatofFibers}
Let $p:X\to Y$ be a centered flat family of reduced curve germs. Then we actually have that $p$ is Lipschitz equisaturated if and only if
$$
X_y^s\simeq X_0^s
\qquad\text{for every }y\in Y.
$$
Assume that $X_y^s\simeq X_0^s$ for every $y$. The generic-projection argument in \cite[proof of Theorem~IV.8]{BrianconGallicoGranger1980DeformationsEqu}, which uses semicontinuity and applies over a smooth base of arbitrary dimension, gives a relative plane projection $X\to X'$ whose fiber maps are generic and whose image family is equisingular. The residual diagonal ideals of this plane family are equimultiple. Since the fiber projections are generic, specialization of integral dependence identifies $X_Y^s$ with $(X')_Y^s$; see \cite[Chap. I, \S5.1 and proof of Proposition~6.4.2]{Teissier1982VarietesPolairesII}. The latter proof also gives
$$
X_Y^s\simeq (X')_Y^s\simeq (X'_0)^s\times Y
\simeq X_0^s\times Y,
$$
with isomorphisms preserving the base.
\end{rmq}

By using the "$1/3$-Theorem" (\ref{Th_HolderSatSmoothness}), we can deduce from the preceding proposition that, for $\alpha>1/3$, the product condition implies the family to be automatically centered.

\begin{cor}\label{Cor_NoSplittingOneThird}
Let $\alpha\in(1/3,1]$ and let $p:X\to Y$ be an $\alpha$-equisaturated family. Then the family is centered and so $\mu(X_y,\sigma(y))$ is constant.
\end{cor}

\begin{proof}
Suppose that $x\in X_y$ is a singular point distinct from $\sigma(y)$. By Proposition \ref{Prop_ProductConsequences}(3), its unique preimage $u$ under $q_y$ is a smooth point of $F$. Moreover, by Proposition \ref{Prop_ProductConsequences}(1), the germ $(F,u)$ is the normalization of $(X_y,x)$. So, applying the inclusions of Proposition \ref{Prop_ProductConsequences}(2), to the germ of the family at $x$ gives
$$
\Ocal_{F,u}\subseteq\Ocal^\alpha_{X_y,x}
\subseteq\overline{\Ocal_{X_y,x}}=\Ocal_{F,u}.
$$
Thus $(X_y,x)^\alpha$ is smooth, so Theorem \ref{Th_HolderSatSmoothness} implies that $(X_y,x)$ is smooth, a contradiction. Hence $\sing(p)\subseteq\sigma(Y)$.

By Proposition \ref{Prop_ProductConsequences}(1), the total delta invariant of every fiber equals $\delta(X_0)>0$. Each fiber therefore has a singular point, necessarily $\sigma(y)$, and consequently $\sing(p)=\sigma(Y)$. The constancy of $\mu(X_y,\sigma(y))$ now follows from Proposition \ref{Prop_ProductConsequences}(4).
\end{proof}

Recall that $p$ is $\alpha$-Hölder equisingular if there is a homeomorphism of germs $\Phi : (X,0)\to(X_0,0)\times(Y,0)$ over $Y$ such that $\Phi$ and $\Phi^{-1}$ are $\alpha$-Hölder. The following proposition shows that $\alpha$-equisaturation implies $\alpha$-Hölder equisingularity, with a trivialization of a particular form.

\begin{prop}\label{Prop_ProductHolderTriviality}
Let $p:X\to Y$ be a $\alpha$-equisaturated family. Then, there exists a homeomorphism of germs $$\Phi:X\longrightarrow X_0\times Y$$ whose restriction to $X_0$ is the identity and such that both $\Phi$ and $\Phi^{-1}$ are $\alpha$-Hölder. More precisely, the morphisms $F\times Y\to X$ and $F\to X_0$ are Lipschitz and their inverses are $\alpha$-Hölder.
\end{prop}

\begin{proof}
Let $q:F\times Y\simeq X_Y^\alpha\to X$ be the natural
morphism, and denote its restriction to the special fiber
by $q_0:F\to X_0$. These maps are finite and bijective by
Proposition~\ref{Prop_ProductConsequences}, hence are
homeomorphisms on sufficiently small representatives. Define
$$
\Phi:=(q_0\times\mathrm{id}_Y)\circ q^{-1}
:X\longrightarrow X_0\times Y.
$$
Thus $\Phi(q(u,y))=(q_0(u),y)$. This is a homeomorphism commuting with the projections to $Y$, and its restriction to $X_0$ is the identity. It remains to establish the Hölder property of $\Phi$ for the metric induced by embeddings $F\subset\C^N$, $Y\subset\C^d$ and $X\subset\C^N\times Y$, with $p$ induced by the projection to $Y$. First, we can see the coordinate functions of $F$ as functions on $F\times Y$, and so as elements of the algebra of $X^\alpha_Y$. Consequently, by applying Proposition~\ref{Prop_RelativeHolderEstimate} to their pullbacks to $\Xbar$, we can consider a constant $C>0$ such that
$$
\|u-v\|\leq C\|q(u,y)-q(v,y)\|^\alpha
$$
for all $u,v\in F$ and $y\in Y$ sufficiently close to the origin. Now consider $(u,y),(v,z)\in F\times Y$. Since $p$ is the restriction of the coordinate projection $\C^N\times Y\to Y$ and $p(q(u,y))=y$, we have $\|y-z\| =\|p(q(u,y))-p(q(v,z))\| \leq \|q(u,y)-q(v,z)\|$. Moreover, since $q$ is analytic, it is locally Lipschitz and we chose a Lipschitz constant for $q$ that we call $L$. Hence we get
$$
\begin{aligned}
\|u-v\|
&\leq C\bigl(\|q(u,y)-q(v,z)\|
            +\|q(v,z)-q(v,y)\|\bigr)^\alpha\\
&\leq C\bigl(\|q(u,y)-q(v,z)\|+L\|y-z\|\bigr)^\alpha\\
&\leq C(1+L)^\alpha \|q(u,y)-q(v,z)\|^\alpha.
\end{aligned}
$$
Shrinking so that $\|q(u,y)-q(v,z)\|\leq1$, and using $0<\alpha\leq1$,
we conclude that
$$
\|(u,y)-(v,z)\|
\leq\|u-v\|+\|y-z\|
\leq\bigl(C(1+L)^\alpha+1\bigr)\|q(u,y)-q(v,z)\|^\alpha.
$$
Thus $q^{-1}$ is $\alpha$-Hölder. Restricting to $y=z=0$
also shows that $q_0^{-1}$ is $\alpha$-Hölder. Finally, $q_0\times\mathrm{id}_Y$ is Lipschitz, while $q_0^{-1}\times\mathrm{id}_Y$ is $\alpha$-Hölder because $0<\alpha\leq1$. Hence, we get that $\Phi=(q_0\times\mathrm{id}_Y)\circ q^{-1}$ and 
$\Phi^{-1}=q\circ(q_0^{-1}\times\mathrm{id}_Y)$ are both $\alpha$-Hölder.
\end{proof}

\begin{ex}\label{Ex_ProductFactorHolderTrivialization}
Consider the family of Example~\ref{Ex_ProductFactorDifferent}.
Writing $F=\{v^2=u^3\}$, its relative saturation morphism is
$$
q:F\times Y\longrightarrow X,
\qquad
q((u,v),\lambda)=(v+\lambda u,u^2,uv,\lambda).
$$
Its restriction to the special fiber is
$q_0:F\to X_0$, $q_0(u,v)=(v,u^2,uv)$.
The trivialization constructed in
Proposition~\ref{Prop_ProductHolderTriviality} is
$\Phi=(q_0\times\operatorname{id}_Y)\circ q^{-1}$, namely
$$
\begin{aligned}
\Phi:X&\longrightarrow X_0\times Y,\\
(t^3+\lambda t^2,t^4,t^5,\lambda)
&\longmapsto\bigl((t^3,t^4,t^5),\lambda\bigr).
\end{aligned}
$$
with $t$ being the parameter used in Example \ref{Ex_ProductFactorDifferent} to describe $X$. It preserves the base, restricts to the identity on $X_0$, and both $\Phi$ and $\Phi^{-1}$ are $1/4$-Hölder for the outer metrics. This example also shows that bi-$1/4$-Hölder triviality does not imply that the $1/4$-Hölder saturations of the fibers are isomorphic.
\end{ex}

\subsection{Lipschitz equisaturation and $\alpha$-equisaturation}\label{Subsec_LipschitzEquisaturation}

In this section, we give some precisions on the link between the classical notion of Lipschitz equisaturation and our notion of $\alpha$-equisaturation. First of all, a Lipschitz equisaturated family is, of course, an $\alpha$-equisaturated, but one can actually say more.
\begin{prop}
    Let $p:X\to Y$ be a Lipschitz-equisaturated family. Then the family is $\alpha$-equisaturated for every $\alpha\in (0,1]$.
\end{prop}

\begin{proof}
    Let $p:X\to Y$ be a Lipschitz-equisaturated family and let $\alpha\in (0,1]$. Then we have $X^s_Y\simeq X_0^s\times Y$. By Proposition \ref{Prop_InteractionsBetweenSaturations}, we know that $X^\alpha_Y = (X^s_Y)^\alpha_Y$. This gives us 
    $$X^\alpha_Y \simeq (X_0^s\times Y)^\alpha_Y \simeq (X_0^s)^\alpha\times Y \simeq X_0^\alpha\times Y$$
    with isomorphisms preserving the base.
\end{proof}
This is remarkable because, in general, if a family is $\beta$-equisaturated, it does not imply that the family is $\alpha$-equisaturated for $\alpha<\beta$.\\

We now investigate what are the conditions for an $\alpha$-equisaturated family to be Lipschitz-equisaturated. We start with the following proposition which reduce the problem to knowing when the $\alpha$-Hölder saturation of the central fiber is equal to its Lipschitz saturation.

\begin{prop}\label{Prop_HolderEquisaturationImpliesLipschitz}
Let $p:X\to Y$ be a family of reduced curve germs over a smooth
base, and let $\alpha\in(0,1]$. Then $p$ is Lipschitz equisaturated if and only if $p$ is $\alpha$-equisaturated and $X_0^\alpha=X_0^s$.
\end{prop}

\begin{proof}
Let $q:F\times Y\simeq X_Y^\alpha\to X$ be the relative saturation
morphism and write $q(u,y)=(q_y(u),y)$ with $q_y : F \to X_y$ for every $y\in Y$. By Proposition~\ref{Prop_ProductConsequences}(2), we have
$$
\Ocal_{X_0}\subseteq\Ocal_{X_0^s}\subseteq\Ocal_F\subseteq\Ocal_{X_0^\alpha}.
$$
If $p$ is Lipschitz equisaturated, it is clear that $X_0^\alpha \simeq X_0^s$ and so the first implication is clear. Conversely, if $X_0^\alpha \simeq X_0^\alpha$ and $\alpha \in (0,1]$, then the previous inclusions imply $\Ocal_F = \Ocal_{X_0^s}$. Since the inclusion $\Ocal_{X_0}\subseteq \Ocal_{F}$ is induced by $q_0^*$, we get that the morphism $q_0:F\to X_0$ is a bi-Lipschitz homeomorphism. In particular, we can consider a constant $L>0$ such that, for all $u,u'\in F$ sufficiently close to $0$,
$$
\|q_0(u)-q_0(u')\|\geq \frac{1}{L}\|u-u'\|.
$$
Choose local coordinates $y_1,\ldots,y_d$ on $Y$. Since $Q(u,y)-Q(u,0)$ is a holomorphic function vanishing when $y=0$, we can write
$$
Q(u,y)-Q(u,0)=\sum_{j=1}^d y_jR_j(u,y),
$$
Moreover, since the $R_j(.,y)$ are holomorphic, they are locally Lipschitz and so one can consider $C>0$ such that 
$$
\|(q_y-q_0)(u) - (q_y-q_0)(u')\| \leq C \|y\| \|u-u'\|
$$
Hence we have 
$$
\begin{aligned}
\|q_y(u)-q_y(u')\|
&\geq \|q_0(u)-q_0(u')\|
-\|(q_y-q_0)(u)-(q_y-q_0)(u')\|\\
&\geq \left(\frac1L-C\|y\|\right)\|u-u'\|.
\end{aligned}
$$
After shrinking $Y$ so that $C\|y\|\leq 1/(2L)$ for example, we obtain
$$
\|q_y(u)-q_y(u')\|\geq \frac{1}{2L}\|u-u'\|.
$$
This exactly means that $q_y :F \to X_y$ is a bi-Lipschitz homeomorphism. Since it is also a finite and birational morphism and since $F$ is Lipschitz saturated, we get $F\simeq X_y^s$ by the universal property of the Lipschitz saturation. Finally, we get $X_0^s\simeq F\simeq X_y^s$, for all $y\in Y$, so by Remark \ref{rmk_LipEquiIFFisoLipSatofFibers}, the family $p$ is Lipschitz equisaturated.
\end{proof}

Since we had already establish in Section \ref{Subsec_BoundsReducedCurves} some lower bounds on $\alpha$ for the $\alpha$-Hölder saturation of a reduced curve to be equal to its Lipschitz saturation. We get the following corollary.

\begin{cor}
    Let $p:X\to Y$ be a family of reduced curves over a smooth base. Write $m_0 = m(X_0)$ and $\mu_0 = \mu(X_0)$. Then
    \begin{itemize}
        \item If $p$ is Lipschitz-equisaturated, it is $\alpha$-equisaturated for all $\alpha\in (0,1]$.
        \item If $p$ is $\alpha$-equisaturated for $\alpha>\frac{\mu_0-m_0+2}{\mu_0-m_0+3}$, then it is Lipschitz equisaturated.
    \end{itemize}
\end{cor}

\subsection{Small exponents and constancy of the Milnor number}

We now study the converse of Proposition~\ref{Prop_ProductConsequences}\textup{(4)}. For a centered family, we will show that constancy of the Milnor number implies $\alpha$-equisaturation for every sufficiently small $\alpha$.

In this subsection, we will use several times the notion of "simultaneous parametrization" of a family, so we recall this notion here. Let $(X,0)$ be a reduced complex curve germ with branches $X_1,\ldots,X_r$. Recall that a \textit{parametrization} of $(X,0)$ is a holomorphic map
$$
\psi:\coprod_{i=1}^r(\C,0_i)\longrightarrow(X,0)
$$
whose restriction to the $i$-th component is a normalization parametrization of $X_i$. In ambient coordinates, this restriction has the form $\psi_i(t_i)=\bigl(\phi_{i,1}(t_i),\ldots,\phi_{i,n}(t_i)\bigr)$.

Let now $p:X\to Y$ be a centered $\mu$-constant family of reduced curve germs, with smooth base $Y$ and section $\sigma$. After shrinking, the family admits a simultaneous normalization $\pi:\Xbar\to X$ and an isomorphism $\Xbar\simeq\overline{X_0}\times Y$ preserving the base and identifying $\pi^{-1}(\sigma(Y))$ with $\coprod_{i=1}^r(\{0_i\}\times Y)$ as sets. Choosing a coordinate $t_i$ on each normalized branch of $X_0$, we obtain a \textit{simultaneous parametrization} along $\sigma$, that is, a holomorphic map
$$
\Psi:\left(\coprod_{i=1}^r(\C,0_i)\right)\times Y \longrightarrow X
$$
such that $p\circ\Psi=\operatorname{pr}_Y$,
$\Psi(0_i,y)=\sigma(y)$ for every $i$, and, for every $y$ near $0$,
the induced map
$$
\Psi_y:\coprod_{i=1}^r(\C,0_i)\longrightarrow(X_y,\sigma(y))
$$
is a normalization parametrization of $X_y$.
In relative ambient coordinates in which $\sigma(y)=(0,y)$, we may write $\Psi_i(t_i,y)=\bigl(\phi_{i,1}(t_i,y),\ldots,\phi_{i,n}(t_i,y),y\bigr)$, where the functions $\phi_{i,j}$ are holomorphic jointly in $t_i$ and $y$ and satisfy $\phi_{i,j}(0,y)=0$.

\subsubsection{The relative seminormalization exponent}

Let $p:X\to Y$ be a family. Let $X_Y^+$ denote the relative seminormalization obtained by replacing $\Ical_\alpha$ with $\sqrt{\Ical}$ in the definition of $X_Y^\alpha$. Let $\pi : \overline{X}\to X$ be the normalization morphism. A holomorphic function $h$ on $\Xbar$ belongs to $\Ocal_{X_Y^+}$ if and only if, for every $(x,x') \in \overline{X}\times_Y \overline{X}$, we have $h(x)=h(x')$ whenever $\pi(x)=\pi(x')$. If $(x,x')$ belongs to $\overline{X}\times\overline{X}$ instead, then the condition $\pi(x)=\pi(x')$ implies $\overline{X}\times_Y\overline{X}$. Consequently, we have $X_Y^+=X^+$.

\begin{definition}\label{Def_RelativeSeminormalThreshold}
The \textit{relative seminormalization exponent} of $p:X\to Y$ is
\[
\alpha_{\mathrm{sn}}(X/Y) := \max\Big\{\alpha\in(0,1]:X_Y^\alpha=X^+\Big\},
\]
where the equality is understood as an equality of the corresponding
algebras inside $\pi_*\Ocal_{\Xbar}$.
\end{definition}

By applying Proposition \ref{Prop_FinitenessOfRationalIdeals} to the relative diagonal ideal $\Ical \subset \Ocal_{\overline{X}\times_Y \overline{X}}$, we know that there exists $e\in \N$ such that $\Ical_\alpha = \sqrt{\Ical}$, for all $0<\alpha\leq 1/e$. Since, for $\alpha\in (0,1]$, the equality $\Ical_\alpha=\sqrt{\Ical}$ implies $X_Y^\alpha = X_Y^+$ by definition, we get that $\alpha_{sn}(X/Y)$ is well defined. Moreover, we have
$$
0<\alpha\leq\alpha_{\mathrm{sn}}(X/Y) \iff X_Y^\alpha=X_Y^+ = X^+.
$$

\begin{rmq}
    When $Y$ is a point and $X$ is irreducible, Corollary \ref{Cor_HolderEqNormaIffBeta_g} gives
    $$
    \alpha_{\mathrm{sn}}(X/\{0\})=\frac{1}{\beta_g(X)}.
    $$
\end{rmq}

In a centered $\mu$-constant family, one can also define $\alpha_{\mathrm{sn}}$ as the maximum ideal for which the rational power of $\Ical$ coincide with its radical.
\begin{lemma}\label{Lem_EquivalenceDefAlphaSN}
    If $p:X\to Y$ is a centered and $\mu$-constant family, then 
    $$
    \Ical_\alpha = \sqrt{\Ical} \iff X^\alpha_Y \simeq X^+_Y
    $$
\end{lemma}

\begin{proof}
The first implication is immediate by definition. To prove the other implication, let us assume that $X^\alpha_Y \simeq X^+_Y$. Since the family is centered and $\mu$-constant, it admits a simultaneous normalization $\pi:\Xbar\to X$, and the number $r$ of branches of the fibers is constant. After shrinking, we may write each irreducible component of $\overline{X}$ on the form $(\C,0)\times Y$. We choose the coordinate $t_i$ on the first factor so that the inverse image of $\sigma(Y)$ in $\Xbar_i$ is given by $t_i=0$. Thus, for each $y$, the function $t_i$ is a local parameter on the normalization of the $i$-th branch of $X_y$, vanishing at the point above $\sigma(y)$.
Let $h_i$ be the holomorphic function on $\Xbar$ equal to $t_i$ on $\Xbar_i$ and to zero on the other components. The normalization is an isomorphism outside $\sigma(Y)$, while all the $h_i$ vanish
at every point above $\sigma(Y)$. Hence each $h_i$ is constant on the fibers of $\pi$, and therefore belongs to $\Ocal_{X_Y^+}$. On the component of $\overline{X}\times_Y\overline{X}$ corresponding to branches $i$ and $j$, write $s,t$ for the two normalization parameters. To identify $\sqrt{\Ical}$, recall that $V(\Ical) = \overline{X}\times_X\overline{X}$ consists of pairs of points of $\Xbar$ having the same image under $\pi$. Since $\pi$ identifies only the origins of distinct branches in each fiber, we have
$$
\sqrt{\Ical}=
\begin{cases}
(s-t),&i=j,\\
(s,t),&i\ne j.
\end{cases}
$$
In the first case, $\Delta h_i=s-t$; in the second, $\Delta h_i=s$ and $\Delta h_j=-t$. Thus the differences $\Delta h_1,\ldots,\Delta h_r$ generate $\sqrt{\Ical}$ on $\overline{X}\times_Y\overline{X}$. Since $X_Y^\alpha\simeq X_Y^+$, each $h_i$ belongs to $\Ocal_{X_Y^\alpha}$, hence $\Delta h_i\in\Ical_\alpha$. Consequently, $\sqrt{\Ical}\subseteq\Ical_\alpha$, and equality follows since $\Ical_\alpha\subseteq\sqrt{\Ical}$ is always true.
\end{proof}

To compute the relative seminormalization exponent, one can use the following expression in terms of orders of the relative diagonal ideal along arcs.

\begin{prop}\label{Prop_RelativeSeminormalThresholdCriteria}
Assume that $p:X\to Y$ is a centered and $\mu$-constant family. Then,
\begin{equation}\label{Eq_RelativeSeminormalThresholdArcs}
\alpha_{\mathrm{sn}}(X/Y) = \inf_\gamma \frac{\ord_u(\gamma^*\sqrt{\Ical})}{\ord_u(\gamma^*\Ical)},
\end{equation}
where $\gamma:(\C,0)\to(W,w)$ ranges over analytic arcs with $w\in\pi^{-1}(x)\times\pi^{-1}(x)$ and $0<\ord_u(\gamma^*\Ical)<+\infty$.
\end{prop}

\begin{proof}
    By Lemma \ref{Lem_EquivalenceDefAlphaSN}, we can find the maximum exponent for which $\Ical_\alpha = \sqrt{I}$. The inclusion $\Ical_\alpha \subseteq \sqrt{I}$ is always true. To get the reverse inclusion, by the valuative criterion, it is enough to have $\ord_u(\gamma^*\sqrt{I}) \geq \alpha \ord_u(\gamma^*\Ical)$ for every arc $\gamma$ and since arcs for which $\ord_u(\gamma^*\Ical)=0$ impose no condition, this is equivalent to $\alpha \leq \inf_\gamma \ord_u(\gamma^*\sqrt{I}) / \ord_u(\gamma^*\Ical)$. 
\end{proof}

\begin{rmq}
Assume that $p:X\to Y$ is centered and $\mu$-constant. Recall from Subsection \ref{Subsec_FinitenessExponents} that, in the normalized blow-up $\overline{\mathrm{Bl_\Ical(W)}}\to W$ of $\Ical$, the finitely many prime components $E_1,\ldots,E_s$ of the divisor defined by the pull-back of $\Ical$ determine its Rees valuations $v_i=\ord_{E_i}$. By the valuative
criterion and Lemma \ref{Lem_EquivalenceDefAlphaSN},
$$
\alpha_{\mathrm{sn}}(X/Y) = \min_{1\leq i\leq s} \frac{v_i(\sqrt{\Ical})}{v_i(\Ical)}.
$$
Thus the infimum over analytic arcs can be computed by finitely many Rees valuations.
\end{rmq}

\subsubsection{Constancy of the milnor number and $\alpha$-equisaturation}

In this subsection we establish the link for a family between being $\mu$-constant and being $\alpha$-equisaturated. We start by showing that being a centered $\mu$-constant family implies $\alpha$-equisaturation for a sufficiently small $\alpha$.

\begin{prop}\label{Prop_RelativeSeminormalizationEquisaturated}
Let $p:X\to Y$ be a centered $\mu$-constant family. Then, for every $0<\alpha\leq\alpha_{\mathrm{sn}}(X/Y)$, the family is $\alpha$-equisaturated with special fiber $X_0^+$. Moreover, one has $X_y^\alpha=X_y^+$ for every $y\in Y$.
\end{prop}

\begin{proof}
The equality $X_Y^\alpha=X^+$ follows from the definition of $\alpha_{\mathrm{sn}}(X/Y)$ and the discussion following it. Choose relative parameters $t_i$ on the components $\Xbar_i$ of the simultaneous normalization, and identify $\Xbar\simeq\overline{X_0}\times Y$. A holomorphic function on $\Xbar$ is then a tuple $(h_1(t_1,y),\ldots,h_r(t_r,y))$. Since $X_0^+$ is obtained by identifying the origins of the normalized branches, a holomorphic function $h$ on this normalization belongs to $\Ocal_{X_0^+\times Y}$ if and only if $h_1(0,y)=\cdots=h_r(0,y)$ for every $y\in Y$. This is exactly the condition for $h$ to belong to $\Ocal_{X^+_Y}$. Thus the two algebras coincide inside $\Ocal_{\Xbar}$, giving an isomorphism $X^+\simeq X_0^+\times Y$ over $Y$. For each fixed $y$, the same condition characterizes $\Ocal_{X_y^+}$, so this isomorphism also identifies $(X^+)_y \simeq X_0^+$ with $X_y^+$. Finally, by Proposition \ref{Prop_ProductConsequences}(2), we have $X_0^+ = X_y^\alpha$ for $y\neq 0$ and $\Ocal_{X_0} \subseteq \Ocal_{X_0^+} \subseteq \Ocal_{X_0^\alpha}$. Since the inclusion $\Ocal_{X_0^\alpha} \subseteq \Ocal_{X_0^+}$ is always true, we conclude that $X^\alpha_y \simeq X_0^+$ for every $y\in Y$.
\end{proof}

\begin{ex}\label{Ex_RelativeSeminormalExponent}
It is interesting to notice that the relative seminormalization exponent $\alpha_{\mathrm{sn}}(X/Y)$ is not equal, in general, to the minimum of the seminormalization exponents of the fibers $\alpha_{\mathrm{sn}}(X_y)$. Let us consider, for example, the family
$$
\Ocal_{Y,0} = \C\{\lambda\} \inj \Ocal_{X,0} = \C\{\tau^3,\tau^5+\lambda\tau^4,\tau^7\} \inj \Ocal_{\overline{X},0} = \C\{\lambda,\tau\}.
$$
This is a centered flat family of branches with $\delta$ and $\mu$ constant. The characteristic sequence of $X_0$ is $(3,5)$ and the one of $X_y$ for $y\neq 0$ is $(3,4)$. By Corollary \ref{Cor_HolderEqNormaIffBeta_g}, we have $\alpha_{\mathrm{sn}}(X_0)=1/5$ and $\alpha_{\mathrm{sn}}(X_y)=1/4$ for $y\neq 0$. So we get $\inf_\lambda\alpha_{\mathrm{sn}}(X_\lambda)=1/5$.

To compute the relative exponent, consider the diagonal ideal on
$W=\Xbar\times_Y\Xbar$:
$$
\Ical= \big(s^3-t^3, s^5-t^5+\lambda(s^4-t^4), s^7-t^7\big).
$$
Its radical is $(s-t)$. Let $\zeta$ be a primitive cube root of unity and consider the arc $\gamma : (\C,0) \to (W,(0,0))$ such that $\gamma(u) = (\lambda(u),s(u),t(u)) = (\zeta^2u, u,\zeta u)$.
The first two generators of $\Ical$ vanish identically along
$\gamma$, since $s^3-t^3=0$ and 
$$s^5-t^5+\lambda(s^4-t^4)
=\bigl(1-\zeta^2+\zeta^2(1-\zeta)\bigr)u^5=0.
$$
The third generator is $(1-\zeta)u^7$, whereas
$s-t=(1-\zeta)u$. Hence $\ord_u(\gamma^*\Ical)=7$ and $\ord_u\bigl(\gamma^*(s-t)\bigr)=1$. The arc criterion therefore gives $\alpha_{\mathrm{sn}}(X/Y)\leq1/7$.
Finally, we have $$ \alpha_{\mathrm{sn}}(X/Y) \leq \frac17 < \frac15 =\inf_\lambda\alpha_{\mathrm{sn}}(X_\lambda).$$
In particular, this implies that, for every $1/7<\alpha\leq1/5$, all fibers satisfy $X_\lambda^\alpha=X_\lambda^+$, while $X_Y^\alpha\ne X_Y^+$.\\
\end{ex}

Recall from Section \ref{Sec_ComputingHolderSat} that, for a germ of irreducible curve, the $\alpha$-Hölder saturation and the seminormalization coincide for $\alpha\leq 1/\beta_g$. In the following theorem, we get a bound of the same spirit linking equiseminormalization and $\alpha$-equisaturation.
 
\begin{theorem}\label{Th_UniformSmallExponent}
Let $p:X\to Y$ be a centered $\mu$-constant family, and put
$\delta=\delta(X_0)$. Then
$$
\alpha_{\mathrm{sn}}(X/Y)\geq\frac{1}{2\delta+1}.
$$
Consequently, for every $0<\alpha\leq1/(2\delta+1)$, the family is $\alpha$-equisaturated and so it admits a bi-$\alpha$-Hölder trivialization for the outer metric.
\end{theorem}

\begin{rmq}
Example \ref{Ex_RelativeSeminormalExponent} shows that this bound can be attained. Indeed, in this example, we have proven that $\alpha_{\mathrm{sn}} \leq 1/7$. The invariant $\delta$ is constant for the considered family and equal to $3$. Hence $1/(2\delta +1) = 1/7$ and the bound of the Theorem is exactly equal to $\alpha_{\mathrm{sn}}$.
\end{rmq}

\begin{proof}
If $\delta=0$, the family is smooth and the assertion follows from the simultaneous parametrization. Suppose that $\delta>0$, put $N=2\delta$, and let $0<\alpha\leq 1/(N+1)$. By Lemma~\ref{Lem_EquivalenceDefAlphaSN}, it is enough to prove $\sqrt{\Ical}\subseteq\Ical_\alpha$, since the reverse inclusion always holds. Let us take a simultaneous parametrization and consider the functions $h_i$ defined in the proof of that lemma.

To get a bound on the order of $\Ical$ along an arc, we first show that $h_i^N,h_i^{N+1}\in\Ocal_X$ for every $i$. For every $y$, write the conductor of $X_y$ in its normalization as $\mathfrak c_y=\bigoplus_{i=1}^r t_i^{c_i(y)}\C\{t_i\}$. The classical conductor bound gives
$$
\sum_{i=1}^r c_i(y) =\dim_\C\big(\Ocal_{\overline{X_y}}/\mathfrak c_y\big) \leq2\delta(X_y,\sigma(y))=N,
$$
thus the restrictions of $h_i^N$ and $h_i^{N+1}$ belong to $\Ocal_{X_y}$ for every $y$.

To pass from the fibers to $X$, consider $\mathcal Q=p_*(\pi_*\Ocal_{\Xbar}/\Ocal_X)$. Since $\pi$ is finite, the quotient $\pi_*\Ocal_{\Xbar}/\Ocal_X$ is a coherent $\Ocal_X$-module. Moreover, since the family is centered, this quotient is supported on $\sigma(Y)$. Hence its direct image $\Qcal$ is a finite $\Ocal_Y$-module. We will denote by $\Qcal_y$ the stalk of $\Qcal$ at $y\in Y$ and $\Qcal(y) = \Qcal_y/\m_{Y,y}\Qcal_y$ the fiber at $y$. By simultaneous normalization, its fiber $\Qcal(y)$ is isomorphic to $\Ocal_{\overline{X_y}}/\Ocal_{X_y}$ and thus has dimension $\delta$. Choose local sections $s_1,\ldots,s_\delta$ of $\Qcal$ whose images in the fiber $\Qcal(0)$ form a basis. By Nakayama's lemma, these sections generate $\Qcal_0$ as an $\Ocal_{Y,0}$-module. So, after shrinking $Y$, the sections $\{s_i\}_i$ generate $\mathcal Q$. Then, their images in every fiber generate $\Qcal(y)$ and since the dimension is constant and equal to $\delta$, they actually form a basis of $\Qcal(y)$ as a $\C$-vector space. Write the classes of $h_i^N$ and $h_i^{N+1}$ in $\Qcal$ in terms of these generators. Because we have proven that their restrictions to any $X_y$ belong to $\Ocal_{X_y}$, their classes in every fiber is zero, and so each coefficient vanishes at every $y$. Since $Y$ is reduced, the coefficients are zero. This proves the claim.

Let $\gamma=(\varphi,\varphi',\eta): (\C,0)\to\Xbar_i\times_Y\Xbar_j$ be an arc where $\eta$ is the common base parameter, and let $v_0=\min\{\ord_u(\varphi),\ord_u(\varphi')\}$. By the arc criterion, we have to show that $\ord_u(\gamma^*\sqrt{\Ical}) \geq \alpha\ord_u(\gamma^*\Ical)$.

Suppose first that $i=j$. By the description of $\sqrt{\Ical}$ in Lemma \ref{Lem_EquivalenceDefAlphaSN}, we have $\ord_u(\gamma^*\sqrt{\Ical})=\ord_u(\varphi-\varphi')$. Moreover, since $h_i^N,h_i^{N+1}\in\Ocal_X$, we have
$$
\ord_u(\gamma^*\Ical) \leq \min\big\{ \ord_u(\varphi^N-(\varphi')^N), \ord_u(\varphi^{N+1}-(\varphi')^{N+1}) \big\}.
$$
If $\ord_u(\varphi-\varphi')>v_0$, then $\varphi$ and $\varphi'$ have the same leading term. As in the proof of Proposition \ref{Prop_ComputeSemigroupOfTheHolderSat}, we get $\ord_u(\varphi^N-(\varphi')^N) =(N-1)v_0+\ord_u(\varphi-\varphi')$. Therefore,
$$
\begin{aligned}
\alpha\ord_u(\gamma^*\Ical)
&\leq \alpha\big((N-1)v_0+\ord_u(\varphi-\varphi')\big)\\
&\leq \alpha N\ord_u(\varphi-\varphi')\\
&\leq \ord_u(\varphi-\varphi').
\end{aligned}
$$
If $\ord_u(\varphi-\varphi')=v_0$, write $\varphi(u)=a u^{v_0}+\cdots$ and $\varphi'(u)=a'u^{v_0}+\cdots$ where $a\ne a'$ and one of these coefficients may be zero. Since $N$ and $N+1$ are coprime, the equalities $a^N=(a')^N$ and $a^{N+1}=(a')^{N+1}$ cannot hold simultaneously. Thus, there exists $k\in\{N,N+1\}$, such that $\ord_u(\varphi^k-(\varphi')^k)=kv_0$, and so
$$
\alpha\ord_u(\gamma^*\Ical) \leq \alpha k v_0 \leq v_0 =\ord_u(\varphi-\varphi').
$$

Now suppose that $i\neq j$. Lemma \ref{Lem_EquivalenceDefAlphaSN} gives $\gamma^*\sqrt{\Ical}=(\varphi,\varphi')$, while $\gamma^*\Delta(h_i^N)=\varphi^N$ and
$\gamma^*\Delta(h_j^N)=-(\varphi')^N$. Therefore
$$
\alpha\ord_u(\gamma^*\Ical) \leq \alpha N v_0\leq v_0=\ord_u(\gamma^*\sqrt{\Ical}).
$$
Hence the required inequality holds for every arc, and the arc
criterion gives $\sqrt{\Ical}\subseteq\Ical_\alpha$. This proves $\alpha_{\mathrm{sn}}(X/Y)\geq 1/(2\delta+1)$. The assertion about Hölder triviality follows from Proposition \ref{Prop_ProductHolderTriviality}.
\end{proof}

The following statement summarize the relation between $\mu$-constancy of $\alpha$-equisaturation.
\begin{cor}\label{Cor_MuSomeAlpha}
Let $p:X\to Y$ be a centered flat family of reduced curve germs
with section $\sigma$. The following conditions are equivalent:
\begin{enumerate}[label=\textup{(\arabic*)},leftmargin=*]
\item The Milnor number $\mu(X_y,\sigma(y))$ is constant.

\item There exists $\varepsilon\in[\frac{1}{2\delta+1},1]$ such that the family is
$\alpha$-equisaturated for every $0<\alpha\leq\varepsilon$.

\item The family is $\alpha$-equisaturated for some
$\alpha\in(0,1]$.
\end{enumerate}
\end{cor}

\begin{proof}
The implication $(1)\Rightarrow(2)$ follows from Theorem~\ref{Th_UniformSmallExponent}. The implication $(2)\Rightarrow(3)$ is immediate. The implication $(3)\Rightarrow(1)$ is exactly Proposition \ref{Prop_ProductConsequences}(4).
\end{proof}

In the special case of a family of irreducible curves with constant semigroup, we get a sharper bound.


\begin{prop}\label{Prop_ConstantSemigroupBound}
Let $p:X\to Y$ be a centered family of reduced unibranch curve germs with constant value semigroup $\Gamma$. Put $\delta=\delta(X_0)$, let $n=\min(\Gamma\setminus\{0\})\geq2$, and denote by $c$ the
conductor exponent of $\Gamma$. For every nontrivial $n$-th
root of unity $\xi$, let $b_\xi:=\min\{k\in\Gamma:\xi^k\ne1\}$ and $b(\Gamma):= \max_{\substack{\xi^n=1,\xi\ne1}}b_\xi$.
Then, we have $b(\Gamma)\leq c+1\leq2\delta+1$ and
$$
\alpha_{\mathrm{sn}}(X/Y)\geq\frac{1}{b(\Gamma)}.
$$
In particular, the family is $\alpha$-equisaturated for every
$0<\alpha\leq1/b(\Gamma)$.
\end{prop}

\begin{rmq}
The number $b(\Gamma)$ is the largest characteristic exponent
of the monomial curve with local ring $\C\{t^k:k\in\Gamma\}$,
obtained by applying the construction of the characteristic sequence
to the exponent set $\Gamma$. Since $\Gamma\subseteq E(X_y)$,
the largest characteristic exponent $\beta_g(X_y)$ of each fiber
satisfies
$$
\beta_g(X_y)\leq b(\Gamma).
$$
By Corollary~\ref{Cor_HolderEqNormaIffBeta_g}, this gives
$X_y^\alpha=\overline{X_y}$ for $\alpha\leq1/b(\Gamma)$. The proposition shows that the same bound holds for the relative saturation, uniformly in the parameter. Of course, the inclusion $\Gamma\subseteq E(X_y)$ is usually strict, so the bound $1/b(\Gamma)$ is not necessarily attained. However, with the present assumptions it is easier to calculate since constancy of $\Gamma$ does not require the characteristic sequences of the fibers to remain constant.
\end{rmq}

\begin{proof}
We use the simultaneous parametrization $\Xbar\simeq(\C,0)\times Y$, with parameter $t$ such that $\pi^{-1}(\sigma(Y))=\{t=0\}$, where $\sigma(Y)$ is the center of the family. First, we may choose $t$ so that $t^n\in\Ocal_X$.
Indeed, lift a function of order $n$ on $X_0$ to a function $f\in\Ocal_X$ and replace $f$ by $f-(f\circ\sigma)\circ p$, so that it vanishes on the center. Since every fiber has multiplicity $n$, the coefficients of $t,\ldots,t^{n-1}$ in its pullback vanish identically. Thus $f(t,y)=t^n U(t,y)$, where $U$ is a unit. Replacing $t$ by $tU(t,y)^{1/n}$ gives the required parameter.

We claim that, for every $k\in\Gamma$ with $k\leq b(\Gamma)$, there exists $f_k\in\Ocal_X$ of the form
$$
f_k(t,y)=t^k+\sum_{j>k}a_{k,j}(y)t^j.
$$
We construct these functions in increasing order of $k$,
starting with $f_0=1$ and $f_n=t^n$. Suppose that the
$f_\ell$ have been constructed for all $\ell\in\Gamma$
with $\ell<k$. Lift a function of order $k$ on $X_0$
to $\Ocal_X$. By successively subtracting suitable
$\Ocal_Y$-multiples of the $f_\ell$, in increasing order
of $\ell$, we cancel its coefficients of degrees
$\ell<k$ belonging to $\Gamma$. If a coefficient of degree less than $k$ remained nonzero on some fiber, the first such degree would be the order
of a function on that fiber without belonging to $\Gamma$. This contradicts the constancy of the value semigroup. Hence all coefficients below degree $k$ vanish identically. Moreover, the coefficients used in the subtractions vanish
at $y=0$, so the coefficient of $t^k$ remains nonzero there. Dividing by this coefficient gives $f_k$, after shrinking $Y$.

Now write $(s,t,y)$ for the coordinates on $W=\Xbar\times_Y\Xbar$. By Lemma~\ref{Lem_EquivalenceDefAlphaSN} and Proposition~\ref{Prop_RelativeSeminormalThresholdCriteria},
it is enough to show that
$$
\ord_u(\psi^*\Ical)
\leq b(\Gamma)\ord_u(\varphi-\varphi')
$$
for every analytic arc
$\psi=(\varphi,\varphi',\lambda):(\C,0)\to(W,0)$ such that
$0<\ord_u(\psi^*\Ical)<+\infty$. Let $v_0,v_1$ and $e$ be defined as in the proof of Proposition~\ref{Prop_ComputeSemigroupOfTheHolderSat}. If $e=1$, then $\ord_u(\varphi-\varphi')=v_1$. By using the fact that $t^n\in\Ocal_X$ and by Lemma \cite[Lemma VI.1.1]{BrianconGallicoGranger1980DeformationsEqu}, we have
$$
\ord_u(\psi^*\Ical)\leq\ord_u\bigl(\varphi^n-(\varphi')^n\bigr)=(n-1)v_0+v_1\leq nv_1
\leq b(\Gamma)\ord_u(\varphi-\varphi').
$$

Suppose now that $e\neq 1$. Then we can choose $k\in\Gamma$ with $k\leq b(\Gamma)$ and $e\nmid k$. By the same lemma, we have 
$\ord_u\bigl(\varphi^k-(\varphi')^k\bigr)=kv_0$ and $\ord_u(\varphi-\varphi')=v_0$. Moreover,
$$
f_k(\varphi,\lambda)-f_k(\varphi',\lambda) =
\varphi^k-(\varphi')^k+
\sum_{j>k}a_{k,j}(\lambda)
\bigl(\varphi^j-(\varphi')^j\bigr)
$$
and since the coefficients $a_{k,j}(\lambda(u))$ are regular, then every term in the sum has order at least $jv_0>kv_0$. So we deduce that
$$
\ord_u(\psi^*\Ical)\leq\ord_u\bigl(f_k(\varphi,\lambda)-f_k(\varphi',\lambda)\bigr)=kv_0\leq b(\Gamma)\ord_u(\varphi-\varphi').
$$
and so $\alpha_{\mathrm{sn}}(X/Y)\geq1/b(\Gamma)$.

Finally, since $c$ and $c+1\in\Gamma$ are coprime, we have $b(\Gamma)\leq c+1
$ by definition of $b(\Gamma)$. The inequality $c+1\leq 2\delta +1$ comes from the conductor bound recalled in the proof of Theorem~\ref{Th_UniformSmallExponent}.
\end{proof}

\begin{ex}\label{Ex_ConstantSemigroupContinuation}
See that $1/b(\Gamma)$ is a more precise bound than $1/(2\delta+1)$. For example, consider the family
$$
\tau\longmapsto(\tau^4,\tau^6+\lambda\tau^7,\tau^8,\tau^9,\tau^{11}).
$$
The value semigroup $\Gamma$ is $\{0,4,6,8,9,10,\ldots\}$ for every $\lambda$ and $\delta=5$, so $1/(2\delta+1) = 1/11$. However $c=8$ and $b(\Gamma)=9$ so $1/b(\Gamma) = 1/9$.
\end{ex}

\subsection{Whitney equisingularity and $\alpha$-equisaturation}

In this subsection, we show that $\alpha$-equisaturation is  Whitney equisingularity if $\alpha$ is large enough, the explicit bound  depending on the multiplicity and the number of  branches of $X_0$. We start with the following technical lemma which we will use to study the multiplicity of the fibers of an $\alpha$-equisaturated family.
\begin{lemma}\label{Lem_MultiplicitySandwich}
Let $p:X\to Y$ be a centered $\alpha$-equisaturated family. Then
$$
m(X_0^\alpha)\leq m(F)\leq m(X_y,\sigma(y))\leq m(X_0).
$$
\end{lemma}

\begin{proof}
By Proposition~\ref{Prop_ProductConsequences}\textup{(2)},
the inclusion $\Ocal_F\subseteq\Ocal_{X_0^\alpha}$ inside
their common normalization gives $m(X_0^\alpha)\leq m(F)$. The same proposition gives $F\simeq X_y^\alpha$ for $y$ in a dense open subset of $Y$. Hence $m(F)=m(X_y^\alpha)\leq m(X_y)$ on this subset. By upper semicontinuity of $y\mapsto m(X_y)$, the inequality $m(F)\leq m(X_y)$ holds for every $y\in Y$. Upper semicontinuity also gives $m(X_y)\leq m(X_0)$ after shrinking $Y$, which proves the result.
\end{proof}

\begin{definition}
    Let $B$ be a irreducible complex analytic curve, let $\beta_0,\beta_1,\ldots,\beta_g$ be its characteristic exponents, and put $e_i=\gcd(\beta_0,\ldots,\beta_i)$. Then we define
    $$
    \alpha_m(B):= \max_{1\leq i\leq g}\frac{\beta_0-e_i}{\beta_i}.
    $$
\end{definition}

\begin{rmq}
    By Proposition \ref{Prop_HolderSatWithSameMultiplicity}, we have $m(B^\alpha)=m(B)$ if and only if $\alpha>\alpha_m(B)$.
\end{rmq}

We now obtain two lower bounds on $\alpha$ for an $\alpha$-equisaturated family to be Whitney equisingular. The first one depends on the characteristic sequence of $X_0$. The second one is less precise but it only depends on the multiplicity and the number of branches of $X_0$.

\begin{theorem}\label{Th_EquisaturationWhitney}
Let $p:X\to Y$ be a centered $\alpha$-equisaturated family. If
$$
\alpha>\max_{B\text{ branch of }X_0}\alpha_m(B),
$$
then $m(X_y)=m(X_0)$ for every $y$, and the family is Whitney equisingular. In particular, writing $r_0=r(X_0)$ and $m_0 = m(X_0)$, it suffices that
$$
\alpha>\frac{m_0-r_0}{m_0-r_0+2}.
$$
For unibranch families, this sufficient bound is $(m_0-1)/(m_0+1)$.
\end{theorem}

\begin{proof}
Write $X_0=B_1\cup\cdots\cup B_r$ and denote by $B_i'$ the branch of $X_0^\alpha$ corresponding to $B_i$. By the hypothesis and Proposition \ref{Prop_HolderSatWithSameMultiplicity}, we have $m(B_i^\alpha)=m(B_i)$, for every $i$. Moreover, by Theorem \ref{Theo_HolderSatReducibleCurve}, we have $\Ocal_{B_i}\subseteq\Ocal_{B_i'} \subseteq\Ocal_{B_i^\alpha}$. This implies $m(B_i')=m(B_i)$ for every $i$, and hence
$$
m(X_0^\alpha) =\sum_{i=1}^rm(B_i') =\sum_{i=1}^rm(B_i) =m(X_0).
$$
Then, using Lemma \ref{Lem_MultiplicitySandwich}, we get $m(X_0)=m(X_0^\alpha)\leq m(X_y)\leq m(X_0)$ for every $y\in Y$ so the multiplicity is constant. Since $p$ is centered and $\alpha$-equisaturated, then Proposition \ref{Prop_ProductConsequences}(4) gives constancy of the Milnor number. Since constancy of the Milnor number and of the multiplicity implies that $p$ is Whitney equisingular, we get the result. Finally, Proposition~\ref{Prop_HolderSatPreservesMultiplicityReduced} gives $m(X_0^\alpha)=m(X_0)$ whenever
$$
\alpha>\frac{m_0-r_0}{m_0-r_0+2},
$$
so the same argument applies. For $r_0=1$, this bound is $(m_0-1)/(m_0+1)$.
\end{proof}

We can finally summarize the relation between Whitney equisingularity and $\alpha$-equisaturation in the following statement. Notice that we don't need the family to be centered in the third implication.

\begin{cor}\label{Cor_WhitneySmallAlphaContinuation}
Let $p:X\to Y$ be a family of curves. Let $\delta_0 = \delta(X_0)$, $m_0 = m(X_0)$, $r_0 = r(X_0)$ and $\alpha_m = \max_{i}\alpha_m(B_i)$ where $B_i$ are the branches of $X_0$. Then
\begin{enumerate}
    \item If $p$ is centered and Whitney equisingular, then $p$ is $\alpha\text{-equisaturated for every } \alpha \in \big(0,\frac{1}{2\delta_0+1}\big].$
    \item If $p$ is centered and $\alpha$-equisaturated for some $\alpha \in \big(\alpha_m,1\big]$, then it is Whitney equisingular.
    \item If $p$ is $\alpha$-equisaturated for some $\alpha \in \big(\frac{m_0-r_0}{m_0-r_0+2},1\big]$, then it is Whitney equisingular.
\end{enumerate}
\end{cor}

\begin{proof}
    The first implication comes from Theorem \ref{Th_UniformSmallExponent}. The implications $2.$ and $3.$ come from Theorem \ref{Th_EquisaturationWhitney}. The only thing to prove is that we don't need the family to be centered in point $3$. Indeed, we have 
    $$\frac{m_0-r_0}{m_0-r_0+2} \geq 1/3 \iff m_0\geq r_0+1$$
    which is the case since at least one of the branches is singular and the conclusion then comes from Corollary \ref{Cor_NoSplittingOneThird}. If all the branches of $X_0$ are smooth, their simultaneous parametrizations show that the branches of every nearby fiber $X_y$ are also smooth. Then the only singularity appears at the crossing of the branches, which is $\sigma(y)$, so the family is centered.
\end{proof}

\begin{rmq}\label{Prop_NoUniversalWhitneyExponent}
There is no exponent, independent of the family, for which Whitney equisingularity implies $\alpha$-equisaturation. More precisely, for every $\alpha\in(0,1]$, we can construct a centered Whitney equisingular one-parameter family of branches which is not $\alpha$-equisaturated.

Fix $\alpha\in(0,1]$ and choose an even integer $n\geq4$
such that $\alpha(n-3)>4$. Consider the family parametrized by
$$
(t,\lambda)\mapsto \big(t^n, t^{2n-1}+\lambda t^{m+2},t^{2n},t^{2n+1},\ldots,t^{3n-1},\lambda\big),
$$
with $p$ induced by projection onto $\lambda$. One can check that this is a centered family with the invariants $\delta(X_\lambda)=2n-3$, $\mu(X_\lambda)=4n-6$ and $m(X_\lambda)=n$ being constant, and so the family is Whitney equisingular. The characteristic sequence of $X_\lambda$ is $(n,2n-1)$ for $\lambda=0$ and $(n,n+2,2n-1)$ for $\lambda\ne0$. Indeed, $n$ is even, $\gcd(n,n+2)=2$, and $2n-1$ is the smallest odd exponent occurring in either algebra. By Proposition \ref{Prop_ComputeSemigroupOfTheHolderSat}, the saturations are monomial and their value semigroups are
$$
\Gamma(X_0^\alpha) = n\N\cup\big\{q\in\N \mid q\geq\alpha(2n-1)\big\}
$$
and, for $\lambda\neq 0$,
$$
\Gamma(X_\lambda^\alpha) =n\N \cup \big\{q\in2\N \mid q\geq\alpha(n+2)\big\} \cup\big\{q\in\N \mid q\geq\alpha(2n-1)\big\}.
$$
The interval $\big[\alpha(n+2),\,\alpha(2n-1)\big)$ has length $\alpha(n-3)>4$, so it contains two consecutive even integers. At least one of them is not divisible by $n$. Hence $\Gamma(X_0^\alpha)\subsetneq\Gamma(X_\lambda^\alpha)$ and therefore $\delta(X_\lambda^\alpha)<\delta(X_0^\alpha)$ and $\lambda\neq 0$. In particular, $X_\lambda^\alpha\not\simeq X_0^\alpha$.

Finally, suppose that $p$ were $\alpha$-equisaturated, with $X_Y^\alpha\simeq F\times Y$ preserving the base. By Proposition \ref{Prop_ProductConsequences}, we have $F\simeq X_\lambda^\alpha$ for generic $\lambda$ and $\Ocal_F\subseteq\Ocal_{X_0^\alpha}\subseteq\Ocal_{\overline{X_0}}$ otherwise. Consequently,
$$
\delta(X_\lambda^\alpha)=\delta(F)\geq\delta(X_0^\alpha)
$$
for generic $\lambda$, contradicting the strict inequality above.
Thus the family is not $\alpha$-equisaturated.
\end{rmq}

\bigskip
\noindent\begin{minipage}{\textwidth}\small
\textsc{Fran\c{c}ois Bernard}, Institute of Mathematics and Informatics, Bulgarian Academy of Sciences\\
\textit{Email address:} \texttt{f.bernard@math.bas.bg}

\medskip
\noindent\textsc{Antoni Rangachev}, Institute of Mathematics and Informatics, Bulgarian Academy of Sciences, and Institut de Math\'ematiques de Jussieu-Paris Rive Gauche, Paris, France\\
\textit{Email address:} \texttt{rangachev@imj-prg.fr}
\end{minipage}

\end{document}